\documentclass[11pt]{article}
\UseRawInputEncoding
\usepackage{amscd}
\usepackage{amssymb,amsfonts,amsmath,psfrag,amscd,cite}
\usepackage{amsmath}
  \usepackage{paralist}
  \usepackage{graphics}
  \usepackage{epsfig}
  \usepackage[colorlinks=true]{hyperref}
  \usepackage[all]{xy}
\usepackage{graphicx}
\newtheorem{theo}{Theorem}
\newtheorem{prop}[theo]{Proposition}

\newtheorem{defi}[theo]{Definition}
\newtheorem{lemm}[theo]{Lemma}
\newtheorem{rema}[theo]{Remark}

\makeatletter \@addtoreset{equation}{section}
\@addtoreset{theo}{section} \makeatother

\begin{document}
\date{}
\title{Hamilton-Jacobi Equations for Controlled Magnetic Hamiltonian \\
System with Nonholonomic Constraint}
\author{Hong Wang  \\
School of Mathematical Sciences and LPMC,\\
Nankai University, Tianjin 300071, P.R.China\\
E-mail: hongwang@nankai.edu.cn}
\date{\emph{ Dedicated to 110th anniversary of the birth of Professor Shiing Shen Chern}\\
June 22, 2022} \maketitle

{\bf Abstract.} In order to describe the impact of different geometric structures
and constraints for the dynamics of a regular controlled Hamiltonian system,
in this paper, we first define a kind of controlled magnetic
Hamiltonian (CMH) system, and give a
good expression of the dynamical vector field of the CMH system,
such that we can describe the magnetic vanishing condition
and the CMH-equivalence,
and derive precisely the geometric constraint conditions
of the magnetic symplectic form for the
dynamical vector field of the CMH system,
which are called the Type I and Type II of Hamilton-Jacobi equation.
Secondly, we prove that the
CMH-equivalence for the CMH systems leaves the
solutions of corresponding to Hamilton-Jacobi equations invariant,
if the associated magnetic Hamiltonian systems are equivalent.
Thirdly, we consider the CMH system with nonholonomic constraint,
and derive a distributional CMH system, which is
determined by a non-degenerate distributional two-form induced
from the magnetic symplectic form. Then we drive precisely
two types of Hamilton-Jacobi equation for the distributional CMH system.
Moreover, we generalize the above results for the nonholonomic
reducible CMH system with symmetry,
and prove two types of Hamilton-Jacobi theorems
for the nonholonomic reduced distributional CMH system.
These research works reveal the deeply internal
relationships of the magnetic symplectic forms,
the nonholonomic constraints, the dynamical
vector fields and controls of the CMH systems.\\

{\bf Keywords:}\;\; controlled magnetic Hamiltonian system,
\;\; geometric constraint condition, \;\; Hamilton-Jacobi equation,
 \;\;\; CMH-equivalence, \;\;\; nonholonomic constraint. \\

{\bf AMS Classification:} 53D20 \;\; 70H20 \;\; 70Q05.
\tableofcontents

\section{Introduction}

It is well-known that Hamilton-Jacobi theory is an important research subject
in mathematics and analytical mechanics,
see Abraham and Marsden \cite{abma78}, Arnold
\cite{ar89} and Marsden and Ratiu \cite{mara99},
and the Hamilton-Jacobi equation is
also fundamental in the study of the quantum-classical relationship
in quantization, and it also plays an important role
in the study of stochastic dynamical systems, see
Woodhouse \cite{wo92}, Ge and Marsden \cite{gema88},
and L\'{a}zaro-Cam\'{i} and Ortega \cite{laor09}.
For these reasons it is described as a useful tool in the study of
Hamiltonian system theory, and has been extensively developed in
past many years and become one of the most active subjects in the
study of modern applied mathematics and analytical mechanics.\\

Just as we have known that Hamilton-Jacobi theory from the
variational point of view is originally developed by Jacobi in 1866,
which state that the integral of Lagrangian of a mechanical system along the
solution of its Euler-Lagrange equation satisfies the
Hamilton-Jacobi equation. The classical description of this problem
from the generating function and the geometrical point of view is
given by Abraham and Marsden in \cite{abma78} as follows:
Let $Q$ be an $n$-dimensional smooth manifold and $TQ$
the tangent bundle, $T^* Q$ the cotangent bundle with a canonical
symplectic form $\omega$ and the projection $\pi_Q: T^* Q
\rightarrow Q $ induces the map $T\pi_{Q}: TT^* Q \rightarrow TQ. $
\begin{theo}
Assume that the triple $(T^*Q,\omega,H)$ is a Hamiltonian system
with Hamiltonian vector field $X_H$, and $W: Q\rightarrow
\mathbb{R}$ is a given generating function. Then the following two assertions
are equivalent:\\
\noindent $(\mathrm{i})$ For every curve $\sigma: \mathbb{R}
\rightarrow Q $ satisfying $\dot{\sigma}(t)= T\pi_Q
(X_H(\mathbf{d}W(\sigma(t))))$, $\forall t\in \mathbb{R}$, then
$\mathbf{d}W \cdot \sigma $ is an integral curve of the Hamiltonian
vector field $X_H$.\\
\noindent $(\mathrm{ii})$ $W$ satisfies the Hamilton-Jacobi equation
$H(q^i,\frac{\partial W}{\partial q^i})=E, $ where $E$ is a
constant.
\end{theo}

From the proof of the above theorem given in
Abraham and Marsden \cite{abma78}, we know that
the assertion $(\mathrm{i})$ with equivalent to
Hamilton-Jacobi equation $(\mathrm{ii})$ by the generating function,
gives a geometric constraint condition of the canonical symplectic form
on the cotangent bundle $T^*Q$
for Hamiltonian vector field of the system.
Thus, the Hamilton-Jacobi equation reveals the deeply internal relationships of
the generating function, the canonical symplectic form
and the dynamical vector field of a Hamiltonian system.\\

On the other hand, the authors in Marsden et al.\cite{mawazh10}
define a regular controlled Hamiltonian (RCH) system, which is a
Hamiltonian system with external force and control.
In general, an RCH system, under the actions of external force and control, is not
Hamiltonian, however, it is a dynamical system closely related to a
Hamiltonian system, and it can be explored and studied by extending
the methods for external force and control in the study of Hamiltonian systems.
Thus, one can emphasize explicitly the impact of external force
and control in the study for the RCH systems.
However, since an RCH system defined on the cotangent bundle
$T^*Q$,  may not be a Hamiltonian system,
and it may have no generating function,
we cannot give the Hamilton-Jacobi theorem for the RCH system
just like same as the above Theorem 1.1. We have to look for a new way.
It is worthy of noting that, in Wang \cite{wa13d} the author derives precisely
the geometric constraint conditions
of canonical symplectic form for the
dynamical vector field of an RCH system. These conditions
are called the Type I and Type II of Hamilton-Jacobi equation,
which are the development of the Type I and Type II of
Hamilton-Jacobi equation for a Hamiltonian system
given in Wang \cite{wa17}.
Moreover, the author proves that the RCH-equivalence for the RCH systems
leaves the solutions of corresponding to Hamilton-Jacobi equations invariant£¬
if the associated Hamiltonian systems are equivalent..\\

It is well known that the different structures of geometry determine
the different Hamiltonian systems and RCH systems, as well as their dynamics.
In order to describe the impact of different structures of geometry
for the dynamics of RCH system and Hamilton-Jacobi equations, we consider
the magnetic symplectic form $\omega^B= \omega-
\pi^*_Q B,$ where $\omega$ is the usual canonical symplectic
form on $T^*Q$, and $B$ is the closed two-form on $Q$.
A magnetic Hamiltonian system is a Hamiltonian system defined by the
magnetic symplectic form, which is a canonical Hamiltonian system
coupling the action of a magnetic field $B$.
A controlled magnetic Hamiltonian (CMH) system on $T^*Q$ is
a magnetic Hamiltonian system $(T^\ast Q,\omega^B,H)$
with external force $F$ and control $W$, where
$F: T^*Q\rightarrow T^*Q$ is the fiber-preserving map,
and $W\subset T^*Q$ is a fiber submanifold of $T^*Q$.
see Wang \cite{wa15a}.
Thus, it is a natural problem how to derive precisely the geometric constraint
conditions of the magnetic symplectic form for the
dynamical vector field of a CMH system,
and how to describe explicitly the relationship between the CMH-equivalence
and the solutions of corresponding Hamilton-Jacobi equations.
These research are one of our goals in this paper.\\

We know that, in mechanics, it is very often that many systems have constraints.
A nonholonomic CMH system is a CMH system with nonholonomic constraint.
Usually, under the restriction given by nonholonomic constraint,
in general, the dynamical vector field of a nonholonomic CMH system
may not be (magnetic) Hamiltonian. Thus, we can not
describe the Hamilton-Jacobi equations for nonholonomic CMH system
from the viewpoint of generating function
as in the classical Hamiltonian case.
In consequence, when a CMH system with nonholonomic constrain,
how to describe the impact of the constrain
for the dynamics of CMH system and its reduced CMH system,
and how to describe Hamilton-Jacobi equations for
the nonholonomic CMH system and the nonholonomic reducible CMH system.
These research are another of our goals in this paper.\\

A brief of outline of this paper is as follows. In the second
section, we first define a kind of controlled magnetic
Hamiltonian (CMH) system by using magnetic symplectic form,
and then give a good expression of the dynamical vector field of the CMH system,
which is the synthetic of magnetic Hamiltonian vector field
and its changes under the actions of the external force and control law,
such that we can describe the magnetic vanishing condition.
In the third section, we first prove a key lemma,
which is an important tool for the proofs of two types of Hamilton-Jacobi
theorems of the CMH system. Then we
derive precisely the geometric constraint conditions of
the magnetic symplectic form for the dynamical vector field
of a CMH system on the cotangent bundle of a
configuration manifold, that is,
the Type I and Type II of Hamilton-Jacobi equation for a CMH system.
Moreover, in the fourth section, we describe the CMH-equivalence
for the CMH system, and prove that the
CMH-equivalence leaves invariant for the
solutions of corresponding to Hamilton-Jacobi equations,
if the associated magnetic Hamiltonian systems are equivalent.
In the fifth section, we consider the CMH system with nonholonomic constraint,
we define a distributional CMH system by analyzing carefully for
the structure of the nonholonomic dynamical vector field,
and this system is determined by a non-degenerate distributional two-form induced
from the magnetic symplectic form. Moreover, in the sixth section,
we derive precisely two types of Hamilton-Jacobi equation for the
distributional CMH system. In the seventh section,
we generalize the above results for the nonholonomic
reducible RCH system with symmetry,
and prove two types of Hamilton-Jacobi theorems
for the nonholonomic reduced distributional RCH system.
These research works reveal the deeply internal
relationships of the geometrical structures of phase spaces, the dynamical
vector fields and controls of the CMH systems,
and make us have much deeper understanding and
recognition for the structures of Hamiltonian system, RCH
system and CMH system.

\section{Controlled Magnetic Hamiltonian System}

In this section, we first define a kind of controlled magnetic
Hamiltonian (CMH) system by using magnetic symplectic form,
and then give a good expression of the dynamical vector field of the CMH system,
by using the vertical lift map of a vector along a fiber. This expression
is the synthetic of magnetic Hamiltonian vector field
and its changes under the actions of the external force and control law,
such that we can describe the magnetic vanishing condition and CMH-equivalence.
We shall follow some of the notations and conventions introduced in Abraham
and Marsden \cite{abma78}, Arnold \cite{ar89}, Marsden et al.
\cite{mamiorpera07, mawazh10}, Marsden and Ratiu \cite{mara99},
Ortega and Ratiu \cite{orra04},
Wang \cite{wa15a} and Wang \cite{wa17}.
For convenience, in this paper,
we assume that all manifolds are real, smooth and finite dimensional
and all actions are smooth left actions.
and all controls appearing in this paper are the admissible controls.\\

In the reduction theory and application of Hamiltonian systems, the
Marsden-Weinstein reduction for a Hamiltonian system
with symmetry and momentum map is very important and foundational,
see Marsden and Weinstein \cite{mawe74},
and Libermann and Marle \cite{lima87}, Marsden \cite{ma92}.
But, from the classification of symplectic
reduced space of the cotangent bundle $T^* Q$,
see Marsden et al. \cite{mamiorpera07}
and Marsden and Perlmutter \cite{mape00},
we know that the set of Hamiltonian systems with symmetries
and momentum maps on the cotangent bundle $T^*Q$
is not complete under the Marsden-Weinstein reduction,
that is, the symplectic reduced system of a
Hamiltonian system with symmetry and momentum map defined on the cotangent bundle
$T^*Q$ may not be a Hamiltonian system on a cotangent bundle.
In consequence, if we define directly a controlled Hamiltonian system
with symmetry on the cotangent bundle $T^*Q$,
then the symplectic reduced controlled
Hamiltonian system may not have definition.\\

In order to describe uniformly RCH systems defined on a cotangent
bundle and on the regular reduced spaces, in this section we
first define an RCH system on a symplectic fiber bundle,
see Marsden et al.\cite{mawazh10}. Then we can
obtain the RCH system and the CMH system on the cotangent bundle of a configuration
manifold as the special cases, and give a
good expression of the dynamical vector field of the CMH system,
which is the synthetic of magnetic Hamiltonian vector field
$X^B_H$ and its changes under the actions of the external force $F$
and control law $u$,
such that we can describe the magnetic vanishing condition,
and discuss CMH-equivalence in fourth section. In
consequence, we can regard the associated Hamiltonian system
and the magnetic Hamiltonian system on the
cotangent bundle as the special cases of the RCH system
and the CMH system without
external forces and controls, such that we can study the RCH system
and the CMH system by extending
the methods for external force and control
in the study of (magnetic) Hamiltonian system.\\

Let $(E,M,\pi )$ be a fiber bundle, and for each point $x \in M$,
assume that the fiber $E_x=\pi^{-1}(x)$ is a smooth submanifold of $E$
and with a symplectic form $\omega_{E}(x)$, that is,
$(E, \omega_E)$ is a symplectic fiber bundle. If for any function
$H: E \rightarrow \mathbb{R}$, we have a Hamiltonian vector field $X_H$,
which satisfies the Hamilton's equation, that is,
$\mathbf{i}_{X_H}\omega_E=\mathbf{d}H$, then $(E, \omega_E, H )$ is a
Hamiltonian system. Moreover, if considering the external force and
control, we can define a kind of regular controlled Hamiltonian
(RCH) system on the symplectic fiber bundle $E$ as
follows.
\begin{defi}
(RCH System) A regular controlled Hamiltonian (RCH) system on $E$ is a 5-tuple
$(E, \omega_E, H, F, W)$, where $(E, \omega_E, H )$ is a
Hamiltonian system, and the function $H: E \rightarrow \mathbb{R}$
is called the Hamiltonian, a fiber-preserving map $F: E\rightarrow
E$ is called the (external) force map, and a fiber submanifold $W$
of $E$ is called the control subset.
\end{defi}
Sometimes, $W$ is also denoted the set of fiber-preserving maps from
$E$ to $W$. When a feedback control law $u: E\rightarrow W$ is
chosen, the 5-tuple $(E, \omega_E, H, F, u)$ is a regular closed-loop
dynamic system. In particular, when $Q$ is an $n$-dimensional smooth manifold, and
$T^\ast Q$ its cotangent bundle with a symplectic form $\omega$ (not
necessarily canonical symplectic form), then $(T^\ast Q, \omega )$
is a symplectic vector bundle. If we take that $E= T^* Q$, from
above definition we can obtain an RCH system on the cotangent bundle
$T^\ast Q$, that is, 5-tuple $(T^\ast Q, \omega, H, F, W)$. Where
the fiber-preserving map $F: T^*Q\rightarrow T^*Q$ is the (external)
force map, which is the reason that the fiber-preserving map $F:
E\rightarrow E$ is called an (external) force map in above
definition.\\

In order to describe the impact of different structures of geometry
for the RCH systems, we shall consider the
magnetic symplectic form on $T^*Q$ as follows:
Assume that $T^*Q$ with the canonical symplectic form $\omega$,
and $B$ is a closed two-form on $Q$,
then $\omega^B= \omega- \pi_Q^*B$ is a symplectic form on $T^*Q$,
where $\pi_Q^*: T^*Q \rightarrow T^*T^*Q $. The
$\omega^B$ is called a magnetic symplectic form, and $\pi_Q^*B$ is called
a magnetic term on $T^*Q$, see Marsden et al. \cite{mamiorpera07}.\\

A magnetic Hamiltonian system is a 3-tuple $(T^\ast Q,\omega^B,H)$,
which is a Hamiltonian system defined by the
magnetic symplectic form $\omega^B$, that is,
a canonical Hamiltonian system
coupling the action of a magnetic field $B$. For a given Hamiltonian $H$,
the dynamical vector field $X^B_H$, which is called
the magnetic Hamiltonian vector field,
satisfies the magnetic Hamilton's equation, that is,
$\mathbf{i}_{X^B_{H} }\omega^B= \mathbf{d}H $.
In canonical cotangent bundle coordinates, for any $q \in Q
, \; (q,p)\in T^* Q, $ we have that
$$
\omega=\sum^n_{i=1} \mathbf{d}q^i \wedge \mathbf{d}p_i ,
\;\;\;\;\;\; B=\sum^n_{i,j=1}B_{ij}\mathbf{d}q^i \wedge
\mathbf{d}q^j ,\;\;\; \mathbf{d}B=0, $$
$$\omega^B= \omega -\pi_Q^*B=\sum^n_{i=1} \mathbf{d}q^i \wedge
\mathbf{d}p_i- \sum^n_{i,j=1}B_{ij}\mathbf{d}q^i \wedge
\mathbf{d}q^j,
$$
and the magnetic Hamiltonian vector field $X^B_H$ with respect to
the magnetic symplectic form $\omega^B$ can be expressed that
$$
X^B_H= \sum^n_{i=1} (\frac{\partial H}{\partial
p_i}\frac{\partial}{\partial q^i} - \frac{\partial H}{\partial
q^i}\frac{\partial}{\partial p_i})-
\sum^n_{i,j=1}B_{ij}\frac{\partial H}{\partial
p_j}\frac{\partial}{\partial p_i}.
$$
See Marsden et al. \cite{mamiorpera07}.\\

Moreover, if considering the external force and
control, we can define a kind of controlled magnetic Hamiltonian
(CMH) system on $T^*Q$ as follows.
\begin{defi}
(CMH System) A controlled magnetic Hamiltonian (CMH) system
on $T^*Q$ is a 5-tuple $(T^*Q, \omega^B, H, F, W)$, which is a
magnetic Hamiltonian system $(T^\ast Q,\omega^B,H)$
with external force $F$ and control $W$, where $F: T^*Q\rightarrow T^*Q$ is
the fiber-preserving map, and $W\subset T^*Q$ is a fiber submanifold,
which is called the control subset.
\end{defi}

From the above Definition 2.1 and Definition 2.2 we know that
a CMH system on $T^*Q$ is also an RCH system on $T^*Q$, but its symplectic structure
is given by a magnetic symplectic form, and the set of
the CMH systems on $T^*Q$ is a subset of the set of the RCH systems on $T^*Q$.
When a feedback control law $u: T^*Q\rightarrow W$ is
chosen, the 5-tuple $(T^*Q, \omega^B, H, F, u)$ is a regular closed-loop
dynamic system.\\

In order to describe the dynamics of the CMH system
$(T^*Q,\omega^B,H,F,W)$ with a control law $u$, we need to give a good
expression of the dynamical vector field of the CMH system. At first, we
introduce a notations of vertical lift maps of a vector along a
fiber, also see Marsden et al. \cite{mawazh10}.
For a smooth manifold $E$, its tangent bundle $TE$ is a
vector bundle, and for the fiber bundle $\pi: E \rightarrow M$, we
consider the tangent mapping $T\pi: TE \rightarrow TM$ and its
kernel $ker (T\pi)=\{\rho\in TE| T\pi(\rho)=0\}$, which is a vector
subbundle of $TE$. Denote by $VE:= ker(T\pi)$, which is called a
vertical bundle of $E$. Assume that there is a metric on $E$, and we
take a Levi-Civita connection $\mathcal{A}$ on $TE$, and denote by
$HE:= ker(\mathcal{A})$, which is called a horizontal bundle of $E$,
such that $TE= HE \oplus VE. $ For any $x\in M, \; a_x, b_x \in E_x,
$ any tangent vector $\rho(b_x)\in T_{b_x}E$ can be split into
horizontal and vertical parts, that is, $\rho(b_x)=
\rho^h(b_x)\oplus \rho^v(b_x)$, where $\rho^h(b_x)\in H_{b_x}E$ and
$\rho^v(b_x)\in V_{b_x}E$. Let $\gamma$ be a geodesic in $E_x$
connecting $a_x$ and $b_x$, and denote by $\rho^v_\gamma(a_x)$ a
tangent vector at $a_x$, which is a parallel displacement of the
vertical vector $\rho^v(b_x)$ along the geodesic $\gamma$ from $b_x$
to $a_x$. Since the angle between two vectors is invariant under a
parallel displacement along a geodesic, then
$T\pi(\rho^v_\gamma(a_x))=0, $ and hence $\rho^v_\gamma(a_x) \in
V_{a_x}E. $ Now, for $a_x, b_x \in E_x $ and tangent vector
$\rho(b_x)\in T_{b_x}E$, we can define the vertical lift map of a
vector along a fiber given by
$$\textnormal{vlift}: TE_x \times E_x \rightarrow TE_x; \;\;
\textnormal{vlift}(\rho(b_x),a_x) = \rho^v_\gamma(a_x). $$
It is easy to check from the basic fact in differential geometry
that this map does not depend on the choice of $\gamma$. If $F: E
\rightarrow E$ is a fiber-preserving map, for any $x\in M$, we have
that $F_x: E_x \rightarrow E_x$ and $TF_x: TE_x \rightarrow TE_x$,
then for any $a_x \in E_x$ and $\rho\in TE_x$, the vertical lift of
$\rho$ under the action of $F$ along a fiber is defined by
$$(\textnormal{vlift}(F_x)\rho)(a_x)
=\textnormal{vlift}((TF_x\rho)(F_x(a_x)), a_x)
= (TF_x\rho)^v_\gamma(a_x), $$
where $\gamma$ is a geodesic in $E_x$ connecting $F_x(a_x)$ and
$a_x$.\\

In particular, when $\pi: E \rightarrow M$ is a vector bundle, for
any $x\in M$, the fiber $E_x=\pi^{-1}(x)$ is a vector space. In this
case, we can choose the geodesic $\gamma$ to be a straight line, and
the vertical vector is invariant under a parallel displacement along
a straight line, that is, $\rho^v_\gamma(a_x)= \rho^v(b_x).$
Moreover, when $E= T^*Q, \; M=Q $, by using the local trivialization
of $TT^*Q$, we have that $TT^*Q\cong TQ \times T^*Q$. Because of
$\pi: T^*Q \rightarrow Q$, and $T\pi: TT^*Q \rightarrow TQ$, then in
this case, for any $\alpha_x, \; \beta_x \in T^*_x Q, \; x\in Q, $
we know that $(0, \beta_x) \in V_{\beta_x}T^*_x Q, $ and hence we
can get that
$$ \textnormal{vlift}((0, \beta_x)(\beta_x), \alpha_x) = (0, \beta_x)(\alpha_x)
= \left.\frac{\mathrm{d}}{\mathrm{d}s}\right|_{s=0}(\alpha_x+s\beta_x),
$$ which is consistent with the definition of vertical lift map
along fiber in Marsden and Ratiu \cite{mara99}.\\

For a given CMH system $(T^\ast Q, \omega^B, H, F, W)$, the dynamical
vector field of the associated magnetic Hamiltonian system $(T^\ast Q,
\omega^B, H) $ is  $X^B_H$, which satisfies the equation
$\mathbf{i}_{X^B_H}\omega^B=\mathbf{d}H$. If considering the
external force $F: T^*Q \rightarrow T^*Q, $ by using the above
notation of vertical lift map of a vector along a fiber, the change
of $X^B_H$ under the action of $F$ is that
$$\textnormal{vlift}(F)X^B_H(\alpha_x)
= \textnormal{vlift}((TFX^B_H)(F(\alpha_x)), \alpha_x)
= (TFX^B_H)^v_\gamma(\alpha_x),$$
where $\alpha_x \in T^*_x Q, \; x\in Q $ and $\gamma$ is a straight
line in $T^*_x Q$ connecting $F_x(\alpha_x)$ and $\alpha_x$. In the
same way, when a feedback control law $u: T^\ast Q \rightarrow W$ is
chosen, the change of $X^B_H$ under the action of $u$ is that
$$\textnormal{vlift}(u)X^B_H(\alpha_x)
= \textnormal{vlift}((TuX^B_H)(u(\alpha_x)), \alpha_x)
= (TuX^B_H)^v_\gamma(\alpha_x).$$
In consequence, we can give an expression of the dynamical vector
field of the CMH system as follows.
\begin{theo}
The dynamical vector field of a CMH system $(T^\ast Q,\omega^B,H,F,W)$
with a control law $u$ is the synthetic of magnetic Hamiltonian vector field
$X^B_H$ and its changes under the actions of the external force $F$
and control law $u$, that is,
\begin{equation}
X_{(T^\ast Q,\omega^B,H,F,u)}(\alpha_x)
= X^B_H(\alpha_x)+ \textnormal{vlift}(F)X^B_H(\alpha_x)
+ \textnormal{vlift}(u)X^B_H(\alpha_x), \;\; \label{2.1}
\end{equation}
 for any $\alpha_x \in T^*_x
Q, \; x\in Q $. For convenience, it is simply written as
\begin{equation}X_{(T^\ast Q,\omega^B,H,F,u)}
=X^B_H +\textnormal{vlift}(F)^B +\textnormal{vlift}(u)^B. \;\; \label{2.2}
\end{equation}
\end{theo}
Where $\textnormal{vlift}(F)^B=\textnormal{vlift}(F)X^B_H$,
and $\textnormal{vlift}(u)^B=\textnormal{vlift}(u)X^B_H.$ are the
changes of $X^B_H$ under the actions of $F$ and $u$.
We also denote that $\textnormal{vlift}(W)^B= \bigcup\{\textnormal{vlift}(u)X^B_H |
\; u\in W\}$. It is worthy of noting that, in order to deduce and calculate
easily, we always use the simple expression of dynamical vector
field $X_{(T^\ast Q,\omega^B,H,F,u)}$. \\

From the expression (2.2) of the dynamical vector
field of a CMH system, we know that under the actions of the external force $F$
and control law $u$, in general, the dynamical vector
field may not be magnetic Hamiltonian, and hence the CMH system may not
be yet a magnetic Hamiltonian system. However,
it is a dynamical system closed relative to a
magnetic Hamiltonian system, and it can be explored and studied by extending
the methods for external force and control
in the study of magnetic Hamiltonian system.\\

For the magnetic Hamiltonian system $(T^\ast Q, \omega^B, H)$, its
magnetic Hamiltonian vector field $X^B_H$ satisfies the equation
$\mathbf{i}_{X^B_H}\omega^B=\mathbf{d}H$, and for the
associated canonical Hamiltonian system $(T^\ast Q,
\omega, H) $, its canonical Hamiltonian vector field $X_H$
satisfies the equation $\mathbf{i}_{X_H}\omega=\mathbf{d}H$.
Denote by the vector field $X^0= X^B_H-X_H, $ and from the
magnetic symplectic form $\omega^B= \omega- \pi_Q^*B$, we have that
$$
\mathbf{i}_{X^0}\omega=\mathbf{i}_{(X^B_H-X_H)}\omega
=\mathbf{i}_{X^B_H}\omega-\mathbf{i}_{X_H}\omega
=\mathbf{i}_{X^B_H}(\omega^B+ \pi_Q^*B)-\mathbf{i}_{X_H}\omega
=\mathbf{i}_{X^B_H}( \pi_Q^*B).
$$
Thus, $X^0$ is called the magnetic vector field and
$\mathbf{i}_{X^0}\omega=\mathbf{i}_{X^B_H}( \pi_Q^*B)$
is called the magnetic equation, which is determined by
the magnetic term $\pi_Q^*B$ on $T^*Q$. When $B=0$,
then $X^0=0$, the magnetic equation holds trivially.
For the CMH system $(T^\ast Q, \omega^B, H, F, W)$,
from the expression (2.2) of its dynamical
vector field, we have that
\begin{equation}X_{(T^\ast Q,\omega^B,H,F,u)}
=X_H + X^0+\textnormal{vlift}(F)^B +\textnormal{vlift}(u)^B. \;\; \label{2.3}
\end{equation}
If we choose the external force $F$ and control law $u$, such that
\begin{equation}
 X^0+\textnormal{vlift}(F)^B +\textnormal{vlift}(u)^B=0, \;\; \label{2.4}
\end{equation}
then from (2.3) we have that $X_{(T^\ast Q,\omega^B,H,F,u)}
=X_H, $ that is, in this case the dynamical vector
field of the CMH system is just the canonical Hamiltonian vector field,
and the motion of the CMH system is just same like the motion of canonical
Hamiltonian system without the actions of magnetic, external force and control.
Thus, the condition (2.4) is called the magnetic vanishing condition for
the CMH system $(T^\ast Q, \omega^B, H, F, W)$.\\

To sum up the above discussion, we have the following theorem.
\begin{theo}
If the external force $F$ and the control law $u$ for a CMH system
$(T^\ast Q, \omega^B, H, F, u)$ satisfy the magnetic vanishing condition
(2.4), then its dynamical vector field $X_{(T^\ast Q,\omega^B,H,F,u)}$
is just the canonical Hamiltonian vector field $X_H$ for the
associated canonical Hamiltonian system $(T^\ast Q, \omega, H) $.
\end{theo}

\section{Two Types of Hamilton-Jacobi Equation for a CMH System }

In this section, we shall derive precisely the geometric constraint conditions of
the magnetic symplectic form for the dynamical vector field
of a CMH system, that is, Type I and Type II of Hamilton-Jacobi equation for
the CMH system. In order to do this, in the following we first give
an important notion and prove a key lemma, which is an important
tool for the proofs of two types of
Hamilton-Jacobi theorem for the CMH system.\\

Denote by $\Omega^i(Q)$ the set of all i-forms on $Q$, $i=1,2.$
For any $\gamma \in \Omega^1(Q),\; q\in Q, $ then $\gamma(q)\in T_q^*Q, $
and we can define a map $\gamma: Q \rightarrow T^*Q, \; q \rightarrow (q, \gamma(q)).$
Hence we say often that the map $\gamma: Q
\rightarrow T^*Q$ is an one-form on $Q$. If the one-form $\gamma$ is closed,
then $\mathbf{d}\gamma(x,y)=0, \; \forall\; x, y \in TQ$.
Note that for any $v, w \in TT^* Q, $ we have that
$\mathbf{d}\gamma(T\pi_{Q}(v),T\pi_{Q}(w))=\pi^*(\mathbf{d}\gamma )(v, w)$
is a two-form on the cotangent bundle $T^*Q$, where
$\pi^*: T^*Q \rightarrow T^*T^*Q.$ Thus,
in the following we can give a weaker notion.
\begin{defi}
The one-form $\gamma$ is called to be closed with respect to $T\pi_{Q}:
TT^* Q \rightarrow TQ, $ if for any $v, w \in TT^* Q, $ we have that
$\mathbf{d}\gamma(T\pi_{Q}(v),T\pi_{Q}(w))=0. $
\end{defi}

From the above definition we know that, if $\gamma$ is a closed one-form,
then it must be closed with respect to $T\pi_{Q}: TT^* Q \rightarrow
TQ. $ Conversely, if $\gamma$ is closed with respect to
$T\pi_{Q}: TT^* Q \rightarrow TQ, $ then it may not be closed. We can
prove a general result as follows, its proof given in Wang \cite{wa17},
which states that the notion that $\gamma$ is closed
with respect to $T\pi_{Q}: TT^* Q \rightarrow TQ, $
is not equivalent to the notion that $\gamma$ is closed.

\begin{prop}
Assume that $\gamma: Q \rightarrow T^*Q$ is an one-form on $Q$ and
it is not closed. we define the set $N$, which is a subset of $TQ$,
such that the one-form $\gamma$ on $N$ satisfies the condition that
for any $x,y \in N, \; \mathbf{d}\gamma(x,y)\neq 0. $ Denote by
$Ker(T\pi_Q)= \{u \in TT^*Q| \; T\pi_Q(u)=0 \}, $ and $T\gamma: TQ
\rightarrow TT^* Q $ is the tangent map of $\gamma: Q \rightarrow T^*Q. $
If $T\gamma(N)\subset Ker(T\pi_Q), $ then
$\gamma$ is closed with respect to $T\pi_{Q}: TT^* Q \rightarrow TQ.
$\end{prop}

For the one-form $\gamma: Q \rightarrow T^*Q$, $\mathbf{d}\gamma$
is a two-form on $Q$. Assume that $B$ is a closed two-form on $Q$,
we say that the $\gamma$ satisfies condition $\mathbf{d}\gamma=-B$,
if for any $ x, y \in TQ$, we have that $(\mathbf{d}\gamma +B)(x,y)=0.$
In the following we can give a new notion.
\begin{defi}
Assume that $\gamma: Q
\rightarrow T^*Q$ is an one-form on $Q$,
we say that the $\gamma$ satisfies condition
$\mathbf{d}\gamma=-B$ with respect to $T\pi_{Q}:
TT^* Q \rightarrow TQ, $ if for any $v, w \in TT^* Q, $ we have that
$(\mathbf{d}\gamma +B)(T\pi_{Q}(v),T\pi_{Q}(w))=0. $
\end{defi}

From the above Definition 3.1 and Definition 3.3, we know that,
when $B=0$, the condition that, $\gamma$ satisfies condition
$\mathbf{d}\gamma=-B$ with respect to $T\pi_{Q}:
TT^* Q \rightarrow TQ, $ become that $\gamma$ is closed with respect to $T\pi_{Q}:
TT^* Q \rightarrow TQ. $  Now, we prove the following lemma,
which  is a generalization of a corresponding to lemma given by
Wang \cite{wa17},
and the lemma is a very important tool for our research.

\begin{lemm}
Assume that $\gamma: Q \rightarrow T^*Q$ is an one-form on $Q$, and
$\lambda=\gamma \cdot \pi_{Q}: T^* Q \rightarrow T^* Q .$
For the magnetic symplectic form $\omega^B= \omega- \pi_Q^*B $ on $T^*Q$,
where $\omega$ is the canonical symplectic form on $T^*Q$,
then we have that the following two assertions hold.\\
\noindent $(\mathrm{i})$ For any $v, w \in
TT^* Q, \; \lambda^*\omega^B(v,w)= -(\mathbf{d}\gamma+B)(T\pi_{Q}(v), \;
T\pi_{Q}(w))$; \\
\noindent $(\mathrm{ii})$ For any $v, w \in TT^* Q, \;
\omega^B(T\lambda \cdot v,w)= \omega^B(v, w-T\lambda \cdot
w)-(\mathbf{d}\gamma+B)(T\pi_{Q}(v), \; T\pi_{Q}(w)). $
\end{lemm}

\noindent{\bf Proof:} We first prove the assertion $(\mathrm{i})$.
Since $\omega$ is the canonical symplectic form on $T^*Q$,
we know that there is an unique canonical one-form $\theta$, such that
$\omega= -\mathbf{d} \theta. $ From the Proposition 3.2.11 in
Abraham and Marsden \cite{abma78}, we have that for the one-form
$\gamma: Q \rightarrow T^*Q, \; \gamma^* \theta= \gamma. $ Then we
can obtain that for any $x, y \in TQ,$
\begin{align*}
\gamma^*\omega(x,y) = \gamma^* (-\mathbf{d} \theta) (x, y) =
-\mathbf{d}(\gamma^* \theta)(x, y)= -\mathbf{d}\gamma (x, y).
\end{align*}
Note that $\lambda=\gamma \cdot \pi_{Q}: T^* Q \rightarrow T^* Q, $
and $\lambda^*= \pi_{Q}^* \cdot \gamma^*: T^*T^* Q \rightarrow
T^*T^* Q, $ then we have that  for any $v, w \in TT^* Q $,
\begin{align*}
\lambda^*\omega(v,w) &= \lambda^* (-\mathbf{d} \theta) (v, w)
=-\mathbf{d}(\lambda^* \theta)(v, w)= -\mathbf{d}(\pi_{Q}^* \cdot
\gamma^* \theta)(v, w)\\ &= -\mathbf{d}(\pi_{Q}^* \cdot\gamma )(v,
w)= -\mathbf{d}\gamma(T\pi_{Q}(v), \; T\pi_{Q}(w)).
\end{align*}
Hence, we have that
\begin{align*}
\lambda^*\omega^B(v,w)& =\lambda^*\omega(v,w)-\lambda^*\cdot \pi_Q^*B(v,w)\\
& =-\mathbf{d}\gamma(T\pi_{Q}(v), \; T\pi_{Q}(w))
-(\pi_Q\cdot \gamma \cdot \pi_{Q})^*B(v,w)\\
& =-\mathbf{d}\gamma(T\pi_{Q}(v), \; T\pi_{Q}(w))- \pi_Q^*B(v,w)\\
& =-(\mathbf{d}\gamma+B)(T\pi_{Q}(v), \; T\pi_{Q}(w)),
\end{align*}
where we have used the relation $\pi_Q\cdot \gamma\cdot \pi_Q= \pi_Q. $
It follows that the assertion $(\mathrm{i})$ holds.\\

Next, we prove the assertion $(\mathrm{ii})$. For any $v, w \in TT^*
Q,$ note that $v- T(\gamma \cdot \pi_Q)\cdot v$ is vertical, because
$$
T\pi_Q(v- T(\gamma \cdot \pi_Q)\cdot v)=T\pi_Q(v)-T(\pi_Q\cdot
\gamma\cdot \pi_Q)\cdot v= T\pi_Q(v)-T\pi_Q(v)=0,
$$
Thus, $\omega(v- T(\gamma \cdot \pi_Q)\cdot v,w- T(\gamma \cdot
\pi_Q)\cdot w)= 0, $ and hence,
$$\omega(T(\gamma \cdot \pi_Q)\cdot v, \; w)=
\omega(v, \; w-T(\gamma \cdot \pi_Q)\cdot w)+ \omega(T(\gamma \cdot
\pi_Q)\cdot v, \; T(\gamma \cdot \pi_Q)\cdot w). $$ However, the
second term on the right-hand side is given by
$$
\omega(T(\gamma \cdot \pi_Q)\cdot v, \; T(\gamma \cdot \pi_Q)\cdot
w)= \gamma^*\omega(T\pi_Q(v), \; T\pi_Q(w))=
-\mathbf{d}\gamma(T\pi_{Q}(v), \; T\pi_{Q}(w)),
$$
It follows that
\begin{align*}
\omega(T\lambda \cdot v,w) &=\omega(T(\gamma \cdot \pi_Q)\cdot v, \;
w)\\ &= \omega(v, \; w-T(\gamma \cdot \pi_Q)\cdot w)-\mathbf{d}\gamma(T\pi_{Q}(v), \; T\pi_{Q}(w))
\\ &= \omega(v,
w-T\lambda \cdot w)-\mathbf{d}\gamma(T\pi_{Q}(v), \; T\pi_{Q}(w)).
\end{align*}
Hence,  we have that
\begin{align*}
& \omega^B(T\lambda \cdot v,w)= \omega(T\lambda \cdot v,w)-\pi_Q^*B(T\lambda \cdot v,w)\\
& =\omega(v, w-T\lambda \cdot w)-\mathbf{d}\gamma(T\pi_{Q}(v), \; T\pi_{Q}(w))
-B(T\pi_Q\cdot T\lambda \cdot v, \; T\pi_{Q}(w))\\
& =\omega^B(v, w-T\lambda \cdot w)+\pi_Q^*B(v, w-T\lambda \cdot w)\\
& \;\;\;\;\; -\mathbf{d}\gamma(T\pi_{Q}(v), \; T\pi_{Q}(w))
-B(T(\pi_Q\cdot \lambda) \cdot v, \; T\pi_{Q}(w))\\
& =\omega^B(v, w-T\lambda \cdot w)+\pi_Q^*B(v, w)-B(T\pi_{Q}(v), \; T\pi_Q\cdot T\lambda \cdot w)\\
& \;\;\;\;\; -\mathbf{d}\gamma(T\pi_{Q}(v), \; T\pi_{Q}(w))
-B(T(\pi_Q\cdot \gamma \cdot \pi_{Q}) \cdot v, \; T\pi_{Q}(w))\\
& =\omega^B(v, w-T\lambda \cdot w)+\pi_Q^*B(v, w)-B(T\pi_{Q}(v), \; T(\pi_Q\cdot \lambda) \cdot w)\\
& \;\;\;\;\; -\mathbf{d}\gamma(T\pi_{Q}(v), \; T\pi_{Q}(w))
-B(T\pi_Q (v), \; T\pi_{Q}(w))\\
& =\omega^B(v, w-T\lambda \cdot w)+\pi_Q^*B(v, w)-B(T\pi_{Q}(v), \; T\pi_{Q}(w))
-(\mathbf{d}\gamma+B)(T\pi_{Q}(v), \; T\pi_{Q}(w))\\
& =\omega^B(v, w-T\lambda \cdot
w)-(\mathbf{d}\gamma+B)(T\pi_{Q}(v), \; T\pi_{Q}(w)).
\end{align*}
Thus, the assertion $(\mathrm{ii})$ holds.
\hskip 0.3cm $\blacksquare$\\

For a given CMH system $(T^*Q,\omega^B,H,F,W)$ on $T^*Q$, by using
the above Lemma 3.4, we can derive precisely the geometric constraint
conditions of the magnetic symplectic form $\omega^B$ for the dynamical
vector field $X_{(T^\ast Q,\omega^B,H,F,u)}$ of the CMH system with a control law $u$,
that is, Type I and Type II of
Hamilton-Jacobi equation for the CMH system.

\begin{theo}
(Type I of Hamilton-Jacobi Theorem for a CMH System)
For the CMH system $(T^*Q,\omega^B,H,F,W)$ with the
magnetic symplectic form $\omega^B= \omega- \pi_Q^*B $ on $T^*Q$,
where $\omega$ is the canonical symplectic form on $T^* Q$
and $B$ is a closed two-form on $Q$,
assume that $\gamma: Q
\rightarrow T^*Q$ is an one-form on $Q$, and
$\tilde{X}^\gamma = T\pi_{Q}\cdot \tilde{X} \cdot \gamma$,
where $\tilde{X}=X_{(T^\ast Q,\omega^B,H,F,u)}$ is the dynamical vector field
of the CMH system $(T^*Q,\omega^B,H,F,W)$ with a control law $u$.
If the one-form $\gamma: Q \rightarrow T^*Q $ satisfies the condition
$\mathbf{d}\gamma=-B $ with respect to $T\pi_{Q}:
TT^* Q \rightarrow TQ, $ then $\gamma$ is a solution of the equation
$T\gamma\cdot \tilde{X}^\gamma= X^B_H\cdot \gamma ,$ where $X^B_H$
is the magnetic Hamiltonian vector field
of the associated magnetic Hamiltonian system $(T^*Q,\omega^B,H),$
and the equation is called the Type I of
Hamilton-Jacobi equation for the CMH system
$(T^*Q,\omega^B,H,F,W)$ with a control law $u$.
Here the maps involved in the theorem are shown
in the following Diagram-1.
\begin{center}
\hskip 0cm \xymatrix{ & T^* Q \ar[d]_{X^B_H}\ar[r]^{\pi_Q}
 & Q \ar[d]_{\tilde{X}^\gamma} \ar[r]^{\gamma} & T^*Q \ar[d]^{\tilde{X}} \\
 & T(T^*Q) & TQ \ar[l]_{T\gamma} & T(T^* Q)\ar[l]_{T\pi_Q}}
\end{center}
$$\mbox{Diagram-1}$$
\end{theo}
\noindent{\bf Proof: }
Since $\tilde{X}=\tilde{X}_{(T^\ast Q,\omega^B,H,F,u)}=X^B_H
+\textnormal{vlift}(F)^B+\textnormal{vlift}(u)^B, $ and
$T\pi_{Q}\cdot \textnormal{vlift}(F)^B=T\pi_{Q}\cdot \textnormal{vlift}(u)^B=0, $
then we have that $T\pi_{Q}\cdot \tilde{X}\cdot \gamma=T\pi_{Q}\cdot X^B_H\cdot \gamma. $
If we take that $v= X^B_H\cdot \gamma \in TT^* Q, $ and for
any $w \in TT^* Q, \; T\pi_{Q}(w)\neq 0, $ from Lemma 3.4(ii) and
$\mathbf{d}\gamma=-B $ with respect to $T\pi_{Q}:
TT^* Q \rightarrow TQ, $ that is,
$(\mathbf{d}\gamma+B)(T\pi_{Q}\cdot X^B_{H}\cdot\gamma, \; T\pi_{Q}\cdot w)=0,$
we have that
\begin{align*}
\omega^B(T\gamma \cdot \tilde{X}^\gamma, \; w)&
=\omega^B(T\gamma \cdot T\pi_{Q} \cdot \tilde{X}\cdot\gamma, \; w)\\
&=\omega^B(T\gamma \cdot T\pi_{Q} \cdot X^B_H\cdot\gamma, \; w)
= \omega^B(T(\gamma \cdot \pi_Q)\cdot X^B_H\cdot \gamma, \; w)\\
&= \omega^B(X^B_H\cdot \gamma, \; w-T(\gamma \cdot \pi_Q)\cdot w)
-(\mathbf{d}\gamma+B)(T\pi_{Q}\cdot X^B_{H}\cdot\gamma, \; T\pi_{Q}\cdot w)\\
& = \omega^B(X^B_H\cdot \gamma, \; w) - \omega^B(X^B_H\cdot \gamma, \;
T\lambda \cdot w).
\end{align*}
Hence, we have that
\begin{equation}
\omega^B(T\gamma \cdot \tilde{X}^\gamma, \; w)- \omega^B(X^B_H\cdot \gamma, \; w)
= -\omega^B(X^B_H\cdot \gamma, \; T\lambda \cdot w). \; \label{3.1}
\end{equation}
If $\gamma$ satisfies the equation $T\gamma\cdot \tilde{X}^\gamma= X^B_H\cdot \gamma ,$
from Lemma 3.4(i) we can obtain that
\begin{align*}
\omega^B(X^B_H\cdot \gamma, \; T\lambda \cdot w) &
= \omega^B(T\gamma \cdot \tilde{X}^\gamma, \; T\lambda \cdot w)\\
&= \omega^B(T\gamma \cdot T\pi_{Q} \cdot \tilde{X}\cdot\gamma, \; T\lambda \cdot w)\\
&= \omega^B(T\gamma \cdot T\pi_{Q} \cdot X^B_H\cdot\gamma, \; T\lambda \cdot w)\\
&= \omega^B(T\lambda \cdot X^B_{H}\cdot\gamma, \; T\lambda \cdot w)\\
&= \lambda^*\omega^B(X^B_{H}\cdot\gamma, \; w)\\
&=-(\mathbf{d}\gamma+B)(T\pi_{Q}\cdot X^B_{H}\cdot\gamma, \; T\pi_{Q}\cdot w)=0,
\end{align*}
since $\gamma: Q \rightarrow T^*Q $ satisfies the condition
$\mathbf{d}\gamma=-B $ with respect to $T\pi_{Q}:
TT^* Q \rightarrow TQ. $
But, because the magnetic symplectic form $\omega^B$ is non-degenerate,
the left side of (3.1) equals zero, only when
$\gamma$ satisfies the equation $T\gamma\cdot \tilde{X}^\gamma= X^B_H\cdot \gamma .$ Thus,
if the one-form $\gamma: Q \rightarrow T^*Q $ satisfies the condition
$\mathbf{d}\gamma=-B$ with respect to $T\pi_{Q}:
TT^* Q \rightarrow TQ, $ then $\gamma$ must be a solution of
the Type I of Hamilton-Jacobi equation
$T\gamma\cdot \tilde{X}^\gamma= X^B_H\cdot \gamma ,$ for
the CMH system $(T^*Q,\omega^B,H,F,W)$ with a control law $u$.
\hskip 0.3cm $\blacksquare$\\

When $B=0$, in this case the magnetic symplectic form $\omega^B$
is just the canonical symplectic form $\omega$ on $T^*Q$,
and the CMH system $(T^*Q,\omega^B,H,F,W)$
is just a canonical RCH system $(T^*Q,\omega,H,F,W)$
and the condition that, the one-form $\gamma: Q \rightarrow T^*Q $ satisfies
the condition $\mathbf{d}\gamma=-B $ with respect to $T\pi_{Q}:
TT^* Q \rightarrow TQ, $ becomes that  $\gamma $ is closed with respect to
$T\pi_Q: TT^* Q \rightarrow TQ.$ Thus, from above Theorem 3.5,
we can obtain Theorem 2.6 in Wang \cite{wa13d}. On the other hand,
it is a natural problem what and how we could do,
if an one-form $\gamma: Q \rightarrow T^*Q $ is not
closed on $Q$ with respect to $T\pi_Q: TT^* Q \rightarrow TQ $
in Theorem 2.6 in Wang \cite{wa13d}, and hence
$\gamma$ is not a solution of the Type I of Hamilton-Jacobi
equation for the canonical RCH system. In this case,
our idea is that we hope to look for a new RCH system,
such that $\gamma$ is a solution of the Type I
of Hamilton-Jacobi equation for the new RCH system.
Note that, if
$\gamma: Q \rightarrow T^*Q $ is not closed on $Q$ with respect to
$T\pi_Q: TT^* Q \rightarrow TQ, $ that is, there exist $v,w \in TT^*Q,$
such that $\mathbf{d}\gamma(T\pi_Q(v), \; T\pi_Q(w))\neq 0,$
and hence $\gamma$ is not yet closed on $Q$.
In this case, we note that
$\mathbf{d}\cdot \mathbf{d}\gamma= \mathbf{d}^2 \gamma =0, $
and hence the $\mathbf{d}\gamma$ is a closed two-form on $Q$.
Thus, we can construct a magnetic symplectic form on $T^*Q$,
that is, $\omega^B= \omega - \pi_Q^* B= \omega+ \pi_Q^*(\mathbf{d}\gamma), $
where $B=- \mathbf{d}\gamma$, and $\omega$ is the canonical symplectic form on $T^*Q$,
and $\pi_Q^*: T^*Q \rightarrow T^*T^*Q $.
In this case, for any $x, y \in TQ, $ we have that $(\mathbf{d}\gamma +B)(x, y)=0$,
and hence for any $v, w \in TT^* Q, $ we have
$(\mathbf{d}\gamma +B)(T\pi_{Q}(v),T\pi_{Q}(w))=0, $
that is, the one-form $\gamma: Q \rightarrow T^*Q $ satisfies the condition
$\mathbf{d}\gamma=-B$ with respect to $T\pi_{Q}:
TT^* Q \rightarrow TQ. $ Thus, we can construct a CMH system
$(T^*Q, \omega^B, H,F,W )$,
its dynamical vector field with a control law $u$ is given by
$X_{(T^\ast Q,\omega^B,H,F,u)}
=X^B_{H}+\textnormal{vlift}(F)^B +\textnormal{vlift}(u)^B$,
where $X^B_{H}$ satisfies the magnetic Hamiltonian equation, that is,
$\mathbf{i}_{X^B_{H} }\omega^B= \mathbf{d}H $,
and $\textnormal{vlift}(F)^B=\textnormal{vlift}(F)X^B_H$,
$\textnormal{vlift}(u)^B=\textnormal{vlift}(u)X^B_H.$
In this case, by using Lemma 3.4 and the dynamical vector field
$X_{(T^\ast Q,\omega^B,H,F,u)}$, from Theorem 3.5
we can obtain the following theorem.
\begin{theo}
For a given RCH system $(T^*Q,\omega,H,F,W)$ with
the canonical symplectic form $\omega$ on $T^*Q$,
assume that the one-form $\gamma: Q
\rightarrow T^*Q$ is not closed with respect to
$T\pi_Q: TT^* Q \rightarrow TQ. $ Construct
a magnetic symplectic form on $T^*Q$,
$\omega^B= \omega - \pi_Q^*B, $ where
$B=- \mathbf{d}\gamma,$
and a CMH system $(T^*Q, \omega^B, H, F, W )$.
Denote $\tilde{X}^\gamma = T\pi_{Q}\cdot \tilde{X} \cdot \gamma$,
where $\tilde{X}= X_{(T^\ast Q,\omega^B,H,F,u)}$ is the dynamical vector field
of the CMH system $(T^*Q,\omega^B,H,F,W)$ with a control law $u$.
Then the one-form $\gamma$ is just a solution of the Type I of
Hamilton-Jacobi equation $T\gamma\cdot \tilde{X}^\gamma= X^B_H\cdot \gamma ,$
for the CMH system with a control law $u$.
\end{theo}

Next, for any symplectic map $\varepsilon: T^* Q \rightarrow T^* Q $
with respect to the magnetic symplectic form $\omega^B$,
we can prove the following Type II of geometric
Hamilton-Jacobi theorem for the CMH system $(T^*Q,\omega^B,H,F,W)$. For convenience,
the maps involved in the following theorem and its proof are shown
in Diagram-2.
\begin{center}
\hskip 0cm \xymatrix{ & T^* Q \ar[r]^{\varepsilon}
& T^*Q \ar[d]_{X^B_{H\cdot \varepsilon}}
\ar[dr]^{\tilde{X}^\varepsilon} \ar[r]^{\pi_Q}
& Q \ar[r]^{\gamma} & T^*Q \ar[d]^{\tilde{X}} \\
&  & T(T^*Q) & TQ \ar[l]_{T\gamma} & T(T^* Q)\ar[l]_{T\pi_Q}}
\end{center}
$$\mbox{Diagram-2}$$
\begin{theo}
(Type II of Hamilton-Jacobi Theorem for a CMH System)
For the CMH system $(T^*Q,\omega^B,H,F,W)$ with the
magnetic symplectic form $\omega^B= \omega- \pi_Q^*B $ on $T^*Q$,
where $\omega$ is the canonical symplectic form on $T^* Q$
and $B$ is a closed two-form on $Q$,
assume that $\gamma: Q \rightarrow T^*Q$ is an one-form on $Q$, and
$\lambda=\gamma\cdot\pi_{Q}: T^* Q \rightarrow T^* Q $, and for any
symplectic map $\varepsilon: T^* Q \rightarrow T^* Q $ with respect to $\omega^B$,
denote by $\tilde{X}^\varepsilon = T\pi_{Q}\cdot \tilde{X} \cdot \varepsilon$,
where $\tilde{X}=X_{(T^\ast Q,\omega^B,H,F,u)}$
is the dynamical vector field of the CMH system
$(T^*Q,\omega^B,H,F,W)$ with a control law $u$.
Then $\varepsilon$ is a solution of the equation
$T\varepsilon\cdot X^B_{H\cdot\varepsilon}= T\lambda \cdot \tilde{X} \cdot \varepsilon,$
if and only if it is a solution of the equation $T\gamma\cdot \tilde{X}^\varepsilon= X^B_H\cdot
\varepsilon, $ where $X^B_H$ and $ X^B_{H\cdot\varepsilon} \in
TT^*Q $ are the magnetic Hamiltonian vector fields of the functions $H$ and $H\cdot\varepsilon:
T^*Q\rightarrow \mathbb{R}, $ respectively.
The equation $T\gamma\cdot \tilde{X}^\varepsilon= X^B_H\cdot
\varepsilon ,$ is called the Type II of Hamilton-Jacobi equation
for the CMH system $(T^*Q,\omega^B,H,F,W)$ with a control law $u$.
\end{theo}

\noindent{\bf Proof: }
Since $\tilde{X}=X_{(T^\ast Q,\omega^B,H,F,u)}
=X^B_H +\textnormal{vlift}(F)^B+\textnormal{vlift}(u)^B, $ and
$T\pi_{Q}\cdot \textnormal{vlift}(F)^B=T\pi_{Q}\cdot \textnormal{vlift}(u)^B=0, $
then we have that $T\pi_{Q}\cdot \tilde{X}\cdot \varepsilon
=T\pi_{Q}\cdot X^B_H\cdot \varepsilon. $
If we take that $v= X^B_H\cdot \varepsilon \in TT^* Q, $ and for
any $w \in TT^* Q, \; T\lambda(w)\neq 0, $ from Lemma 3.4(ii) we have that
\begin{align*}
&\omega^B(T\gamma \cdot \tilde{X}^\varepsilon, \; w)
= \omega^B(T\gamma \cdot T\pi_Q\cdot \tilde{X}\cdot \varepsilon, \; w)
= \omega^B(T\gamma \cdot T\pi_Q\cdot X^B_H\cdot \varepsilon, \; w)\\ &
= \omega^B(T(\gamma \cdot \pi_Q)\cdot X^B_H\cdot \varepsilon, \; w)
= \omega^B(X^B_H\cdot \varepsilon, \; w-T(\gamma \cdot \pi_Q)\cdot w)
-(\mathbf{d}\gamma+B)(T\pi_{Q}(X^B_H\cdot \varepsilon), \; T\pi_{Q}(w))\\
& =\omega^B(X^B_H\cdot \varepsilon, \; w) - \omega^B(X^B_H\cdot \varepsilon, \;
T\lambda \cdot w)+\lambda^*\omega^B(X^B_H\cdot \varepsilon, \; w)\\
& =\omega^B(X^B_H\cdot \varepsilon, \; w) - \omega^B(X^B_H\cdot \varepsilon, \;
T\lambda \cdot w)+ \omega^B(T\lambda \cdot X^B_H\cdot \varepsilon, \; T\lambda \cdot w).
\end{align*}
Because $\varepsilon: T^* Q
\rightarrow T^* Q $ is symplectic with respect to $\omega^B$,
and hence $ X^B_H\cdot \varepsilon= T\varepsilon \cdot X^B_{H\cdot\varepsilon}, $
along $\varepsilon$. Note that
$T\lambda \cdot X^B_H\cdot \varepsilon=T\gamma \cdot
T\pi_Q\cdot X^B_H\cdot \varepsilon=T\gamma \cdot
T\pi_Q\cdot \tilde{X}\cdot \varepsilon=T\lambda\cdot \tilde{X}\cdot \varepsilon.$
From the above arguments, we can obtain that
\begin{align*}
&\omega^B(T\gamma \cdot \tilde{X}^\varepsilon, \; w)- \omega^B(X^B_H\cdot \varepsilon, \; w)\\
& =- \omega^B(X^B_H\cdot \varepsilon, \; T\lambda \cdot w)
+ \omega^B(T\lambda \cdot X^B_H\cdot \varepsilon, \; T\lambda \cdot w)\\
& =-\omega^B(T\varepsilon \cdot X^B_{H\cdot\varepsilon}, \; T\lambda \cdot w)
+ \omega^B(T\lambda \cdot \tilde{X}\cdot \varepsilon, \; T\lambda \cdot w)\\
& = \omega^B(T\lambda \cdot \tilde{X}\cdot \varepsilon
-T\varepsilon \cdot X^B_{H\cdot\varepsilon}, \; T\lambda \cdot w).
\end{align*}
Because the magnetic symplectic form $\omega^B$ is non-degenerate,
it follows that $T\gamma\cdot \tilde{X}^\varepsilon= X^B_H\cdot
\varepsilon ,$ is equivalent to $T\varepsilon \cdot X^B_{H\cdot\varepsilon}
= T\lambda\cdot \tilde{X}\cdot \varepsilon $.
Thus, $\varepsilon$ is a solution of the equation
$T\varepsilon\cdot X^B_{H\cdot\varepsilon}= T\lambda \cdot \tilde{X} \cdot\varepsilon,$
if and only if it is a solution of the Type II of Hamilton-Jacobi equation
$T\gamma\cdot \tilde{X}^\varepsilon= X^B_H\cdot \varepsilon .$
\hskip 0.3cm $\blacksquare$

\begin{rema}
It is worthy of noting that,
the Type I of Hamilton-Jacobi equation
$T\gamma\cdot \tilde{X}^\gamma= X^B_H \cdot \gamma ,$
is the equation of the differential one-form $\gamma$; and
the Type II of Hamilton-Jacobi equation $T\gamma\cdot \tilde{X}^\varepsilon
= X^B_H \cdot \varepsilon ,$ is the equation of
the symplectic diffeomorphism map $\varepsilon$.
If the external force and control of a CMH
system $(T^*Q,\omega^B,H,F,W)$ are both zeros, that is, $F=0 $
and $W=\emptyset$, in this case the CMH system
is just a magnetic Hamiltonian system $(T^*Q,\omega^B,H)$, and from the proofs of
the above Theorem 3.5-3.7, we can obtain two types of Hamilton-Jacobi
equation for the associated magnetic Hamiltonian system, that is,  Theorem 3.5-3.7
in Wang \cite{wa21c}.
Thus, Theorem 3.5-3.7 can be regarded as an extension of two types
of Hamilton-Jacobi equation for a magnetic Hamiltonian system
to that for the system with external force and control. Moreover,
if $B=0$, in this case the magnetic symplectic form $\omega^B$
is just the canonical symplectic form $\omega$ on $T^*Q$, and the
condition that the one-form $\gamma: Q \rightarrow T^*Q $ satisfies the condition
$\mathbf{d}\gamma=-B $ with respect to $T\pi_{Q}:
TT^* Q \rightarrow TQ, $ becomes that  $\gamma $ is closed with respect to
$T\pi_Q: TT^* Q \rightarrow TQ.$ Thus, from above Theorem 3.5 and Theorem 3.7,
we can obtain Theorem 2.5 and Theorem 2.6 in Wang \cite{wa17}.
Thus, Theorem 3.5 and Theorem 3.7 can be regarded as an extension of two types of Hamilton-Jacobi
equation for a canonical Hamiltonian system to that for the system
with magnetic, external force and control.
\end{rema}

\section{CMH-equivalence and the Solutions of Hamilton-Jacobi Equations}

In the following we first give the definition of CMH-equivalence
for the CMH systems, then prove that the solutions of corresponding to
Hamilton-Jacobi equations leave invariant under the conditions of
CMH-equivalence, if the associated
magnetic Hamiltonian systems are equivalent.
This result describes the relationship between the CMH-equivalence
for the CMH systems and the solutions of the associated Hamilton-Jacobi
equations.\\

For two given Hamiltonian systems $(T^\ast
Q_i,\omega_i,H_i),$ $ i= 1,2,$ we say them to be
equivalent, if there exists a
diffeomorphism $\varphi: Q_1\rightarrow Q_2$, such that
their Hamiltonian vector fields $X_{H_i}, \; i=1,2 $ satisfy
the condition $X_{H_1}\cdot \varphi^\ast
=T(\varphi^\ast) X_{H_2}$, where the map
$\varphi^\ast= T^\ast \varphi:T^\ast Q_2\rightarrow T^\ast Q_1$
is the cotangent lifted map of $\varphi$, and the
map $T(\varphi^\ast):TT^\ast Q_2\rightarrow TT^\ast Q_1$ is the
tangent map of $\varphi^\ast$. From Marsden and Ratiu \cite{mara99},
we know that the condition $X_{H_1}\cdot \varphi^\ast
=T(\varphi^\ast) X_{H_2}$ is equivalent the fact that
the map $\varphi^\ast:T^\ast Q_2\rightarrow T^\ast Q_1$
is symplectic with respect to the canonical symplectic forms
on $T^*Q_i, \; i=1,2.$ In the same way, for two given magnetic
Hamiltonian systems $(T^\ast Q_i,\omega^B_i,H_i),$ $ i= 1,2,$
we say them to be equivalent, if there exists a
diffeomorphism $\varphi: Q_1\rightarrow Q_2$, which
is symplectic with respect to their magnetic symplectic forms,
such that their magnetic Hamiltonian vector fields $X^B_{H_i}, \; i=1,2 $ satisfy
the condition $X^B_{H_1}\cdot \varphi^\ast =T(\varphi^\ast) X^B_{H_2}$. \\

For two given CMH systems $(T^\ast
Q_i,\omega^B_i,H_i,F_i,W_i),$ $ i= 1,2,$ we also want to define
their equivalence, that is, to look for a diffeomorphism
$\varphi: Q_1\rightarrow Q_2$, such that
$X_{(T^\ast Q_1,\omega^B_1,H_1,F_1,W_1)}\cdot \varphi^\ast
=T(\varphi^\ast) X_{(T^\ast Q_2,\omega^B_2,H_2,F_2,W_2)}$. But,
it is worthy of noting that, when a CMH system is given,
the force map $F$ is determined,
but the feedback control law $u: T^\ast Q\rightarrow W$
could be chosen. In order to describe the feedback control law to
modify the structure of the CMH system, the controlled magnetic Hamiltonian matching
conditions and CMH-equivalence are induced as follows.
\begin{defi}
(CMH-equivalence) Suppose that we have two CMH systems $(T^\ast
Q_i,\omega^B_i,H_i,F_i,W_i),$ $ i= 1,2,$ we say them to be
CMH-equivalent, or simply, $(T^\ast
Q_1,\omega^B_1,H_1,F_1,W_1)\stackrel{CMH}{\sim}\\ (T^\ast
Q_2,\omega^B_2,H_2,F_2,W_2)$, if there exists a
diffeomorphism $\varphi: Q_1\rightarrow Q_2$, such that the
following controlled magnetic Hamiltonian matching conditions hold:

\noindent {\bf CMH-1:} The control subsets $W_i, \; i=1,2$ satisfy
the condition $W_1=\varphi^\ast (W_2),$ where the map
$\varphi^\ast= T^\ast \varphi:T^\ast Q_2\rightarrow T^\ast Q_1$
is cotangent lifted map of $\varphi$.

\noindent {\bf CMH-2:} For each control law $u_1:
T^\ast Q_1 \rightarrow W_1, $ there exists the control law $u_2:
T^\ast Q_2 \rightarrow W_2, $  such that the two
closed-loop dynamical systems produce the same dynamical vector fields, that is,
$X_{(T^\ast Q_1,\omega^B_1,H_1,F_1,u_1)}\cdot \varphi^\ast
=T(\varphi^\ast) X_{(T^\ast Q_2,\omega^B_2,H_2,F_2,u_2)}$,
where the map $T(\varphi^\ast):TT^\ast Q_2\rightarrow TT^\ast Q_1$
is the tangent map of $\varphi^\ast$.
\end{defi}

From the expression (2.1) of the dynamical vector field of the CMH system
and the condition $X_{(T^\ast Q_1,\omega^B_1,H_1,F_1,u_1)}\cdot \varphi^\ast
=T(\varphi^\ast) X_{(T^\ast Q_2,\omega^B_2,H_2,F_2,u_2)}$, we have that
$$
(X^B_{H_1}+\textnormal{vlift}(F_1)X^B_{H_1}+\textnormal{vlift}(u_1)X^B_{H_1})\cdot \varphi^\ast
= T(\varphi^\ast)[X^B_{H_2}+\textnormal{vlift}(F_2)X^B_{H_2}+\textnormal{vlift}(u_2)X^B_{H_2}].
$$
By using the notation of vertical lift map of a vector along a fiber,
for $\alpha_x \in T_x^\ast Q_2, \; x \in Q_2$, we have that
\begin{align*}
T(\varphi^\ast)\textnormal{vlift}(F_2)X^B_{H_2}(\alpha_x)
&=T(\varphi^\ast)\textnormal{vlift}((TF_2X^B_{H_2})(F_2(\alpha_x)), \alpha_x)\\
&=\textnormal{vlift}(T(\varphi^\ast)\cdot TF_2\cdot T(\varphi_\ast)X^B_{H_2}
(\varphi^\ast F_2 \varphi_\ast(\varphi^\ast \alpha_x)),\varphi^\ast \alpha)\\
&=\textnormal{vlift}(T(\varphi^\ast F_2 \varphi_\ast)X^B_{H_2}
(\varphi^\ast F_2 \varphi_\ast(\varphi^\ast \alpha_x)),\varphi^\ast \alpha)\\
&=\textnormal{vlift}(\varphi^\ast
F_2\varphi_\ast)X^B_{H_2}(\varphi^\ast \alpha_x),
\end{align*}
where the map $\varphi_\ast=(\varphi^{-1})^\ast: T^\ast Q_1\rightarrow T^\ast Q_2$.
In the same way, we have that
$T(\varphi^\ast)\textnormal{vlift}(u_2)X^B_{H_2}=\textnormal{vlift}(\varphi^\ast
u_2\varphi_\ast)X^B_{H_2}\cdot \varphi^\ast$. Note that
$\textnormal{vlift}(F)^B=\textnormal{vlift}(F)X^B_H$,
and $\textnormal{vlift}(u)^B=\textnormal{vlift}(u)X^B_H.$,
and hence we have that the explicit relation
between the two control laws $u_i \in W_i, \; i=1,2$ in {\bf RCH-2} is given by
\begin{align}
& (\textnormal{vlift}(u_1)^B -\textnormal{vlift}(\varphi^\ast u_2\varphi_\ast)^B)\cdot \varphi^\ast \nonumber \\
& = -X^B_{H_1}\cdot \varphi^\ast +T(\varphi^\ast) (X^B_{H_2})+
(-\textnormal{vlift}(F_1)^B+\textnormal{vlift}(\varphi^\ast F_2
\varphi_\ast)^B)\cdot \varphi^\ast.
\label{4.1}\end{align}
From the above relation (4.1) we know that, when two CMH systems $(T^\ast
Q_i,\omega^B_i,H_i,F_i,W_i),$ $ i= 1,2,$ are CMH-equivalent with respect to $\varphi^*$, the
associated magnetic Hamiltonian systems $(T^\ast Q_i,\omega^B_i,H_i),$
$ i= 1,2,$ may not be equivalent with respect to $\varphi^*$.\\

On the other hand,
note that the magnetic vector field $X^0= X^B_H-X_H, $
from (2.3) and (4.1) we have that
\begin{align*}
& (\textnormal{vlift}(u_1)^B -\textnormal{vlift}(\varphi^\ast u_2\varphi_\ast)^B)\cdot \varphi^\ast \\
& = -(X_{H_1}+X^0_1)\cdot \varphi^\ast +T(\varphi^\ast) (X_{H_2}+X^0_2)+
(-\textnormal{vlift}(F_1)^B+\textnormal{vlift}(\varphi^\ast F_2
\varphi_\ast)^B)\cdot \varphi^\ast.
\end{align*}
and hence we have that
\begin{align}
& (X^0_1+\textnormal{vlift}(F_1)^B+\textnormal{vlift}(u_1)^B )\cdot \varphi^\ast \nonumber \\
& = -X_{H_1}\cdot \varphi^\ast +T(\varphi^\ast) X_{H_2}+T(\varphi^\ast) (X^0_2
+\textnormal{vlift}(F_2)^B+\textnormal{vlift}(u_2)^B).
\label{4.2}\end{align}
If the associated canonical Hamiltonian systems $(T^\ast
Q_i,\omega_i,H_i),$ $ i= 1,2,$ are also equivalent with respect to $\varphi^*$, that is,
$T(\varphi^*) \cdot X_{H_2}= X_{H_1}\cdot \varphi^*$. In this case,
from (2.4) and (4.2), we have the following theorem.
\begin{theo}
Suppose that two CMH systems $(T^\ast Q_i,\omega^B_i,H_i,F_i,W_i)$,
$i=1,2,$ are CMH-equivalent with respect to $\varphi^*$,
and the associated canonical Hamiltonian systems $(T^\ast
Q_i,\omega_i,H_i),$ $ i= 1,2,$ are also equivalent with respect to $\varphi^*$.
Then we have the following fact that, if one system satisfies the magnetic
vanishing condition, then another CMH-equivalent system must satisfy
the associated magnetic vanishing condition.
\end{theo}

Moreover, if considering the CMH-equivalence of the CMH systems,
we can prove the following Theorem 4.3, which states that
the solutions of two types of Hamilton-Jacobi equations for the CMH systems leave
invariant under the conditions of CMH-equivalence, if the associated
magnetic Hamiltonian systems are equivalent.
\begin{theo}
Suppose that two CMH systems $(T^\ast Q_i,\omega^B_i,H_i,F_i,W_i)$,
$i=1,2,$ are CMH-equivalent with an equivalent map $\varphi: Q_1
\rightarrow Q_2 $, and the associated magnetic Hamiltonian systems $(T^\ast
Q_i,\omega^B_i,H_i),$ $ i= 1,2,$ are also equivalent with respect to $\varphi^*$,
under the hypotheses and notations of Theorem 3.5,
Theorem 3.7, we have that\\

\noindent $(\mathrm{i})$ If the one-form $\gamma_2: Q_2 \rightarrow T^* Q_2$
satisfies the condition that $\mathbf{d}\gamma_2=-B_2 $ with
respect to $T\pi_{Q_2}: TT^* Q_2 \rightarrow TQ_2, $ then $\gamma_1=
\varphi^* \cdot \gamma_2\cdot \varphi: Q_1 \rightarrow T^* Q_1 $ satisfies also the condition that
$\mathbf{d}\gamma_1=-B_1 $ with respect to $T\pi_{Q_1}:
TT^* Q_1 \rightarrow TQ_1, $ and hence it is
a solution of the Type I of Hamilton-Jacobi equation for the CMH system
$(T^*Q_1,\omega^B_1,H_1,F_1,W_1). $ Vice versa;\\

\noindent $(\mathrm{ii})$ If the symplectic map
$\varepsilon_2: T^*Q_2\rightarrow T^* Q_2$ with respect to $\omega^B_2$
is a solution of the Type II of Hamilton-Jacobi equation for the CMH system
$(T^*Q_2,\omega^B_2,H_2, F_2,W_2)$, then $\varepsilon_1=
\varphi^* \cdot \varepsilon_2\cdot \varphi_*: T^*Q_1 \rightarrow
T^* Q_1 $ is a symplectic map with respect to $\omega^B_1$,
and it is a solution of the Type II of Hamilton-Jacobi equation
for the CMH system $(T^*Q_1,\omega^B_1,H_1,F_1,W_1). $ Vice versa.
\end{theo}

\noindent{\bf Proof: }
We first prove the assertion $(\mathrm{i})$.
If two given CMH systems $(T^\ast Q_i,\omega^B_i,H_i,F_i,W_i)$, $i=1,2,$ are
CMH-equivalent with an equivalent map $\varphi: Q_1 \rightarrow Q_2
$, from the definition of CMH-equivalence, we know that
for each control law $u_1:
T^\ast Q_1 \rightarrow W_1, $ there exists the control law $u_2:
T^\ast Q_2 \rightarrow W_2, $  such that the two
closed-loop dynamical systems produce the same dynamical vector fields, that is,
$\tilde{X}_1 \cdot \varphi^\ast =T(\varphi^\ast)\cdot \tilde{X}_2, $
where $\tilde{X}_i= X_{(T^\ast
Q_i,\omega^B_i, H_i, F_i, u_i)}, \; i=1,2.$ From the following
commutative Diagram-3:
\[
\begin{CD}
 Q_1 @> \gamma_1 >> T^* Q_1 @> \tilde{X}_1 >> TT^* Q_1 @> T\pi_{Q_1} >> TQ_1 \\
@V \varphi VV @A \varphi^* AA @A T\varphi^* AA @V T\varphi VV \\
 Q_2 @> \gamma_2 >> T^* Q_2 @> \tilde{X}_2 >> TT^* Q_2 @> T\pi_{Q_2} >> TQ_2
\end{CD}
\]
$$\mbox{Diagram-3}$$
we have that $\gamma_1= \varphi^* \cdot \gamma_2\cdot \varphi, $
$\mathbf{d}\gamma_1 = \varphi^* \cdot \mathbf{d}\gamma_2 \cdot
\varphi, $ $B_1= \varphi^* \cdot B_2 \cdot \varphi, $ and
$T\varphi \cdot T\pi_{Q_1} \cdot T\varphi^*= T\pi_{Q_2}. $
For $x\in Q_1, $ and $v, \; w \in TT^*_x Q_1, $ then $\varphi(x) \in Q_2 $
and $T\varphi_*(v), \; T\varphi_*(w) \in TT^*_{\varphi(x)}Q_2. $
Since the one-form $\gamma_2: Q_2 \rightarrow T^* Q_2$ satisfies
the condition that $\mathbf{d}\gamma_2=-B_2 $ with
respect to $T\pi_{Q_2}: TT^* Q_2 \rightarrow TQ_2, $ then
$$(\mathbf{d}\gamma_2+B_2) (T\pi_{Q_2}\cdot T\varphi_*(v), \;
T\pi_{Q_2}\cdot T\varphi_*(w) )(\varphi(x))=0.$$
Thus,
\begin{align*}
& (\mathbf{d}\gamma_1+B_1) (T\pi_{Q_1}(v), \; T\pi_{Q_1}(w)) (x) \\
&= \varphi^* \cdot (\mathbf{d}\gamma_2+B_2) \cdot
\varphi (T\pi_{Q_1}(v), \; T\pi_{Q_1}(w))(x)\\
&= ( \mathbf{d}\gamma_2+B_2) (T\varphi \cdot T\pi_{Q_1}(v),\; T\varphi \cdot T\pi_{Q_1}(w))(\varphi(x))\\
&= (\mathbf{d}\gamma_2+B_2) (T\varphi \cdot T\pi_{Q_1} \cdot T\varphi^*\cdot T(\varphi^{-1})^*(v), \;
T\varphi \cdot T\pi_{Q_1} \cdot T\varphi^*\cdot T(\varphi^{-1})^*(w)) (\varphi(x))\\
&= (\mathbf{d}\gamma_2+B_2) (T\pi_{Q_2}\cdot T\varphi_*(v), \;
T\pi_{Q_2}\cdot T\varphi_*(w) )(\varphi(x))=0,
\end{align*}
that is, the one-form $\gamma_1=
\varphi^* \cdot \gamma_2\cdot \varphi: Q_1 \rightarrow T^* Q_1 $ satisfies
the condition that $\mathbf{d}\gamma_1=-B_1 $ with respect to $T\pi_{Q_1}:
TT^* Q_1 \rightarrow TQ_1. $ Moreover, from Theorem 3.5 we know that,
the one-form $\gamma_2$ is
a solution of the Type I of Hamilton-Jacobi equation for the CMH system
$(T^*Q_2,\omega^B_2, H_2, F_2, W_2), $ that is,
$T\gamma_2\cdot \tilde{X}_2^{\gamma_2}= X^B_{H_2}\cdot \gamma_2 ,$
where $\tilde{X}_i^{\gamma_i}= T\pi_{Q_i}\cdot \tilde{X}_i \cdot \gamma_i, \; i=1,2.$
Hence,
\begin{align*}
T\gamma_1\cdot \tilde{X}_1^{\gamma_1}& =T(\varphi^* \cdot \gamma_2\cdot \varphi)\cdot
T\pi_{Q_1}\cdot \tilde{X}_1 \cdot \gamma_1 \\
& =T(\varphi^*) \cdot T\gamma_2\cdot T\varphi \cdot T\pi_{Q_1}\cdot \tilde{X}_1 \cdot (\varphi^* \cdot \gamma_2\cdot \varphi) \\
& =T(\varphi^*) \cdot T\gamma_2\cdot T\varphi \cdot T\pi_{Q_1}\cdot (T(\varphi^\ast) \cdot \tilde{X}_2) \cdot \gamma_2\cdot \varphi \\
& =T(\varphi^*) \cdot T\gamma_2\cdot (T\pi_{Q_2}\cdot \tilde{X}_2 \cdot \gamma_2)\cdot \varphi
=T(\varphi^*) \cdot T\gamma_2\cdot \tilde{X}_2^{\gamma_2}\cdot \varphi \\
& =T(\varphi^*) \cdot X^B_{H_2}\cdot \gamma_2 \cdot \varphi =X^B_{H_1}\cdot \varphi^*\cdot \gamma_2 \cdot \varphi
=X^B_{H_1}\cdot \gamma_1,
\end{align*}
where we have used that $T(\varphi^*) \cdot X^B_{H_2}= X^B_{H_1}\cdot \varphi^*$,
because the associated magnetic Hamiltonian systems $(T^\ast
Q_i,\omega^B_i,H_i),$ $ i= 1,2,$ are equivalent with respect to $\varphi^*$.
Thus, the one-form $\gamma_1=
\varphi^* \cdot \gamma_2\cdot \varphi $ is
a solution of the Type I of Hamilton-Jacobi equation for the CMH system
$(T^*Q_1,\omega^B_1,H_1,F_1,W_1). $
Note that the map $\varphi: Q_1 \rightarrow Q_2 $ is a diffeomorphism,
and $\varphi^*: T^* Q_2 \rightarrow T^* Q_1 $ is a symplectic isomorphisms,
vice versa. It follows that the assertion $(\mathrm{i})$ of Theorem 4.3 holds.\\

Next, we prove the assertion $(\mathrm{ii})$.
From the following commutative Diagram-4:
\[
\begin{CD}
 Q_1 @> \gamma_1 >> T^* Q_1 @> \varepsilon_1 >> T^* Q_1 @> \tilde{X}_1 >> TT^* Q_1 @> T\pi_{Q_1} >> TQ_1 \\
@V \varphi VV @V \varphi_* VV @ A \varphi^* AA @A T\varphi^* AA @V T\varphi VV \\
 Q_2 @> \gamma_2 >> T^* Q_2 @> \varepsilon_2 >> T^* Q_2 @> \tilde{X}_2 >> TT^* Q_2 @> T\pi_{Q_2} >> TQ_2
\end{CD}
\]
$$\mbox{Diagram-4}$$
we have that $\varepsilon_1=
\varphi^* \cdot \varepsilon_2\cdot \varphi_*: T^*Q_1 \rightarrow
T^* Q_1 .$ Since
$\varepsilon_2: T^*Q_2\rightarrow T^* Q_2$ is symplectic with respect to $\omega^B_2$, then
for $x\in Q_1 $, $v, \; w \in TT^*_x Q_1, $ and $\varphi(x) \in Q_2 $,
$T\varphi_*(v), \; T\varphi_*(w) \in TT^*_{\varphi(x)}Q_2, $ we have that
$\varepsilon_2^*\cdot \omega^B_2(T\varphi_*(v), \; T\varphi_*(w))(\varphi(x))=\omega^B_2 (T\varphi_*(v), \; T\varphi_*(w))(\varphi(x)).$
Note that the associated magnetic Hamiltonian systems $(T^\ast
Q_i,\omega^B_i,H_i),$ $ i= 1,2,$ are equivalent with respect to $\varphi^*$, and hence
$\varphi^\ast: T^\ast Q_2\rightarrow T^\ast Q_1$ is symplectic
with respect to their magnetic symplectic forms, that is,
$(\varphi^*)^*\omega^B_1(v,w)(x)= \omega^B_2 (T\varphi_*(v), \; T\varphi_*(w))(\varphi(x))$,
then we have that
\begin{align*}
\varepsilon_1^*\cdot \omega^B_1(v,w)(x)& = (\varphi^* \cdot \varepsilon_2\cdot \varphi_*)^*\omega^B_1(v,w)(x)
= (\varphi_*)^* \cdot \varepsilon_2^*\cdot (\varphi^*)^*\omega^B_1(v,w)(x)\\
& =(\varphi_*)^* \cdot \varepsilon_2^*\cdot \omega^B_2(T\varphi_*(v), \; T\varphi_*(w))(\varphi(x))
=((\varphi^{-1})^*)^*\cdot \omega^B_2(T\varphi_*(v), \; T\varphi_*(w))(\varphi(x))\\
& =\omega^B_1(T(\varphi^{-1})_*\cdot T\varphi_*(v), \; T(\varphi^{-1})_*\cdot T\varphi_*(w))(\varphi^{-1}\cdot \varphi(x))
= \omega^B_1(v,w)(x),
\end{align*}
that is, the map $\varepsilon_1=
\varphi^* \cdot \varepsilon_2\cdot \varphi_*: T^*Q_1 \rightarrow
T^* Q_1 $ is symplectic map with respect to $\omega^B_1$. Moreover, because
the symplectic map $\varepsilon_2: T^*Q_2\rightarrow T^* Q_2$ is
a solution of the Type II of Hamilton-Jacobi equation for the CMH system
$(T^*Q_2,\omega^B_2,H_2, F_2,W_2)$, that is,
$T\gamma_2\cdot \tilde{X}_2^{\varepsilon_2}= X^B_{H_2}\cdot
\varepsilon_2 ,$ where $\tilde{X}_i^{\varepsilon_i}= T\pi_{Q_i}\cdot \tilde{X}_i \cdot \varepsilon_i, \; i=1,2.$
Hence, we have that
\begin{align*}
T\gamma_1\cdot \tilde{X}_1^{\varepsilon_1}& =T(\varphi^* \cdot \gamma_2\cdot \varphi)\cdot
T\pi_{Q_1}\cdot \tilde{X}_1 \cdot \varepsilon_1 \\
& =T(\varphi^*) \cdot T\gamma_2\cdot T\varphi \cdot T\pi_{Q_1}\cdot \tilde{X}_1 \cdot (\varphi^* \cdot \varepsilon_2\cdot \varphi_*) \\
& =T(\varphi^*) \cdot T\gamma_2\cdot T\varphi \cdot T\pi_{Q_1}\cdot (T(\varphi^\ast) \cdot \tilde{X}_2) \cdot \varepsilon_2\cdot \varphi_* \\
& =T(\varphi^*) \cdot T\gamma_2\cdot (T\pi_{Q_2}\cdot \tilde{X}_2 \cdot \varepsilon_2)\cdot \varphi_*
=T(\varphi^*) \cdot T\gamma_2\cdot \tilde{X}_2^{\varepsilon_2}\cdot \varphi_* \\
& =T(\varphi^*) \cdot X^B_{H_2}\cdot \varepsilon_2 \cdot \varphi_*
=X^B_{H_1}\cdot \varphi^*\cdot \varepsilon_2 \cdot \varphi_*
=X^B_{H_1}\cdot \varepsilon_1,
\end{align*}
that is, the symplectic map $\varepsilon_1=
\varphi^* \cdot \varepsilon_2\cdot \varphi_*$
is a solution of the Type II of Hamilton-Jacobi equation for the CMH
system $(T^*Q_1,\omega^B_1,H_1,F_1,W_1). $
In the same way, because the map $\varphi: Q_1 \rightarrow Q_2 $ is a diffeomorphism,
and $\varphi^\ast: T^\ast Q_2\rightarrow T^\ast Q_1$ is a symplectic isomorphisms,
vice versa. Hence we prove the assertion $(\mathrm{ii})$ of Theorem 4.3.
\hskip 0.3cm  $\blacksquare$\\

\section{Nonholonomic CMH System and Distributional CMH System}

In order to describe the impact of nonholonomic constraints
for Hamilton-Jacobi theory of the dynamics of a CMH system,
in this section we first give some definitions and basic facts
about the nonholonomic constraint and the nonholonomic CMH system.
Moreover, by analyzing carefully the structure for the nonholonomic dynamical
vector field, we give a geometric formulation of distributional CMH system,
which is determined by a non-degenerate distributional two-form induced
from the magnetic symplectic form,
and which will be used in subsequent section.\\

In order to describe the nonholonomic CMH system,
in the following we first give the completeness and regularity
conditions for nonholonomic constraints of a mechanical system,
see Le\'{o}n and Wang \cite{lewa15} and Wang \cite{wa21c, wa22a}. In fact,
in order to describe the dynamics of a nonholonomic mechanical system,
we need some restriction conditions for nonholonomic constraints of
the system. At first, we note that the set of Hamiltonian vector fields
forms a Lie algebra with respect to the Lie bracket, since
$X_{\{f,g\}}=-[X_f, X_g]. $ But, the Lie bracket operator, in
general, may not be closed on the restriction of a nonholonomic
constraint. Thus, we have to give the following completeness
condition for nonholonomic constraints of a system.\\

{\bf $\mathcal{D}$-completeness } Let $Q$ be a smooth manifold and
$TQ$ its tangent bundle. A distribution $\mathcal{D} \subset TQ$ is
said to be {\bf completely nonholonomic} (or bracket-generating) if
$\mathcal{D}$ along with all of its iterated Lie brackets
$[\mathcal{D},\mathcal{D}], [\mathcal{D}, [\mathcal{D},\mathcal{D}]],
\cdots ,$ spans the tangent bundle $TQ$. Moreover, we consider a
nonholonomic mechanical system on $Q$, which is
given by a Lagrangian function $L: TQ \rightarrow \mathbb{R}$
subject to constraints determined by a nonholonomic
distribution $\mathcal{D}\subset TQ$ on the configuration manifold $Q$.
Then the nonholonomic system is said to be {\bf completely nonholonomic},
if the distribution $\mathcal{D} \subset TQ$ determined by the nonholonomic
constraints is completely nonholonomic.\\

{\bf $\mathcal{D}$-regularity } In the following we always assume
that $Q$ is an $n$-dimensional smooth manifold with coordinates $(q^i)$, and $TQ$ its
tangent bundle with coordinates $(q^i,\dot{q}^i)$, and $T^\ast Q$
its cotangent bundle with coordinates $(q^i,p_j)$, which are the
canonical cotangent coordinates of $T^\ast Q$ and $\omega= \sum^n_{i=1}
dq^{i}\wedge dp_{i}$ is canonical symplectic form on $T^{\ast}Q$. If
the Lagrangian $L: TQ \rightarrow \mathbb{R}$ is hyperregular, that
is, the Hessian matrix
$(\partial^2L/\partial\dot{q}^i\partial\dot{q}^j)$ is nondegenerate
everywhere, then the Legendre transformation $FL: TQ \rightarrow T^*
Q$ is a diffeomorphism. In this case the Hamiltonian $H: T^* Q
\rightarrow \mathbb{R}$ is given by $H(q,p)=\dot{q}\cdot
p-L(q,\dot{q}) $ with Hamiltonian vector field $X_H$,
which is defined by the Hamilton's equation
$\mathbf{i}_{X_H}\omega=\mathbf{d}H$, and
$\mathcal{M}=\mathcal{F}L(\mathcal{D})$ is a constraint submanifold
in $T^* Q$. In particular, for the nonholonomic constraint
$\mathcal{D}\subset TQ$, the Lagrangian $L$ is said to be {\bf
$\mathcal{D}$-regular}, if the restriction of Hessian matrix
$(\partial^2L/\partial\dot{q}^i\partial\dot{q}^j)$ on $\mathcal{D}$
is nondegenerate everywhere. Moreover, a nonholonomic system is said
to be {\bf $\mathcal{D}$-regular}, if its Lagrangian $L$ is {\bf
$\mathcal{D}$-regular}. Note that the restriction of a positive
definite symmetric bilinear form to a subspace is also positive
definite, and hence nondegenerate. Thus, for a simple nonholonomic
mechanical system, that is, whose Lagrangian is the total kinetic
energy minus potential energy, it is {\bf $\mathcal{D}$-regular }
automatically.\\

A nonholonomic magnetic Hamiltonian system is the 4-tuple
$(T^\ast Q,\omega^B,\mathcal{D},H)$,
which is a magnetic Hamiltonian system with a
$\mathcal{D}$-completely and $\mathcal{D}$-regularly nonholonomic
constraint $\mathcal{D} \subset TQ$.
A nonholonomic CMH system is the 6-tuple
$(T^\ast Q,\omega^B,\mathcal{D},H, F,W)$,
which is a nonholonomic magnetic Hamiltonian system with
external force $F$ and control $W$, where
$F: T^*Q\rightarrow T^*Q$ is the fiber-preserving map,
and $W\subset T^*Q$ is a fiber submanifold of $T^*Q$.
Under the restriction given by constraint, in general, the dynamical
vector field of a nonholonomic CMH system may not be
magnetic Hamiltonian, however the system is a dynamical system
closely related to a magnetic Hamiltonian system.
In the following we shall derive a distributional CMH system of the nonholonomic
CMH system $(T^*Q,\omega^B,\mathcal{D},H,F,W)$,
by analyzing carefully the structure for the nonholonomic
dynamical vector field similar to the method used in
Le\'{o}n and Wang \cite{lewa15} and Wang \cite{wa21c, wa22a}.\\

We consider that the constraint submanifold
$\mathcal{M}=\mathcal{F}L(\mathcal{D})\subset T^*Q$ and
$i_{\mathcal{M}}: \mathcal{M}\rightarrow T^*Q $ is the inclusion,
the symplectic form $\omega^B_{\mathcal{M}}= i_{\mathcal{M}}^* \omega^B $,
is induced from the magnetic symplectic form $\omega^B$ on $T^* Q$.
We define the distribution $\mathcal{F}$ as the pre-image of the nonholonomic
constraints $\mathcal{D}$ for the map $T\pi_Q: TT^* Q \rightarrow TQ$,
that is, $\mathcal{F}=(T\pi_Q)^{-1}(\mathcal{D})\subset TT^*Q,
$ which is a distribution along $\mathcal{M}$, and
$\mathcal{F}^\circ:=\{\alpha \in T^*T^*Q | <\alpha,v>=0, \; \forall
v\in TT^*Q \}$ is the annihilator of $\mathcal{F}$ in
$T^*T^*Q_{|\mathcal{M}}$. We consider the following nonholonomic
constraints condition
\begin{align} (\mathbf{i}_X \omega^B -\mathbf{d}H) \in \mathcal{F}^\circ,
\;\;\;\;\;\; X \in T \mathcal{M},
\label{5.1} \end{align} from Cantrijn et al.
\cite{calemama99}, we know that there exists an unique nonholonomic
vector field $X_n$ satisfying the above condition $(5.1)$, if the
admissibility condition $\mathrm{dim}\mathcal{M}=
\mathrm{rank}\mathcal{F}$ and the compatibility condition
$T\mathcal{M}\cap \mathcal{F}^\bot= \{0\}$ hold, where
$\mathcal{F}^\bot$ denotes the magnetic symplectic orthogonal of
$\mathcal{F}$ with respect to the magnetic symplectic form
$\omega^B$ on $T^*Q$. In particular, when we consider the Whitney sum
decomposition $T(T^*Q)_{|\mathcal{M}}=T\mathcal{M}\oplus
\mathcal{F}^\bot$ and the canonical projection $P:
T(T^*Q)_{|\mathcal{M}} \rightarrow T\mathcal{M}$,
then we have that $X_n= P(X^B_H)$.\\

From the condition (5.1) we know that the nonholonomic vector field,
in general, may not be magnetic Hamiltonian, because of the restriction
of nonholonomic constraints. But, we hope to study the dynamical
vector field of nonholonomic CMH system by using the similar
method of studying magnetic Hamiltonian vector field.
From Le\'{o}n and Wang \cite{lewa15} and Bates
and $\acute{S}$niatycki \cite{basn93}, by using the
similar method, we can define the
distribution $ \mathcal{K}=\mathcal {F}\cap T\mathcal{M}.$ and
$\mathcal{K}^\bot=\mathcal {F}^\bot\cap T\mathcal{M}, $ where
$\mathcal{K}^\bot$ denotes the magnetic symplectic orthogonal of
$\mathcal{K}$ with respect to the magnetic symplectic form
$\omega^B$, and the admissibility condition $\mathrm{dim}\mathcal{M}=
\mathrm{rank}\mathcal{F}$ and the compatibility condition
$T\mathcal{M}\cap \mathcal{F}^\bot= \{0\}$ hold, then we know that the
restriction of the symplectic form $\omega^B_{\mathcal{M}}$ on
$T^*\mathcal{M}$ fibrewise to the distribution $\mathcal{K}$, that
is, $\omega^B_\mathcal{K}= \tau_{\mathcal{K}}\cdot
\omega^B_{\mathcal{M}}$ is non-degenerate, where $\tau_{\mathcal{K}}$
is the restriction map to distribution $\mathcal{K}$. It is worthy
of noting that $\omega^B_\mathcal{K}$ is not a true two-form on a
manifold, so it does not make sense to speak about it being closed.
We call $\omega^B_\mathcal{K}$ as a distributional two-form to avoid
any confusion. Because $\omega^B_\mathcal{K}$ is non-degenerate as a
bilinear form on each fibre of $\mathcal{K}$, there exists a vector
field $X^B_{\mathcal{K}}$ on $\mathcal{M}$ which takes values in the
constraint distribution $\mathcal{K}$,
such that the distributional magnetic Hamiltonian equation holds, that is,
\begin{align}
\mathbf{i}_{X^B_\mathcal{K}}\omega^B_{\mathcal{K}}=\mathbf{d}H_\mathcal
{K}, \label{5.2} \end{align}
 where $\mathbf{d}H_\mathcal{K}$ is the restriction of
$\mathbf{d}H_\mathcal{M}$ to $\mathcal{K}$,
and the function $H_{\mathcal{K}}$ satisfies
$\mathbf{d}H_{\mathcal{K}}= \tau_{\mathcal{K}}\cdot \mathbf{d}H_{\mathcal {M}}$,
and $H_\mathcal{M}= \tau_{\mathcal{M}}\cdot H$ is the restriction of $H$ to
$\mathcal{M}$. Moreover, from the distributional magnetic Hamiltonian equation (5.2),
we have that $X^B_{\mathcal{K}}= \tau_{\mathcal{K}}\cdot X^B_H.$\\

Moreover, if considering the external force $F$ and control subset $W$,
and define $F^B_\mathcal{K}=\tau_{\mathcal{K}}\cdot \textnormal{vlift}(F_{\mathcal{M}})X^B_H,$
and for a control law $u\in W$,
$u^B_\mathcal{K}= \tau_{\mathcal{K}}\cdot  \textnormal{vlift}(u_{\mathcal{M}})X^B_H,$
where $F_\mathcal{M}= \tau_{\mathcal{M}}\cdot F$ and
$u_\mathcal{M}= \tau_{\mathcal{M}}\cdot u$ are the restrictions of
$F$ and $u$ to $\mathcal{M}$, that is, $F^B_\mathcal{K}$ and $u^B_\mathcal{K}$
are the restrictions of the changes of magnetic Hamiltonian vector field $X^B_H$
under the actions of $F_\mathcal{M}$ and $u_\mathcal{M}$ to $\mathcal{K}$.
Then the 5-tuple $(\mathcal{K},\omega^B_{\mathcal{K}},
H_\mathcal{K}, F^B_\mathcal{K}, u^B_\mathcal{K})$
is a distributional CMH system of the nonholonomic
CMH system $(T^*Q,\omega,\mathcal{D},H,F,W)$ with a control law $u$.
Thus, the geometric formulation of the distributional CMH
system may be summarized as follows.

\begin{defi} (Distributional CMH System)
Assume that the 6-tuple $(T^*Q,\omega^B, \mathcal{D},H,F,W)$ is a
nonholonomic CMH system, where the magnetic symplectic form
$\omega^B= \omega- \pi_Q^*B $ on $T^*Q$,
and $\omega$ is the canonical symplectic form on $T^* Q$
and $B$ is a closed two-form on $Q$, and
$\mathcal{D}\subset TQ$ is a
$\mathcal{D}$-completely and $\mathcal{D}$-regularly nonholonomic
constraint of the system, and the external force
$F: T^*Q\rightarrow T^*Q$ is the fiber-preserving map,
and the control subset $W\subset T^*Q$ is a fiber submanifold of $T^*Q$.
For a control law $u\in W,$ if there exist a distribution $\mathcal{K}$,
an associated non-degenerate distributional two-form
$\omega^B_{\mathcal{K}}$ induced by the magnetic symplectic form $\omega^B$
and a vector field $X^B_\mathcal {K}$ on the
constraint submanifold $\mathcal{M}=\mathcal{F}L(\mathcal{D})\subset
T^*Q$, such that the distributional magnetic Hamiltonian equation
$\mathbf{i}_{X^B_\mathcal{K}}\omega^B_{\mathcal{K}}=\mathbf{d}H_\mathcal
{K}$ holds, where $\mathbf{d}H_\mathcal{K}$ is the restriction of
$\mathbf{d}H_\mathcal{M}$ to $\mathcal{K}$,  and
the function $H_{\mathcal{K}}$ satisfies
$\mathbf{d}H_{\mathcal{K}}= \tau_{\mathcal{K}}\cdot \mathbf{d}H_{\mathcal {M}},$
and $F^B_{\mathcal{K}}=\tau_{\mathcal{K}}\cdot \textnormal{vlift}(F_{\mathcal{M}})X^B_H$,
and $u^B_{\mathcal{K}}=\tau_{\mathcal{K}}\cdot \textnormal{vlift}(u_{\mathcal{M}})X^B_H$
as defined above,
then the 5-tuple $(\mathcal{K},\omega^B_{\mathcal{K}},H_{\mathcal{K}},
F^B_{\mathcal{K}},u^B_{\mathcal{K}})$
is called a distributional CMH system of the nonholonomic
CMH system $(T^*Q,\omega^B,\mathcal{D},H,F,u)$, and $X^B_\mathcal
{K}$ is called a nonholonomic dynamical vector field.
Denote by
\begin{align}\tilde{X}=X^B_{(\mathcal{K},\omega^B_{\mathcal{K}},
H_{\mathcal{K}}, F^B_{\mathcal{K}}, u^B_{\mathcal{K}})}
=X^B_\mathcal {K}+ F^B_{\mathcal{K}}+u^B_{\mathcal{K}}
\label{5.3} \end{align}
is the dynamical vector field of the distributional CMH system
$(\mathcal{K},\omega^B_{\mathcal {K}},H_{\mathcal{K}},
F^B_{\mathcal{K}},u^B_{\mathcal{K}})$, which is the synthetic
of the nonholonomic dynamical vector field $X^B_{\mathcal{K}}$ and
the vector fields $F^B_{\mathcal{K}}$ and $u^B_{\mathcal{K}}$.
Under the above circumstances, we refer to
$(T^*Q,\omega^B,\mathcal{D},H,F,u)$ as a nonholonomic CMH system
with an associated distributional CMH system
$(\mathcal{K},\omega^B_{\mathcal {K}},H_{\mathcal{K}},
F^B_{\mathcal{K}},u^B_{\mathcal{K}})$.
\end{defi}

\begin{rema}
It is worthy of noting that, when $B=0$, in this case the magnetic symplectic form $\omega^B$
is just the canonical symplectic form $\omega$ on $T^*Q$,
and the distributional CMH system
$(\mathcal{K},\omega^B_{\mathcal {K}},H_{\mathcal {K}}, F^B_{\mathcal {K}}, u^B_{\mathcal {K}})$
becomes the distributional RCH system
$(\mathcal{K},\omega_{\mathcal {K}},H_{\mathcal {K}}, F_{\mathcal {K}}, u_{\mathcal {K}})$,
which is given in Wang \cite{wa22a}. Moreover,
if the external force and control of a distributional CMH
system $(\mathcal{K},\omega^B_{\mathcal {K}},H_{\mathcal {K}},
F^B_{\mathcal {K}}, u^B_{\mathcal {K}})$ are both zeros,
that is, $B=0$, $F^B_{\mathcal {K}}=0 $ and $u^B_{\mathcal {K}}=0$, in this case,
the distributional CMH system is just a distributional Hamiltonian system
$(\mathcal{K},\omega_{\mathcal {K}},H_{\mathcal {K}})$,
which is given in Le\'{o}n and Wang \cite{lewa15}.
Thus, the distributional CMH system
can be regarded as an extension of the distributional Hamiltonian system to
the system with magnetic, external force and control.
\end{rema}

For the nonholonomic CMH system $(T^*Q,\omega^B,\mathcal{D},H,F,u)$
with an associated distributional CMH system
$(\mathcal{K},\omega^B_{\mathcal {K}},H_{\mathcal {K}},
F^B_{\mathcal {K}},u^B_{\mathcal {K}})$,
the magnetic vector field $X^0= X^B_H-X_H, $
which is determined by the magnetic equation
$\mathbf{i}_{X^0}\omega=\mathbf{i}_{X^B_H}( \pi_Q^*B)$
on $T^*Q$. Denote by $X^0_\mathcal {K}=\tau_\mathcal {K}(X^0)
= \tau_\mathcal {K}(X^B_H)- \tau_\mathcal {K}(X_H)=X^B_\mathcal {K}- X_\mathcal {K},$
from the expression (5.3) of the dynamical
vector field of the distributional CMH system
$(\mathcal{K},\omega^B_{\mathcal {K}},H_{\mathcal{K}},
F^B_{\mathcal{K}},u^B_{\mathcal{K}})$, we have that
\begin{align}\tilde{X}
=X^B_\mathcal {K}+ F^B_{\mathcal{K}}+u^B_{\mathcal{K}}
=X_\mathcal {K}+ X^0_{\mathcal{K}}+ F^B_{\mathcal{K}}+u^B_{\mathcal{K}}.
\label{5.4} \end{align}

If the vector fields $F^B_{\mathcal{K}}$ and $u^B_{\mathcal{K}}$ satisfy the following condition
\begin{equation}
 X^0_{\mathcal{K}}+ F^B_{\mathcal{K}}+u^B_{\mathcal{K}}=0, \;\; \label{5.5}
\end{equation}
then from (5.4) we have that $X^B_{(\mathcal{K},\omega^B_{\mathcal{K}},
H_{\mathcal{K}}, F^B_{\mathcal{K}}, u^B_{\mathcal{K}})}
=X_{\mathcal{K}}, $ that is, in this case the dynamical vector
field of the distributional CMH system is just the dynamical
vector field of the canonical distributional Hamiltonian system
without the actions of magnetic, external force and control.
Thus, the condition (5.5) is called the magnetic vanishing condition for
the distributional CMH system $(\mathcal{K},\omega^B_{\mathcal {K}},H_{\mathcal{K}},
F^B_{\mathcal{K}},u^B_{\mathcal{K}})$.\\

\section{Hamilton-Jacobi Equations for Distributional CMH System }

In this section we shall derive precisely
the geometric constraint conditions of the induced distributional two-form
for the dynamical vector field of distributional CMH system,
that is, the two types of Hamilton-Jacobi equation
for the distributional CMH system.
In order to do this, in the following we first give
some important notions and prove a key lemma, which is an important
tool for the proofs of two types of
Hamilton-Jacobi theorem for the distributional CMH system.\\

Assume that $\mathcal{D}\subset TQ$ is a $\mathcal{D}$-regularly nonholonomic
constraint, and the constraint submanifold
$\mathcal{M}=\mathcal{F}L(\mathcal{D})\subset T^*Q$,
the distribution
$\mathcal{F}=(T\pi_Q)^{-1}(\mathcal{D})\subset TT^*Q,$
where the projection $\pi_Q: T^* Q \rightarrow Q $ induces
the map $T\pi_{Q}: TT^* Q \rightarrow TQ. $
Assume that $\gamma: Q \rightarrow T^*Q$ is an one-form on $Q$, and
$B$ is a closed two-form on $Q$. But, note that
$\mathbf{d}\gamma $ is a two-form on $Q$, and
for any $v, w \in TT^* Q, $ we have that
$\mathbf{d}\gamma(T\pi_{Q}(v),T\pi_{Q}(w))=\pi^*(\mathbf{d}\gamma )(v, w)$
is a two-form on the cotangent bundle $T^*Q$, where
$\pi^*: T^*Q \rightarrow T^*T^*Q.$ Thus,
in the following we first introduce two weaker notions.

\begin{defi}
\noindent $(\mathrm{i})$ The one-form $\gamma$ is called to be closed
on $\mathcal{D}$ with respect to $T\pi_{Q}:
TT^* Q \rightarrow TQ, $ if for any $v, w \in TT^* Q, $
and $T\pi_{Q}(v), \; T\pi_{Q}(w) \in \mathcal{D},$  we have
that $\mathbf{d}\gamma(T\pi_{Q}(v),T\pi_{Q}(w))=0; $\\

\noindent $(\mathrm{ii})$ The one-form $\gamma: Q
\rightarrow T^*Q$ is called that satisfies condition that
$\mathbf{d}\gamma=-B$ on $\mathcal{D}$ with respect to $T\pi_{Q}:
TT^* Q \rightarrow TQ, $ if for any $v, w \in TT^* Q, $
and $T\pi_{Q}(v), \; T\pi_{Q}(w) \in \mathcal{D},$ we have
that  $(\mathbf{d}\gamma +B)(T\pi_{Q}(v),T\pi_{Q}(w))=0. $
\end{defi}

From the above Definition 6.1, we know that,
when $B=0$, the notion that, $\gamma$ satisfies condition
that $\mathbf{d}\gamma=-B$ on $\mathcal{D}$ with respect to $T\pi_{Q}:
TT^* Q \rightarrow TQ, $ become the notion that
$\gamma$ is closed on $\mathcal{D}$ with respect to $T\pi_{Q}:
TT^* Q \rightarrow TQ. $ On the other hand, it is worthy of noting that
the notion that $\gamma$ satisfies condition
that $\mathbf{d}\gamma=-B$ on $\mathcal{D}$ with respect to $T\pi_{Q}:
TT^* Q \rightarrow TQ, $ is weaker than the notion that $\gamma$
satisfies condition $\mathbf{d}\gamma=-B$ on $\mathcal{D},$
that is, $(\mathbf{d}\gamma+B)(x,y)=0, \; \forall\; x, y \in \mathcal{D}$.
In fact, if $\gamma$  satisfies condition
$\mathbf{d}\gamma=-B$ on $\mathcal{D}$,
then it must satisfy condition that $\mathbf{d}\gamma=-B$
on $\mathcal{D}$ with respect to $T\pi_{Q}: TT^* Q \rightarrow TQ. $
Conversely, if $\gamma$  satisfies condition
that $\mathbf{d}\gamma=-B$ on $\mathcal{D}$ with respect to
$T\pi_{Q}: TT^* Q \rightarrow TQ, $ then it may not satisfy condition
$\mathbf{d}\gamma=-B$ on $\mathcal{D}$.
We can prove a general result as follows, which states that
the notion that, the $\gamma$ satisfies condition that
$\mathbf{d}\gamma=-B$ on $\mathcal{D}$
with respect to $T\pi_{Q}: TT^* Q \rightarrow TQ, $
is not equivalent to the notion that $\gamma$ satisfies condition
$\mathbf{d}\gamma=-B$ on $\mathcal{D}$, also see Wang \cite{wa21c}.

\begin{prop}
Assume that $\gamma: Q \rightarrow T^*Q$ is an one-form on $Q$ and
it doesn't satisfy condition $\mathbf{d}\gamma=-B$ on $\mathcal{D}$.
We define the set $N$, which is a subset of $TQ$,
such that the one-form $\gamma$ on $N$ satisfies the condition that
for any $x,y \in N, \; (\mathbf{d}\gamma+B)(x,y)\neq 0. $ Denote
$Ker(T\pi_Q)= \{u \in TT^*Q| \; T\pi_Q(u)=0 \}, $ and $T\gamma: TQ
\rightarrow TT^* Q .$ If $T\gamma(N)\subset Ker(T\pi_Q), $ then
$\gamma$ satisfies condition $\mathbf{d}\gamma=-B$
with respect to $T\pi_{Q}: TT^* Q \rightarrow TQ.$
and hence $\gamma$ satisfies condition $\mathbf{d}\gamma=-B$
on $\mathcal{D}$ with respect to
$T\pi_{Q}: TT^* Q \rightarrow TQ.$
\end{prop}

\noindent{\bf Proof: } If the $\gamma: Q \rightarrow T^*Q$
doesn't satisfy condition $\mathbf{d}\gamma=-B$ on $\mathcal{D}$,
then it doesn't yet satisfy condition $\mathbf{d}\gamma=-B$.
For any $v, w \in TT^* Q, $ if
$T\pi_{Q}(v) \notin N, $ or $T\pi_{Q}(w))\notin N, $ then by the
definition of $N$, we know that
$(\mathbf{d}\gamma+B)(T\pi_{Q}(v),T\pi_{Q}(w))=0; $ If $T\pi_{Q}(v)\in
N, $ and $T\pi_{Q}(w))\in N, $ from the condition $T\gamma(N)\subset
Ker(T\pi_Q), $ we know that $T\pi_{Q}\cdot T\gamma \cdot
T\pi_{Q}(v)= T\pi_{Q}(v)=0, $ and $T\pi_{Q}\cdot T\gamma \cdot
T\pi_{Q}(w)= T\pi_{Q}(w)=0, $ where we have used the
relation $\pi_Q\cdot \gamma\cdot \pi_Q= \pi_Q, $ and hence
$(\mathbf{d}\gamma+B)(T\pi_{Q}(v),T\pi_{Q}(w))=0. $ Thus, for any $v, w
\in TT^* Q, $ we have always that
$(\mathbf{d}\gamma+B)(T\pi_{Q}(v),T\pi_{Q}(w))=0, $
In particular, for any $v, w \in TT^* Q, $
and $T\pi_{Q}(v), \; T\pi_{Q}(w) \in \mathcal{D},$  we have
$(\mathbf{d}\gamma+B)(T\pi_{Q}(v),T\pi_{Q}(w))=0. $
that is, $\gamma$ satisfies condition that $\mathbf{d}\gamma=-B$ on $\mathcal{D}$
with respect to $T\pi_{Q}: TT^* Q \rightarrow TQ. $
\hskip 0.3cm $\blacksquare$\\

Now, we prove the following Lemma 6.3. It is worthy of noting that
this lemma and Lemma 3.4 given in \S3 are the
important tool for the proofs of the two types of Hamilton-Jacobi
theorems for the distributional CMH system and the nonholonomic
reduced distributional CMH system.

\begin{lemm}
Assume that $\gamma: Q \rightarrow T^*Q$ is an one-form on $Q$, and
$\lambda=\gamma \cdot \pi_{Q}: T^* Q \rightarrow T^* Q ,$ and
$\omega$ is the canonical symplectic form on $T^*Q$, and
$\omega^B= \omega- \pi_Q^*B $
is the magnetic symplectic form on $T^*Q$.
If the Lagrangian $L$ is $\mathcal{D}$-regular, and
$\textmd{Im}(\gamma)\subset \mathcal{M}=\mathcal{F}L(\mathcal{D}), $
then we have that $ X^B_{H}\cdot \gamma \in \mathcal{F}$ along
$\gamma$, and $ X^B_{H}\cdot \lambda \in \mathcal{F}$ along
$\lambda$, that is, $T\pi_{Q}(X^B_H\cdot\gamma(q))\in
\mathcal{D}_{q}, \; \forall q \in Q $, and $T\pi_{Q}(X^B_H\cdot\lambda(q,p))\in
\mathcal{D}_{q}, \; \forall q \in Q, \; (q,p) \in T^* Q. $
Moreover, if a symplectic map $\varepsilon: T^* Q \rightarrow T^* Q $
with respect to the magnetic symplectic form $\omega^B$ satisfies the
condition $\varepsilon(\mathcal{M})\subset \mathcal{M},$ then
we have that $ X^B_{H}\cdot \varepsilon \in \mathcal{F}$ along
$\varepsilon. $
\end{lemm}

\noindent{\bf Proof:} Under the canonical cotangent bundle coordinates, for any $q \in Q
, \; (q,p)\in T^* Q, $ we have that
$$
X^B_H\cdot \gamma(q)= (\sum^n_{i=1} (\frac{\partial H}{\partial
p_i}\frac{\partial}{\partial q^i} - \frac{\partial H}{\partial
q^i}\frac{\partial}{\partial p_i})-
\sum^n_{i,j=1}B_{ij}\frac{\partial H}{\partial
p_j}\frac{\partial}{\partial p_i})\gamma(q).
$$
and
$$
X^B_H\cdot \lambda(q,p)= (\sum^n_{i=1} (\frac{\partial H}{\partial
p_i}\frac{\partial}{\partial q^i} - \frac{\partial H}{\partial
q^i}\frac{\partial}{\partial p_i})-
\sum^n_{i,j=1}B_{ij}\frac{\partial H}{\partial
p_j}\frac{\partial}{\partial p_i})\gamma\cdot \pi_Q(q,p).
$$
Then,
$$
T\pi_Q(X^B_H\cdot \gamma(q))=T\pi_Q(X^B_H\cdot \lambda(q,p))
=\sum^n_{i=1}(\frac{\partial H}{\partial
p_i}\frac{\partial}{\partial q^i})\gamma(q) \in T_q Q.
$$
Since $\textmd{Im}(\gamma)\subset \mathcal{M}, $ and
$\gamma(q)\in \mathcal{M}_{(q,p)}=\mathcal{F}L(\mathcal{D}_q), $ from the Lagrangian $L$ is
$\mathcal{D}$-regular, and $\mathcal{F}L$ is a diffeomorphism, then
there exists a point $(q, \; v_q)\in \mathcal{D}_q, $ such that
$\mathcal{F}L(q, \; v_q)=\gamma(q). $ Thus,
$$
T\pi_Q(X^B_H\cdot \gamma(q))=T\pi_Q(X^B_H\cdot
\lambda(q,p))=\mathcal{F}L(q, \; v_q)\sum^n_{i=1}(\frac{\partial H}{\partial
p_i}\frac{\partial}{\partial q^i}) \in
\mathcal{D}_q,
$$
it follows that $ X^B_{H}\cdot \gamma \in \mathcal{F}$ along
$\gamma$, and $ X^B_{H}\cdot \lambda \in \mathcal{F}$ along $\lambda$.
Moreover, for the symplectic map $\varepsilon: T^* Q \rightarrow T^* Q $
with respect to the magnetic symplectic form $\omega^B$, we have that
$$
X^B_H\cdot \varepsilon (q,p)= (\sum^n_{i=1} (\frac{\partial H}{\partial
p_i}\frac{\partial}{\partial q^i} - \frac{\partial H}{\partial
q^i}\frac{\partial}{\partial p_i})-
\sum^n_{i,j=1}B_{ij}\frac{\partial H}{\partial
p_j}\frac{\partial}{\partial p_i})\varepsilon (q,p).
$$
If $\varepsilon$ satisfies the
condition $\varepsilon(\mathcal{M})\subset \mathcal{M},$
then for any $(q,p)\in \mathcal{M}_{(q,p)}$, we have that
$\varepsilon(q,p)\in \mathcal{M}_{(q,p)},$
and there exists a point $(q,\; v_q)\in \mathcal{D}_q, $ such that
$\mathcal{F}L(q,\; v_q)=\varepsilon (q,p). $ Thus,
$$
T\pi_Q(X^B_H\cdot \varepsilon(q,p))
= \sum^n_{i=1}(\frac{\partial H}{\partial p_i}\frac{\partial}{\partial q^i})
\varepsilon(q, p)
=\mathcal{F}L(q, \; v_q)\sum^n_{i=1}(\frac{\partial H}{\partial
p_i}\frac{\partial}{\partial q^i}) \in
\mathcal{D}_q,
$$
it follows that $ X^B_{H}\cdot \varepsilon \in \mathcal{F}$ along
$\varepsilon$.
\hskip 0.3cm $\blacksquare$\\

We note that for a nonholonomic CMH system,
under the restriction given by nonholonomic constraint,
in general, the dynamical vector field of a nonholonomic CMH
system may not be Hamiltonian.
On the other hand, since the distributional CMH system
is determined by a non-degenerate distributional two-form
induced from the magnetic symplectic form, but, the non-degenerate distributional two-form
is not a "true two-form" on a manifold, and hence the leading
distributional CMH system can not be Hamiltonian.
Thus, we can not describe the Hamilton-Jacobi equations for the distributional
CMH system from the viewpoint of generating function
as in the classical Hamiltonian case, that is,
we cannot prove the Hamilton-Jacobi
theorem for the distributional CMH system,
just like same as the above Theorem 1.1.
In the following by using Lemma 3.4, Lemma 6.3,
and the non-degenerate
distributional two-form $\omega^B_{\mathcal{K}}$ and the
dynamical vector field $X^B_{(\mathcal{K},\omega^B_{\mathcal{K}},
H_{\mathcal{K}}, F^B_{\mathcal{K}}, u^B_{\mathcal{K}})}$ given
in \S 5 for the distributional CMH system,
we can derive precisely the geometric constraint conditions of
the non-degenerate distributional two-form $\omega^B_{\mathcal{K}}$
for the dynamical vector field $X^B_{(\mathcal{K},\omega^B_{\mathcal{K}},
H_{\mathcal{K}}, F^B_{\mathcal{K}}, u^B_{\mathcal{K}})}$,
that is, the two types of Hamilton-Jacobi equation for the distributional
CMH system $(\mathcal{K},\omega^B_{\mathcal {K}},H_{\mathcal {K}}, F^B_{\mathcal {K}}, u^B_{\mathcal {K}})$.
At first, we  prove the following Type I of
Hamilton-Jacobi theorem for the distributional CMH system.

\begin{theo} (Type I of Hamilton-Jacobi Theorem for the Distributional CMH System)
For the nonholonomic CMH system $(T^*Q,\omega^B,\mathcal{D},H,F,u)$
with an associated distributional CMH system
$(\mathcal{K},\omega^B_{\mathcal {K}},H_{\mathcal {K}}, F^B_{\mathcal {K}},u^B_{\mathcal {K}})$,
assume that $\gamma: Q \rightarrow T^*Q$ is an one-form on $Q$,
and $\tilde{X}^\gamma = T\pi_{Q}\cdot \tilde{X}\cdot \gamma$,
where $\tilde{X}=X_{(\mathcal{K},\omega^B_{\mathcal{K}},
H_{\mathcal{K}}, F^B_{\mathcal{K}}, u^B_{\mathcal{K}})}
=X^B_\mathcal {K}+ F^B_{\mathcal{K}}+u^B_{\mathcal{K}}$
is the dynamical vector field of the distributional CMH system.
Moreover, assume that $\textmd{Im}(\gamma)\subset
\mathcal{M}=\mathcal{F}L(\mathcal{D}), $ and $
\textmd{Im}(T\gamma)\subset \mathcal{K}. $ If the
one-form $\gamma: Q \rightarrow T^*Q $ satisfies the condition,
$\mathbf{d}\gamma=-B $ on $\mathcal{D}$ with respect to
$T\pi_Q: TT^* Q \rightarrow TQ, $ then $\gamma$ is a
solution of the equation $T\gamma \cdot
\tilde{X}^\gamma= X^B_{\mathcal{K}} \cdot \gamma. $ Here
$X^B_{\mathcal{K}}$ is the nonholonomic dynamical vector field
of the distributional CMH system
$(\mathcal{K},\omega^B_{\mathcal {K}},
H_{\mathcal {K}}, F^B_{\mathcal {K}}, u^B_{\mathcal {K}})$.
The equation $T\gamma \cdot \tilde{X}^\gamma
= X^B_{\mathcal{K}} \cdot \gamma $ is called the Type I of
Hamilton-Jacobi equation for the distributional CMH system
$(\mathcal{K},\omega^B_{\mathcal {K}},H_{\mathcal {K}}, F^B_{\mathcal {K}}, u^B_{\mathcal {K}})$.
Here the maps involved in the theorem are shown in
the following Diagram-5.
\begin{center}
\hskip 0cm \xymatrix{& \mathcal{M} \ar[d]_{X^B_{\mathcal{K}}}
\ar[r]^{i_{\mathcal{M}}} & T^* Q \ar[d]_{X^B_{H}}
 \ar[r]^{\pi_Q}
& Q \ar[d]_{{\tilde{X}}^{\gamma}} \ar[r]^{\gamma} & T^*Q \ar[d]^{\tilde{X}} \\
& \mathcal{K}  & T(T^*Q) \ar[l]_{\tau_{\mathcal{K}}} & TQ
\ar[l]_{T\gamma} & T(T^* Q)\ar[l]_{T\pi_Q}}
\end{center}
$$\mbox{Diagram-5}$$
\end{theo}

\noindent{\bf Proof: } From Definition 5.1 we have that
$\tilde{X}=X^B_{(\mathcal{K},\omega^B_{\mathcal{K}},
H_{\mathcal{K}}, F^B_{\mathcal{K}}, u^B_{\mathcal{K}})}
=X^B_\mathcal {K}+ F^B_{\mathcal{K}}+u^B_{\mathcal{K}}$,
and $F^B_{\mathcal{K}}=\tau_{\mathcal{K}}\cdot \textnormal{vlift}(F_{\mathcal{M}})X^B_H$,
and $u^B_{\mathcal{K}}=\tau_{\mathcal{K}}\cdot \textnormal{vlift}(u_{\mathcal{M}})X^B_H$,
note that $T\pi_{Q}\cdot \textnormal{vlift}(F_{\mathcal{M}})X^B_H=T\pi_{Q}\cdot \textnormal{vlift}(u_{\mathcal{M}})X^B_H=0, $
then we have that $T\pi_{Q}\cdot F^B_{\mathcal{K}}=T\pi_{Q}\cdot u^B_{\mathcal{K}}=0,$
and hence $T\pi_{Q}\cdot \tilde{X}\cdot \gamma=T\pi_{Q}\cdot X^B_{\mathcal{K}}\cdot \gamma. $
On the other hand, we note that
$\textmd{Im}(\gamma)\subset \mathcal{M}, $ and
$\textmd{Im}(T\gamma)\subset \mathcal{K}, $ in this case,
$\omega^B_{\mathcal{K}}\cdot
\tau_{\mathcal{K}}=\tau_{\mathcal{K}}\cdot \omega^B_{\mathcal{M}}=
\tau_{\mathcal{K}}\cdot i_{\mathcal{M}}^* \cdot \omega^B, $ along
$\textmd{Im}(T\gamma)$. Moreover, from the distributional magnetic Hamiltonian equation (5.2),
we have that $X^B_{\mathcal{K}}= \tau_{\mathcal{K}}\cdot X^B_H,$
and $\tau_{\mathcal{K}}\cdot X^B_{H}\cdot \gamma = X^B_{\mathcal{K}}\cdot \gamma $.
Thus, using the non-degenerate
distributional two-form $\omega^B_{\mathcal{K}}$, from Lemma 3.4(ii) and Lemma 6.3,
if we take that $v= X^B_{\mathcal{K}}\cdot \gamma \in \mathcal{K} (\subset \mathcal{F}),$
and for any $w \in \mathcal{F}, \; T\lambda(w)\neq 0, $ and
$\tau_{\mathcal{K}}\cdot w \neq 0, $ then we have that
\begin{align*}
& \omega^B_{\mathcal{K}}(T\gamma \cdot \tilde{X}^\gamma, \;
\tau_{\mathcal{K}}\cdot w)=
\omega^B_{\mathcal{K}}(\tau_{\mathcal{K}}\cdot T\gamma \cdot
T\pi_{Q}\cdot \tilde{X}\cdot \gamma, \; \tau_{\mathcal{K}}\cdot w)\\ & =
\tau_{\mathcal{K}}\cdot i_{\mathcal{M}}^* \cdot \omega^B(T\gamma \cdot
T\pi_{Q}\cdot X^B_{\mathcal{K}}\cdot \gamma, \; w ) = \tau_{\mathcal{K}}\cdot
i_{\mathcal{M}}^* \cdot \omega^B (T(\gamma \cdot \pi_Q)\cdot X^B_{\mathcal{K}} \cdot \gamma, \; w)\\
& =\tau_{\mathcal{K}}\cdot i_{\mathcal{M}}^* \cdot
(\omega^B (X^B_{\mathcal{K}} \cdot \gamma, \; w-T(\gamma \cdot \pi_Q)\cdot w)
-(\mathbf{d}\gamma+B)(T\pi_{Q}\cdot X^B_{\mathcal{K}}\cdot\gamma, \; T\pi_{Q}\cdot w))\\
& = \tau_{\mathcal{K}}\cdot i_{\mathcal{M}}^* \cdot \omega^B (X^B_{\mathcal{K}} \cdot
\gamma, \; w) - \tau_{\mathcal{K}}\cdot i_{\mathcal{M}}^* \cdot
\omega^B (X^B_{\mathcal{K}} \cdot \gamma, \; T(\gamma
\cdot \pi_Q) \cdot w) \\
& \;\;\;\;\;\; - \tau_{\mathcal{K}}\cdot i_{\mathcal{M}}^* \cdot
(\mathbf{d}\gamma+B)(T\pi_{Q}\cdot X^B_{\mathcal{K}}\cdot\gamma, \; T\pi_{Q}\cdot w)\\
& = \omega^B_{\mathcal{K}}( \tau_{\mathcal{K}}\cdot X^B_{\mathcal{K}}\cdot \gamma,
\; \tau_{\mathcal{K}}\cdot w) -
\omega^B_{\mathcal{K}}(\tau_{\mathcal{K}}\cdot X^B_{\mathcal{K}} \cdot \gamma, \;
\tau_{\mathcal{K}}\cdot T(\gamma \cdot \pi_Q) \cdot w)\\
& \;\;\;\;\;\; - \tau_{\mathcal{K}}\cdot i_{\mathcal{M}}^* \cdot
(\mathbf{d}\gamma+B)(T\pi_{Q}\cdot X^B_{\mathcal{K}}\cdot\gamma, \; T\pi_{Q}\cdot w)\\
& = \omega^B_{\mathcal{K}}(X^B_{\mathcal{K}}\cdot \gamma, \;
\tau_{\mathcal{K}} \cdot w) -
\omega^B_{\mathcal{K}}(X^B_{\mathcal{K}} \cdot \gamma, \;
 T\gamma \cdot T\pi_{Q}(w))\\
& \;\;\;\;\;\; - \tau_{\mathcal{K}}\cdot i_{\mathcal{M}}^* \cdot
(\mathbf{d}\gamma+B)(T\pi_{Q}\cdot X^B_{\mathcal{K}}\cdot\gamma, \; T\pi_{Q}\cdot w),
\end{align*}
where we have used that $ \tau_{\mathcal{K}}\cdot T\gamma= T\gamma, $ and
$\tau_{\mathcal{K}}\cdot X^B_{\mathcal{K}}\cdot \gamma = X^B_{\mathcal{K}}\cdot
\gamma, $ since $\textmd{Im}(T\gamma)\subset \mathcal{K}. $
Note that $X^B_{\mathcal{K}}\cdot \gamma, \; w \in \mathcal{F},$ and
$T\pi_{Q}(X^B_{\mathcal{K}}\cdot \gamma), \; T\pi_{Q}(w) \in \mathcal{D}. $
If the one-form $\gamma: Q \rightarrow T^*Q $ satisfies the condition,
$\mathbf{d}\gamma=-B $ on $\mathcal{D}$ with respect to
$T\pi_Q: TT^* Q \rightarrow TQ, $ then
$(\mathbf{d}\gamma+B)(T\pi_{Q}\cdot X^B_{\mathcal{K}}\cdot\gamma, \; T\pi_{Q}\cdot w)=0,$
and hence
$$
\tau_{\mathcal{K}}\cdot i_{\mathcal{M}}^* \cdot(\mathbf{d}\gamma+B)
(T\pi_{Q}(X^B_{\mathcal{K}}\cdot \gamma), \; T\pi_{Q}(w))=0,
$$
Thus, we have that
\begin{equation}
\omega^B_{\mathcal{K}}(T\gamma \cdot \tilde{X}^\gamma, \;
\tau_{\mathcal{K}}\cdot w)- \omega^B_{\mathcal{K}}(X^B_{\mathcal{K}}\cdot \gamma, \;
\tau_{\mathcal{K}} \cdot w)
= -\omega^B_{\mathcal{K}}(X^B_{\mathcal{K}}\cdot
\gamma, \; T\gamma \cdot T\pi_{Q}(w)).
\label{6.1} \end{equation}
If $\gamma$ satisfies the equation $T\gamma \cdot \tilde{X}^\gamma= X^B_{\mathcal{K}}\cdot \gamma ,$
from Lemma 3.4(i) we know that the right side of (6.1) becomes that
\begin{align*}
 -\omega^B_{\mathcal{K}}(X^B_{\mathcal{K}} \cdot \gamma, \;
T\gamma \cdot T\pi_{Q}(w))
& = -\omega^B_{\mathcal{K}}(T\gamma\cdot \tilde{X}^\gamma, \;
T\gamma \cdot T\pi_{Q}(w))\\
& = -\omega^B_{\mathcal{K}}(\tau_{\mathcal{K}}\cdot T\gamma\cdot T\pi_{Q}\cdot
\tilde{X}\cdot \gamma, \; \tau_{\mathcal{K}}\cdot T\gamma \cdot T\pi_{Q}(w))\\
& = -\tau_{\mathcal{K}}\cdot
i_{\mathcal{M}}^* \cdot \omega^B(T\gamma
\cdot T\pi_{Q}(X^B_{\mathcal{K}}\cdot\gamma), \; T\gamma \cdot T\pi_{Q}(w))\\
& = -\tau_{\mathcal{K}}\cdot
i_{\mathcal{M}}^* \cdot \lambda^* \omega^B (X^B_{\mathcal{K}}\cdot\gamma, \; w)\\
& = \tau_{\mathcal{K}}\cdot
i_{\mathcal{M}}^* \cdot (\mathbf{d}\gamma+B)
(T\pi_{Q}\cdot X^B_{\mathcal{K}}\cdot\gamma, \; T\pi_{Q}\cdot w)=0.
\end{align*}
Because the distributional two-form $\omega^B_{\mathcal{K}}$ is non-degenerate,
the left side of (6.1) equals zero, only when
$\gamma$ satisfies the equation $T\gamma\cdot \tilde{X}^\gamma= X^B_{\mathcal{K}}\cdot \gamma .$ Thus,
if the one-form $\gamma: Q \rightarrow T^*Q $ satisfies the condition that
$\mathbf{d}\gamma=-B $ on $\mathcal{D}$ with respect to
$T\pi_Q: TT^* Q \rightarrow TQ, $ then $\gamma$ must be a solution of the Type I of Hamilton-Jacobi equation
$T\gamma \cdot \tilde{X}^\gamma= X^B_{\mathcal{K}}\cdot \gamma ,$
for the distributional CMH system
$(\mathcal{K},\omega^B_{\mathcal {K}},H_{\mathcal {K}}, F^B_{\mathcal {K}}, u^B_{\mathcal {K}})$.
\hskip 0.3cm $\blacksquare$\\

It is worthy of noting that, when $B=0$, in this case the magnetic symplectic form $\omega^B$
is just the canonical symplectic form $\omega$ on $T^*Q$,
and the nonholonomic CMH system
$(T^*Q,\omega^B,\mathcal{D},H, F, u)$
becomes the nonholonomic RCH system $(T^*Q,\omega,\mathcal{D},H, F, u)$
with the canonical symplectic form $\omega$,
and the distributional CMH system
$(\mathcal{K},\omega^B_{\mathcal {K}},H_{\mathcal {K}}, F^B_{\mathcal {K}}, u^B_{\mathcal {K}})$
becomes the distributional RCH system
$(\mathcal{K},\omega_{\mathcal {K}},H_{\mathcal {K}}, F_{\mathcal {K}}, u_{\mathcal {K}})$,
and the condition that the one-form $\gamma: Q \rightarrow T^*Q $ satisfies the condition that
$\mathbf{d}\gamma=-B $ on $\mathcal{D}$ with respect to $T\pi_{Q}:
TT^* Q \rightarrow TQ, $ becomes that  $\gamma $ is closed on $\mathcal{D}$ with respect to
$T\pi_Q: TT^* Q \rightarrow TQ.$ Thus, from above Theorem 6.4,
we can obtain Theorem 3.4 given in Wang \cite{wa22a}, that is.
the Type I of Hamilton-Jacobi theorem for the distributional RCH system.
On the other hand, from the proofs of Theorem 3.4 given
in Wang \cite{wa22a}, we know that, if the one-form
$\gamma: Q \rightarrow T^*Q $ is not closed on $\mathcal{D}$ with respect to
$T\pi_Q: TT^* Q \rightarrow TQ, $ then $\gamma$ is not a solution of the Type I
of Hamilton-Jacobi equation for the distributional RCH system. In this case,
our idea is that we hope to look for a nonholonomic CMH system,
such that $\gamma$ is a solution of the Type I of Hamilton-Jacobi
equation for the associated distributional CMH system. Since
$\gamma: Q \rightarrow T^*Q $ is not closed on $\mathcal{D}$ with respect to
$T\pi_Q: TT^* Q \rightarrow TQ, $ then
the $\gamma$ is not yet closed on $\mathcal{D}$,
that is,  $\mathbf{d}\gamma(x,y)\neq 0, \; \forall\;
x, y \in \mathcal{D}$, and hence $\gamma$ is not yet closed on $Q$.
However, in this case, we note that
$\mathbf{d}\cdot \mathbf{d}\gamma= \mathbf{d}^2 \gamma =0, $
and hence the $\mathbf{d}\gamma$ is a closed two-form on $Q$. Thus, we can construct
a magnetic symplectic form on $T^*Q$, $\omega^B= \omega - \pi_Q^* B
= \omega+ \pi_Q^*(\mathbf{d}\gamma), $
where $B=- \mathbf{d}\gamma,$.
Moreover,  we can also construct a nonholonomic CMH system
$(T^*Q, \omega^B,\mathcal{D}, H, F, u )$
with an associated distributional CMH system
$(\mathcal{K},\omega^B_{\mathcal{K}},H_{\mathcal {K}}, F^B_{\mathcal {K}}, u^B_{\mathcal {K}})$,
which satisfies the distributional magnetic Hamiltonian equation (5.2),
$\mathbf{i}_{X^B_{\mathcal{K}} }\omega^B_{\mathcal{K}}= \mathbf{d}H_{\mathcal{K}},$
and $F^B_{\mathcal{K}}=\tau_{\mathcal{K}}\cdot \textnormal{vlift}(F_{\mathcal{M}})X^B_H,$
and $u^B_{\mathcal{K}}=\tau_{\mathcal{K}}\cdot \textnormal{vlift}(u_{\mathcal{M}})X^B_H$.
In this case, the one-form $\gamma: Q \rightarrow T^*Q $ satisfies the condition that
$\mathbf{d}\gamma=-B$ with respect to $T\pi_{Q}: TT^* Q \rightarrow TQ, $
and hence it satisfies also the condition that $\mathbf{d}\gamma=-B$
on $\mathcal{D}$ with respect to $T\pi_{Q}: TT^* Q \rightarrow TQ. $
By using Lemma 3.4, Lemma 6.3, and the non-degenerate
distributional two-form $\omega^B_{\mathcal{K}}$ and the
dynamical vector field $X^B_{(\mathcal{K},\omega^B_{\mathcal{K}},
H_{\mathcal{K}}, F^B_{\mathcal{K}}, u^B_{\mathcal{K}})}$,
from Theorem 6.4 we can obtain the following Theorem 6.5.
\begin{theo}
For a given nonholonomic RCH system $(T^*Q,\omega,\mathcal{D},H, F, u)$ with
the canonical symplectic form $\omega$ on $T^*Q$ and
$\mathcal{D}$-completely and $\mathcal{D}$-regularly nonholonomic
constraint $\mathcal{D} \subset TQ$,
and assume that the one-form $\gamma: Q
\rightarrow T^*Q$ is not closed on $\mathcal{D}$ with respect to
$T\pi_Q: TT^* Q \rightarrow TQ. $ Then one can construct
a magnetic symplectic form on $T^*Q$,
$\omega^B= \omega+ \pi_Q^*(\mathbf{d}\gamma), $ where $B=- \mathbf{d}\gamma,$
and a nonholonomic CMH system $(T^*Q, \omega^B, \mathcal{D},H, F, u )$
with an associated distributional CMH system
$(\mathcal{K},\omega^B_{\mathcal{K}},H_{\mathcal {K}}, F^B_{\mathcal {K}}, u^B_{\mathcal {K}})$.
Denote $\tilde{X}^\gamma = T\pi_{Q}\cdot \tilde{X}\cdot \gamma$,
where $\tilde{X}=X_{(\mathcal{K},\omega^B_{\mathcal{K}},
H_{\mathcal{K}}, F^B_{\mathcal{K}}, u^B_{\mathcal{K}})}
=X^B_\mathcal {K}+ F^B_{\mathcal{K}}+u^B_{\mathcal{K}}$
is the dynamical vector field of the distributional CMH system.
Moreover, assume that $\textmd{Im}(\gamma)\subset
\mathcal{M}=\mathcal{F}L(\mathcal{D}), $ and $
\textmd{Im}(T\gamma)\subset \mathcal{K}. $
Then $\gamma$ is a solution of the Type I of
Hamilton-Jacobi equation $T\gamma \cdot \tilde{X}^\gamma= X^B_{\mathcal{K}}\cdot \gamma ,$
for the distributional CMH system
$(\mathcal{K},\omega^B_{\mathcal {K}},H_{\mathcal {K}}, F^B_{\mathcal {K}}, u^B_{\mathcal {K}})$.
\end{theo}

Next, for any symplectic map $\varepsilon: T^* Q \rightarrow T^* Q $
with respect to the magnetic symplectic form $\omega^B$,
we can prove the following Type II of
Hamilton-Jacobi theorem for the distributional CMH system.
For convenience, the
maps involved in the following theorem and its proof are shown in
Diagram-6.

\begin{center}
\hskip 0cm \xymatrix{& \mathcal{M} \ar[d]_{X^B_{\mathcal{K}}}
\ar[r]^{i_{\mathcal{M}}} & T^* Q \ar[d]_{X^B_{H\cdot \varepsilon}}
\ar[dr]^{X^\varepsilon} \ar[r]^{\pi_Q}
& Q \ar[r]^{\gamma} & T^*Q \ar[d]^{X^B_H} \\
& \mathcal{K}  & T(T^*Q) \ar[l]_{\tau_{\mathcal{K}}} & TQ
\ar[l]_{T\gamma} & T(T^* Q)\ar[l]_{T\pi_Q}}
\end{center}
$$\mbox{Diagram-6}$$

\begin{theo} (Type II of Hamilton-Jacobi Theorem for a Distributional CMH System)
For the nonholonomic CMH system $(T^*Q,\omega^B,\mathcal{D},H,F,u)$
with an associated distributional CMH system
$(\mathcal{K},\omega^B_{\mathcal {K}},H_{\mathcal {K}}, F^B_{\mathcal {K}}, u^B_{\mathcal {K}})$,
assume that $\gamma: Q \rightarrow T^*Q$ is an one-form on $Q$,
and $\lambda=\gamma \cdot \pi_{Q}: T^* Q \rightarrow T^* Q, $ and for any
symplectic map $\varepsilon: T^* Q \rightarrow T^* Q $ with respect to
the magnetic symplectic form $\omega^B$, denote by
$\tilde{X}^\varepsilon = T\pi_{Q}\cdot \tilde{X}\cdot \varepsilon$,
where $\tilde{X}=X_{(\mathcal{K},\omega^B_{\mathcal{K}},
H_{\mathcal{K}}, F^B_{\mathcal{K}}, u^B_{\mathcal{K}})}
=X^B_{\mathcal {K}}+ F^B_{\mathcal{K}}+u^B_{\mathcal{K}}$
is the dynamical vector field of the distributional CMH system.
Moreover, assume that $\textmd{Im}(\gamma)\subset
\mathcal{M}=\mathcal{F}L(\mathcal{D}), $ and $\varepsilon(\mathcal{M})\subset \mathcal{M}$,
and $\textmd{Im}(T\gamma)\subset \mathcal{K}. $
Then $\varepsilon$ is a solution of the equation
$\tau_{\mathcal{K}}\cdot T\varepsilon(X^B_{H\cdot\varepsilon})
= T\lambda \cdot \tilde{X}\cdot\varepsilon,$
if and only if it is a solution of the equation
$T\gamma \cdot \tilde{X}^\varepsilon= X^B_{\mathcal{K}} \cdot \varepsilon $.
Here $ X^B_{H\cdot\varepsilon}$ is the magnetic Hamiltonian vector field of the function
$H \cdot \varepsilon: T^* Q\rightarrow \mathbb{R}, $ and $X^B_{\mathcal{K}}$
is the nonholonomic dynamical vector field of the distributional CMH system
$(\mathcal{K},\omega^B_{\mathcal {K}},H_{\mathcal {K}}, F^B_{\mathcal {K}}, u^B_{\mathcal {K}})$.
The equation $T\gamma \cdot \tilde{X}^\varepsilon
= X^B_{\mathcal{K}} \cdot \varepsilon,$ is called the Type II of
Hamilton-Jacobi equation for the distributional CMH system
$(\mathcal{K},\omega^B_{\mathcal {K}},H_{\mathcal {K}}, F^B_{\mathcal {K}}, u^B_{\mathcal {K}})$.
\end{theo}

\noindent{\bf Proof: } In the same way,
from Definition 5.1 we have that
$\tilde{X}=X^B_{(\mathcal{K},\omega^B_{\mathcal{K}},
H_{\mathcal{K}}, F^B_{\mathcal{K}}, u^B_{\mathcal{K}})}
=X^B_\mathcal {K}+ F^B_{\mathcal{K}}+u^B_{\mathcal{K}}$,
and $F^B_{\mathcal{K}}=\tau_{\mathcal{K}}\cdot \textnormal{vlift}(F_{\mathcal{M}})X^B_H$,
and $u^B_{\mathcal{K}}=\tau_{\mathcal{K}}\cdot \textnormal{vlift}(u_{\mathcal{M}})X^B_H$,
note that $T\pi_{Q}\cdot \textnormal{vlift}(F_{\mathcal{M}})X^B_H=T\pi_{Q}\cdot \textnormal{vlift}(u_{\mathcal{M}})X^B_H=0, $
then we have that $T\pi_{Q}\cdot F^B_{\mathcal{K}}=T\pi_{Q}\cdot u^B_{\mathcal{K}}=0,$
and hence $T\pi_{Q}\cdot \tilde{X}\cdot \varepsilon =T\pi_{Q}\cdot X^B_{\mathcal{K}}\cdot \varepsilon. $
On the other hand, we note that
$\textmd{Im}(\gamma)\subset \mathcal{M}, $ and
$\textmd{Im}(T\gamma)\subset \mathcal{K}, $ in this case,
$\omega^B_{\mathcal{K}}\cdot
\tau_{\mathcal{K}}=\tau_{\mathcal{K}}\cdot \omega^B_{\mathcal{M}}=
\tau_{\mathcal{K}}\cdot i_{\mathcal{M}}^* \cdot \omega^B, $ along
$\textmd{Im}(T\gamma)$. Moreover, from the distributional magnetic Hamiltonian equation (5.2),
we have that $X^B_{\mathcal{K}}= \tau_{\mathcal{K}}\cdot X^B_H,$
and $\tau_{\mathcal{K}}\cdot X^B_{H}\cdot \varepsilon = X^B_{\mathcal{K}}\cdot \varepsilon $.
Note that $\varepsilon(\mathcal{M})\subset \mathcal{M},$ and
$T\pi_{Q}(X^B_H\cdot \varepsilon(q,p))\in
\mathcal{D}_{q}, \; \forall q \in Q, \; (q,p) \in \mathcal{M}(\subset T^* Q), $
and hence $X^B_H\cdot \varepsilon \in \mathcal{F}$ along $\varepsilon$,
and $X^B_{\mathcal{K}}\cdot \varepsilon
=\tau_{\mathcal{K}}\cdot X^B_{H}\cdot \varepsilon
\in \mathcal{K} (\subset \mathcal{F})$.
Thus, using the non-degenerate
distributional two-form $\omega^B_{\mathcal{K}}$, from Lemma 3.4 and Lemma 6.3, if
we take that $v= X^B_{\mathcal{K}}\cdot
\varepsilon \in \mathcal{K} (\subset \mathcal{F}), $ and for any $w
\in \mathcal{F}, \; T\lambda(w)\neq 0, $ and
$\tau_{\mathcal{K}}\cdot w \neq 0, $ then we have that
\begin{align*}
& \omega^B_{\mathcal{K}}(T\gamma \cdot \tilde{X}^\varepsilon, \;
\tau_{\mathcal{K}}\cdot w)=
\omega^B_{\mathcal{K}}(\tau_{\mathcal{K}}\cdot T\gamma \cdot
\tilde{X}^\varepsilon, \; \tau_{\mathcal{K}}\cdot w)\\ & =
\tau_{\mathcal{K}}\cdot i_{\mathcal{M}}^* \cdot\omega^B(T\gamma \cdot
T\pi_{Q}\cdot \tilde{X}\cdot \varepsilon, \; w ) = \tau_{\mathcal{K}}\cdot
i_{\mathcal{M}}^* \cdot\omega^B(T(\gamma \cdot \pi_Q)\cdot X^B_{\mathcal{K}}\cdot \varepsilon, \; w)\\
& =\tau_{\mathcal{K}}\cdot i_{\mathcal{M}}^* \cdot(\omega^B(X^B_{\mathcal{K}}\cdot
\varepsilon, \; w-T(\gamma \cdot \pi_Q)\cdot w)
-(\mathbf{d}\gamma+B)(T\pi_{Q}(X^B_{\mathcal{K}}\cdot \varepsilon), \; T\pi_{Q}(w)))\\
& = \tau_{\mathcal{K}}\cdot i_{\mathcal{M}}^* \cdot\omega^B(X^B_{\mathcal{K}}\cdot
\varepsilon, \; w) - \tau_{\mathcal{K}}\cdot i_{\mathcal{M}}^* \cdot
\omega^B(X^B_{\mathcal{K}}\cdot \varepsilon, \; T\lambda \cdot w)\\
& \;\;\;\;\;\;
-\tau_{\mathcal{K}}\cdot i_{\mathcal{M}}^* \cdot (\mathbf{d}\gamma+B)
(T\pi_{Q}(X^B_{\mathcal{K}}\cdot \varepsilon), \; T\pi_{Q}(w))\\
& = \omega^B_{\mathcal{K}}( \tau_{\mathcal{K}}\cdot X^B_{\mathcal{K}}\cdot \varepsilon,
\; \tau_{\mathcal{K}}\cdot w) -
\omega^B_{\mathcal{K}}(\tau_{\mathcal{K}}\cdot X^B_{\mathcal{K}}\cdot \varepsilon, \;
\tau_{\mathcal{K}}\cdot T\lambda \cdot w)\\
& \;\;\;\;\;\;
+\tau_{\mathcal{K}}\cdot i_{\mathcal{M}}^* \cdot \lambda^* \omega^B
(X^B_{\mathcal{K}}\cdot \varepsilon, \; w)\\
& = \omega^B_{\mathcal{K}}(X^B_{\mathcal{K}}\cdot \varepsilon, \;
\tau_{\mathcal{K}} \cdot w) -
\omega^B_{\mathcal{K}}(X^B_{\mathcal{K}}\cdot \varepsilon,
\; T\lambda \cdot w)+ \omega^B_{\mathcal{K}}(T\lambda\cdot X^B_{\mathcal{K}}\cdot \varepsilon,
\; T\lambda \cdot w),
\end{align*}
where we have used that $ \tau_{\mathcal{K}}\cdot T\gamma
= T\gamma, \; \tau_{\mathcal{K}}\cdot T\lambda= T\lambda, $ and
$\tau_{\mathcal{K}}\cdot X^B_{\mathcal{K}}\cdot \varepsilon = X^B_{\mathcal{K}}\cdot
\varepsilon, $ since $\textmd{Im}(T\gamma)\subset \mathcal{K}. $
From the distributional magnetic Hamiltonian equation (5.2),
$\mathbf{i}_{X^B_\mathcal{K}}\omega^B_{\mathcal{K}}=\mathbf{d}H_\mathcal
{K}$, we have that $X^B_{\mathcal{K}}=\tau_{\mathcal{K}}\cdot X^B_H. $
Note that $\varepsilon: T^* Q \rightarrow T^* Q $ is symplectic
with respect to the magnetic symplectic form $\omega^B$, and $
X^B_H\cdot \varepsilon = T\varepsilon \cdot X^B_{H\cdot\varepsilon}, $ along
$\varepsilon$, and hence $X^B_{\mathcal{K}}\cdot \varepsilon
=\tau_{\mathcal{K}}\cdot X^B_H \cdot \varepsilon=
\tau_{\mathcal{K}}\cdot T\varepsilon \cdot X^B_{H \cdot \varepsilon}, $ along $\varepsilon$.
Note that
$T\lambda \cdot X^B_\mathcal{K}\cdot \varepsilon=T\gamma \cdot
T\pi_Q\cdot X^B_\mathcal{K}\cdot \varepsilon=T\gamma \cdot
T\pi_Q\cdot \tilde{X}\cdot \varepsilon=T\lambda\cdot \tilde{X}\cdot \varepsilon.$
Then we have that
\begin{align*}
& \omega^B_{\mathcal{K}}(T\gamma \cdot \tilde{X}^\varepsilon, \;
\tau_{\mathcal{K}}\cdot w)-
\omega^B_{\mathcal{K}}(X^B_{\mathcal{K}}\cdot \varepsilon, \;
\tau_{\mathcal{K}} \cdot w) \nonumber \\
& = -\omega^B_{\mathcal{K}}(X^B_{\mathcal{K}}\cdot \varepsilon,
\; T\lambda \cdot w)+ \omega^B_{\mathcal{K}}
(T\lambda\cdot X^B_{\mathcal{K}}\cdot \varepsilon, \; T\lambda \cdot w) \\
& = - \omega^B_{\mathcal{K}}(\tau_{\mathcal{K}}\cdot X^B_H \cdot \varepsilon, \;
 T\lambda \cdot w)+ \omega^B_{\mathcal{K}}(T\lambda\cdot \tilde{X}\cdot \varepsilon,
\; T\lambda \cdot w)\\
&= \omega^B_{\mathcal{K}}(T\lambda\cdot \tilde{X}\cdot \varepsilon
-\tau_{\mathcal{K}}\cdot T\varepsilon \cdot X^B_{H \cdot \varepsilon},
\; T\lambda \cdot w).
\end{align*}
Because the induced distributional two-form
$\omega^B_{\mathcal{K}}$ is non-degenerate, it follows that the equation
$T\gamma\cdot \tilde{X}^\varepsilon= X^B_{\mathcal{K}}\cdot
\varepsilon ,$ is equivalent to the equation
$\tau_{\mathcal{K}}\cdot T\varepsilon \cdot X^B_{H\cdot\varepsilon}
= T\lambda\cdot \tilde{X}\cdot \varepsilon $.
Thus, $\varepsilon$ is a solution of the equation
$\tau_{\mathcal{K}}\cdot T\varepsilon\cdot X^B_{H\cdot\varepsilon}
= T\lambda \cdot \tilde{X} \cdot\varepsilon,$
if and only if it is a solution of
the Type II of Hamilton-Jacobi equation
$T\gamma\cdot \tilde{X}^\varepsilon= X^B_{\mathcal{K}}\cdot
\varepsilon .$
\hskip 0.3cm $\blacksquare$

\begin{rema}
It is worthy of noting that, the Type I of Hamilton-Jacobi equation
$T\gamma \cdot \tilde{X}^\gamma= X^B_{\mathcal{K}}\cdot \gamma ,$
is the equation of the differential one-form $\gamma$; and
the Type II of Hamilton-Jacobi equation $T\gamma \cdot \tilde{X}^\varepsilon
= X^B_{\mathcal{K}}\cdot \varepsilon ,$ is the equation of
the symplectic diffeomorphism map $\varepsilon$.
If the nonholonomic CMH system we considered has not any constrains,
in this case, the distributional CMH system is just the CMH system itself.
From the above Type I and Type II of Hamilton-Jacobi theorems, that is,
Theorem 6.4- 6.6, we can get the Theorem 3.5- 3.7.
It shows that Theorem 6.4- 6.6 can be regarded as an extension of two types of
Hamilton-Jacobi theorem for the CMH system to the system with
nonholonomic context. On the other hand,
if both the external force and control of a nonholonomic CMH
system $(T^*Q,\omega^B,\mathcal{D},H,F,W)$ are zero, that is, $F=0 $
and $W=\emptyset$, in this case the nonholonomic CMH system
is just a nonholonomic magnetic Hamiltonian system $(T^*Q,\omega^B,\mathcal{D},H)$,
and from the proofs of
the above Theorem 6.4- 6.6 , we can obtain two types of Hamilton-Jacobi
equation for the associated distributional magnetic Hamiltonian system,
that is, Theorem 4.4- 4.6 given in Wang \cite{wa21c}.
Thus, Theorem 6.4- 6.6 can be regarded as an extension of two types of Hamilton-Jacobi
equation for a nonholonomic magnetic Hamiltonian system to
that for the system with external force and control.
In particular, in this case, if $B=0$, then
the magnetic symplectic form $\omega^B$
is just the canonical symplectic form $\omega$ on $T^*Q$, and
the distributional magnetic Hamiltonian system is just the distributional Hamiltonian system itself.
From the above Type I and Type II of Hamilton-Jacobi theorems, that is,
Theorem 6.4 and Theorem 6.6, we can get the Theorem 3.5 and Theorem 3.6
given in Le\'{o}n and Wang in \cite{lewa15}.
It shows that Theorem 6.4 and Theorem 6.6 can be regarded as an extension of two types of
Hamilton-Jacobi theorem for the distributional Hamiltonian system to that for the
system with magnetic, external force and control.
\end{rema}

\section{Hamilton-Jacobi Equations for a Nonholonomic Reduced Distributional CMH System }

It is well-known that the reduction of nonholonomically constrained mechanical systems
is a very important subject in geometric mechanics, and it is
regarded as a useful tool for simplifying and studying
concrete nonholonomic systems, see
Bates and $\acute{S}$niatycki \cite{basn93}, Cendra et al. \cite{cemara01},
Cushman et al. \cite{cudusn10}, Koiller \cite{ko92},
and Le\'{o}n and Wang \cite{lewa15} and so on,
for more details and development.\\

In this section, we shall consider the nonholonomic reduction and
Hamilton-Jacobi theory of a nonholonomic CMH
system with symmetry. We first give the definition of
a nonholonomic CMH system with symmetry.
Then, by using the similar method in Le\'{o}n and Wang \cite{lewa15}
and Bates and $\acute{S}$niatycki \cite{basn93}.
and by analyzing carefully the structure of dynamical vector field
of the nonholonomic CMH system with symmetry,
we give a geometric formulation of
the nonholonomic reduced distributional CMH system,
Moreover, we derive precisely the geometric constraint conditions of
the non-degenerate, and nonholonomic reduced distributional two-form
for the nonholonomic reducible dynamical vector field,
that is, the two types of Hamilton-Jacobi equation for the
nonholonomic reduced distributional CMH system,
which are an extension of the above two types of Hamilton-Jacobi equation
for the distributional CMH system
given in section 6 under nonholonomic reduction.\\

Assume that the Lie group $G$ acts smoothly on the manifold $Q$ by the left,
and we also consider the natural lifted actions on $TQ$ and $T^* Q$,
and assume that the cotangent lifted left action $\Phi^{T^\ast}:
G\times T^\ast Q\rightarrow T^\ast Q$ is free, proper and
symplectic with respect to the magnetic symplectic form
$\omega^B$ on $T^* Q$.
The orbit space $T^* Q/ G$ is a smooth manifold and the
canonical projection $\pi_{/G}: T^* Q \rightarrow T^* Q /G $ is
a surjective submersion. For the cotangent lifted left action
$\Phi^{T^\ast}: G\times T^\ast Q\rightarrow T^\ast Q$,
assume that $H: T^*Q \rightarrow \mathbb{R}$ is a
$G$-invariant Hamiltonian, and the fiber-preserving map
$F:T^\ast Q\rightarrow T^\ast Q$ and the control subset
$W$ of\; $T^\ast Q$ are both $G$-invariant,
and the $\mathcal{D}$-completely and
$\mathcal{D}$-regularly nonholonomic constraint $\mathcal{D}\subset
TQ$ is a $G$-invariant distribution for the tangent lifted left action $\Phi^{T}:
G\times TQ\rightarrow TQ$, that is, the tangent of group action maps
$\mathcal{D}_q$ to $\mathcal{D}_{gq}$ for any
$q\in Q $. A nonholonomic CMH system with symmetry
is 7-tuple $(T^*Q,G,\omega^B,\mathcal{D},H, F,W)$, which is an
CMH system with symmetry and $G$-invariant
nonholonomic constraint $\mathcal{D}$.\\

In the following we first consider the nonholonomic reduction
of a nonholonomic CMH system with symmetry
$(T^*Q,G,\omega^B,\mathcal{D},H, F,W)$.
Note that the Legendre transformation $\mathcal{F}L: TQ
\rightarrow T^*Q$ is a fiber-preserving map,
and $\mathcal{D}\subset TQ$ is $G$-invariant
for the tangent lifted left action $\Phi^{T}: G\times TQ\rightarrow TQ, $
then the constraint submanifold
$\mathcal{M}=\mathcal{F}L(\mathcal{D})\subset T^*Q$ is
$G$-invariant for the cotangent lifted left action $\Phi^{T^\ast}:
G\times T^\ast Q\rightarrow T^\ast Q$,
For the nonholonomic CMH system with symmetry
$(T^*Q,G, \omega^B,\mathcal{D},H,F,W)$,
in the same way, we define the distribution $\mathcal{F}$, which is the pre-image of the
nonholonomic constraints $\mathcal{D}$ for the map $T\pi_Q: TT^* Q
\rightarrow TQ$, that is, $\mathcal{F}=(T\pi_Q)^{-1}(\mathcal{D})$,
and the distribution $\mathcal{K}=\mathcal{F} \cap T\mathcal{M}$.
Moreover, we can also define the distributional two-form $\omega^B_\mathcal{K}$,
which is induced from the magnetic symplectic form $\omega^B$ on $T^* Q$, that is,
$\omega^B_\mathcal{K}= \tau_{\mathcal{K}}\cdot \omega^B_{\mathcal{M}},$ and
$\omega^B_{\mathcal{M}}= i_{\mathcal{M}}^* \omega^B $.
If the admissibility condition $\mathrm{dim}\mathcal{M}=
\mathrm{rank}\mathcal{F}$ and the compatibility condition
$T\mathcal{M}\cap \mathcal{F}^\bot= \{0\}$ hold, then
$\omega^B_\mathcal{K}$ is non-degenerate as a
bilinear form on each fibre of $\mathcal{K}$, there exists a vector
field $X^B_\mathcal{K}$ on $\mathcal{M}$ which takes values in the
constraint distribution $\mathcal{K}$, such that for the function $H_\mathcal{K}$,
the following distributional magnetic Hamiltonian equation holds, that is,
\begin{align}
\mathbf{i}_{X^B_\mathcal{K}}\omega^B_\mathcal{K}
=\mathbf{d}H_\mathcal{K},
\label{7.1} \end{align}
where the function $H_{\mathcal{K}}$ satisfies
$\mathbf{d}H_{\mathcal{K}}= \tau_{\mathcal{K}}\cdot \mathbf{d}H_{\mathcal {M}}$,
and $H_\mathcal{M}= \tau_{\mathcal{M}}\cdot H$
is the restriction of $H$ to $\mathcal{M}$, and
from the equation (7.1), we have that
$X^B_{\mathcal{K}}=\tau_{\mathcal{K}}\cdot X^B_H $.\\

In the following we define that the quotient space
$\bar{\mathcal{M}}=\mathcal{M}/G$ of the $G$-orbit in $\mathcal{M}$
is a smooth manifold with projection $\pi_{/G}:
\mathcal{M}\rightarrow \bar{\mathcal{M}}( \subset T^* Q /G),$ which
is a surjective submersion. The reduced magnetic symplectic form
$\omega^B_{\bar{\mathcal{M}}}= \pi^*_{/G} \cdot \omega^B_{\mathcal{M}}$
on $\bar{\mathcal{M}}$ is induced from the magnetic symplectic form $\omega^B_{\mathcal{M}}
= i_{\mathcal{M}}^* \omega^B $ on $\mathcal{M}$.
Since $G$ is the symmetry group of the system
$(T^*Q,G,\omega^B,\mathcal{D},H, F,W)$, all intrinsically
defined vector fields and distributions are pushed down to
$\bar{\mathcal{M}}$. In particular, the vector field $X^B_\mathcal{M}$
on $\mathcal{M}$ is pushed down to a vector field
$X^B_{\bar{\mathcal{M}}}=T\pi_{/G}\cdot X^B_\mathcal{M}$, and the
distribution $\mathcal{K}$ is pushed down to a distribution
$T\pi_{/G}\cdot \mathcal{K}$ on $\bar{\mathcal{M}}$, and the
Hamiltonian $H$ is pushed down to $h_{\bar{\mathcal{M}}}$, such that
$h_{\bar{\mathcal{M}}}\cdot \pi_{/G}=
\tau_{\mathcal{M}}\cdot H$. However, $\omega^B_\mathcal{K}$ need not
to be pushed down to a distributional two-form defined on $T\pi_{/G}\cdot
\mathcal{K}$, despite of the fact that $\omega^B_\mathcal{K}$ is
$G$-invariant. This is because there may be infinitesimal symmetry
$\eta_{\mathcal{K}}$ that lies in $\mathcal{M}$, such that
$\mathbf{i}_{\eta_\mathcal{K}} \omega^B_\mathcal{K}\neq 0$. From Bates
and $\acute{S}$niatycki \cite{basn93}, we know that in order to eliminate
this difficulty, $\omega^B_\mathcal{K}$ is restricted to a
sub-distribution $\mathcal{U}$ of $\mathcal{K}$ defined by
$$\mathcal{U}=\{u\in\mathcal{K} \; | \; \omega^B_\mathcal{K}(u,v)
=0,\quad \forall \; v \in \mathcal{V}\cap \mathcal{K}\},$$ where
$\mathcal{V}$ is the distribution on $\mathcal{M}$ tangent to the
orbits of $G$ in $\mathcal{M}$ and it is spanned by the infinitesimal
symmetries. Clearly, $\mathcal{U}$ and $\mathcal{V}$ are both
$G$-invariant, project down to $\bar{\mathcal{M}}$ and
$T\pi_{/G}\cdot \mathcal{V}=0$, and define the distribution $\bar{\mathcal{K}}$ by
$\bar{\mathcal{K}}= T\pi_{/G}\cdot \mathcal{U}$. Moreover, we take
that $\omega^B_\mathcal{U}= \tau_{\mathcal{U}}\cdot
\omega^B_{\mathcal{M}}$ is the restriction of the induced magnetic symplectic form
$\omega^B_{\mathcal{M}}$ on $T^*\mathcal{M}$ fibrewise to the
distribution $\mathcal{U}$, where $\tau_{\mathcal{U}}$ is the
restriction map to distribution $\mathcal{U}$, and the
$\omega^B_{\mathcal{U}}$ is pushed down to a
distributional two-form $\omega^B_{\bar{\mathcal{K}}}$ on
$\bar{\mathcal{K}}$, such that $\pi_{/G}^*
\omega^B_{\bar{\mathcal{K}}}= \omega^B_{\mathcal{U}}$.
We know that distributional two-form
$\omega^B_{\bar{\mathcal{K}}}$ is not a "true two-form"
on a manifold, which is called the nonholonomic reduced
distributional two-form to avoid any confusion.\\

From the above construction we know that,
if the admissibility condition $\mathrm{dim}\bar{\mathcal{M}}=
\mathrm{rank}\bar{\mathcal{F}}$ and the compatibility condition
$T\bar{\mathcal{M}} \cap \bar{\mathcal{F}}^\bot= \{0\}$ hold, where
$\bar{\mathcal{F}}^\bot$ denotes the symplectic orthogonal of
$\bar{\mathcal{F}}$ with respect to the reduced magnetic symplectic form
$\omega^B_{\bar{\mathcal{M}}}$, then the nonholonomic reduced
distributional two-form
$\omega^B_{\bar{\mathcal{K}}}$ is non-degenerate as a bilinear form on
each fibre of $\bar{\mathcal{K}}$, and hence there exists a vector field
$X^B_{\bar{\mathcal{K}}}$ on $\bar{\mathcal{M}}$ which takes values in
the constraint distribution $\bar{\mathcal{K}}$, such that the
reduced distributional magnetic Hamiltonian equation holds, that is,
\begin{align}
\mathbf{i}_{X^B_{\bar{\mathcal{K}}}}\omega^B_{\bar{\mathcal{K}}}
=\mathbf{d}h_{\bar{\mathcal{K}}},
\label{7.2} \end{align}
where $\mathbf{d}h_{\bar{\mathcal{K}}}$ is the restriction of
$\mathbf{d}h_{\bar{\mathcal{M}}}$ to $\bar{\mathcal{K}}$ and
the function $h_{\bar{\mathcal{K}}}:\bar{M}(\subset T^* Q/G)\rightarrow \mathbb{R}$ satisfies
$\mathbf{d}h_{\bar{\mathcal{K}}}= \tau_{\bar{\mathcal{K}}}\cdot \mathbf{d}h_{\bar{\mathcal{M}}}$,
and $h_{\bar{\mathcal{M}}}\cdot \pi_{/G}= H_{\mathcal{M}}$ and
$H_{\mathcal{M}}$ is the restriction of the Hamiltonian function $H$
to $\mathcal{M}$, and the function
$h_{\bar{\mathcal{M}}}:\bar{M}(\subset T^* Q/G)\rightarrow \mathbb{R}$.
In addition, from the distributional magnetic Hamiltonian equation (7.1),
$\mathbf{i}_{X^B_\mathcal{K}}\omega^B_\mathcal{K}=\mathbf{d}H_\mathcal
{K},$ we have that $X^B_{\mathcal{K}}=\tau_{\mathcal{K}}\cdot X^B_H, $
and from the reduced distributional magnetic Hamiltonian equation (7.2),
$\mathbf{i}_{X^B_{\bar{\mathcal{K}}}}\omega^B_{\bar{\mathcal{K}}}
=\mathbf{d}h_{\bar{\mathcal{K}}}$, we have that
$X^B_{\bar{\mathcal{K}}}
=\tau_{\bar{\mathcal{K}}}\cdot X^B_{h_{\bar{\mathcal{K}}}},$
where $ X^B_{h_{\bar{\mathcal{K}}}}$ is the magnetic Hamiltonian vector field of
the function $h_{\bar{\mathcal{K}}}$ with respect to the reduced magnetic symplectic
form $\omega^B_{\bar{\mathcal{M}}}$,
and the vector fields $X^B_{\mathcal{K}}$
and $X^B_{\bar{\mathcal{K}}}$ are $\pi_{/G}$-related,
that is, $X^B_{\bar{\mathcal{K}}}\cdot \pi_{/G}=T\pi_{/G}\cdot X^B_{\mathcal{K}}.$ \\

Moreover, if considering the external force $F$ and control subset $W$,
and we define the vector fields $F^B_\mathcal{K}
=\tau_{\mathcal{K}}\cdot \textnormal{vlift}(F_{\mathcal{M}})X^B_H,$
and for a control law $u\in W$,
$u^B_\mathcal{K}= \tau_{\mathcal{K}}\cdot  \textnormal{vlift}(u_{\mathcal{M}})X^B_H,$
where $F_\mathcal{M}= \tau_{\mathcal{M}}\cdot F$ and
$u_\mathcal{M}= \tau_{\mathcal{M}}\cdot u$ are the restrictions of
$F$ and $u$ to $\mathcal{M}$, that is, $F^B_\mathcal{K}$ and $u^B_\mathcal{K}$
are the restrictions of the changes of magnetic Hamiltonian vector field $X^B_H$
under the actions of $F_\mathcal{M}$ and $u_\mathcal{M}$ to $\mathcal{K}$,
then the 5-tuple $(\mathcal{K},\omega^B_{\mathcal{K}},
H_\mathcal{K}, F^B_\mathcal{K}, u^B_\mathcal{K})$
is a distributional CMH system corresponding to the nonholonomic CMH system with symmetry
$(T^*Q,G,\omega^B,\mathcal{D},H,F,u)$,
and the dynamical vector field of the distributional CMH system
can be expressed by
\begin{align}
\tilde{X}=X^B_{(\mathcal{K},\omega^B_{\mathcal{K}},
H_{\mathcal{K}}, F^B_{\mathcal{K}}, u^B_{\mathcal{K}})}
=X^B_\mathcal {K}+ F^B_{\mathcal{K}}+u^B_{\mathcal{K}},
\label{7.3} \end{align}
which is the synthetic
of the nonholonomic dynamical vector field $X^B_{\mathcal{K}}$ and
the vector fields $F^B_{\mathcal{K}}$ and $u^B_{\mathcal{K}}$.
Assume that the vector fields $F^B_\mathcal{K}$ and $u^B_\mathcal{K}$
on $\mathcal{M}$ are pushed down to the vector fields
$f^B_{\bar{\mathcal{M}}}= T\pi_{/G}\cdot F^B_\mathcal{K}$ and
$u^B_{\bar{\mathcal{M}}}=T\pi_{/G}\cdot u^B_\mathcal{K}$ on $\bar{\mathcal{M}}$.
Then we define that $f^B_{\bar{\mathcal{K}}}=T\tau_{\bar{\mathcal{K}}}\cdot f^B_{\bar{\mathcal{M}}}$ and
$u^B_{\bar{\mathcal{K}}}=T\tau_{\bar{\mathcal{K}}}\cdot u^B_{\bar{\mathcal{M}}},$
that is, $f^B_{\bar{\mathcal{K}}}$ and
$u^B_{\bar{\mathcal{K}}}$ are the restrictions of
$f^B_{\bar{\mathcal{M}}}$ and $u^B_{\bar{\mathcal{M}}}$ to $\bar{\mathcal{K}}$,
where $\tau_{\bar{\mathcal{K}}}$
is the restriction map to distribution $\bar{\mathcal{K}}$,
and $T\tau_{\bar{\mathcal{K}}}$ is the tangent map of $\tau_{\bar{\mathcal{K}}}$.
Then the 5-tuple $(\bar{\mathcal{K}},\omega^B_{\bar{\mathcal{K}}},
h_{\bar{\mathcal{K}}}, f^B_{\bar{\mathcal{K}}}, u^B_{\bar{\mathcal{K}}})$
is a nonholonomic reduced distributional CMH system of the nonholonomic
reducible CMH system with symmetry $(T^*Q,G,\omega^B,\mathcal{D},H,F,W)$,
as well as with a control law $u \in W$.
Thus, the geometrical formulation of a nonholonomic reduced distributional
CMH system may be summarized as follows.

\begin{defi} (Nonholonomic Reduced Distributional CMH System)
Assume that the 7-tuple \\ $(T^*Q,G,\omega^B,\mathcal{D},H,F,W)$ is a nonholonomic
reducible CMH system with symmetry, where $\omega^B$ is the magnetic
symplectic form on $T^* Q$, and $\mathcal{D}\subset TQ$ is a
$\mathcal{D}$-completely and $\mathcal{D}$-regularly nonholonomic
constraint of the system, and $\mathcal{D}$, $H, F$ and $W$ are all
$G$-invariant. If there exists a nonholonomic reduced distribution $\bar{\mathcal{K}}$,
an associated non-degenerate  and nonholonomic reduced
distributional two-form $\omega^B_{\bar{\mathcal{K}}}$
and a vector field $X^B_{\bar{\mathcal {K}}}$ on the reduced constraint
submanifold $\bar{\mathcal{M}}=\mathcal{M}/G, $ where
$\mathcal{M}=\mathcal{F}L(\mathcal{D})\subset T^*Q$, such that the
nonholonomic reduced distributional magnetic Hamiltonian equation
$ \mathbf{i}_{X^B_{\bar{\mathcal{K}}}}\omega^B_{\bar{\mathcal{K}}} =
\mathbf{d}h_{\bar{\mathcal{K}}}, $ holds,
where $\mathbf{d}h_{\bar{\mathcal{K}}}$ is the restriction of
$\mathbf{d}h_{\bar{\mathcal{M}}}$ to $\bar{\mathcal{K}}$ and
the function $h_{\bar{\mathcal{K}}}$ satisfies
$\mathbf{d}h_{\bar{\mathcal{K}}}= \tau_{\bar{\mathcal{K}}}\cdot \mathbf{d}h_{\bar{\mathcal{M}}}$
and $h_{\bar{\mathcal{M}}}\cdot \pi_{/G}= H_{\mathcal{M}}$,
and the vector fields $f^B_{\bar{\mathcal{K}}}=T\tau_{\bar{\mathcal{K}}}\cdot f^B_{\bar{\mathcal{M}}}$ and
$u^B_{\bar{\mathcal{K}}}=T\tau_{\bar{\mathcal{K}}}\cdot u^B_{\bar{\mathcal{M}}}$ as defined above.
Then the 5-tuple $(\bar{\mathcal{K}},\omega^B_{\bar{\mathcal {K}}},h_{\bar{\mathcal{K}}},
f^B_{\bar{\mathcal{K}}}, u^B_{\bar{\mathcal{K}}})$
is called a nonholonomic reduced distributional CMH system
of the nonholonomic reducible CMH system $(T^*Q,G,\omega^B,\mathcal{D},H,F,W)$
with a control law $u \in W$, and $X^B_{\bar{\mathcal {K}}}$ is
called a nonholonomic reduced dynamical vector field.
Denote by
\begin{align}
\hat{X}=X^B_{(\bar{\mathcal{K}},\omega^B_{\bar{\mathcal{K}}},
h_{\bar{\mathcal{K}}}, f^B_{\bar{\mathcal{K}}}, u^B_{\bar{\mathcal{K}}})}
=X^B_{\bar{\mathcal{K}}}+ f^B_{\bar{\mathcal{K}}}+u^B_{\bar{\mathcal{K}}}
\label{7.4} \end{align}
is the dynamical vector field of the
nonholonomic reduced distributional CMH system
$(\bar{\mathcal{K}},\omega^B_{\bar{\mathcal{K}}},h_{\bar{\mathcal{K}}},
f^B_{\bar{\mathcal{K}}}, \\ u^B_{\bar{\mathcal{K}}})$, which is the synthetic
of the nonholonomic reduced dynamical vector field $X^B_{\bar{\mathcal{K}}}$ and
the vector fields $F^B_{\bar{\mathcal{K}}}$ and $u^B_{\bar{\mathcal{K}}}$.
Under the above
circumstances, we refer to $(T^*Q,G,\omega^B,\mathcal{D},H,F,u)$ as a
nonholonomic reducible CMH system with the associated
distributional CMH system
$(\mathcal{K},\omega^B_{\mathcal {K}},H_{\mathcal{K}}, F^B_{\mathcal{K}}, u^B_{\mathcal{K}})$
and the nonholonomic reduced distributional CMH system
$(\bar{\mathcal{K}},\omega^B_{\bar{\mathcal{K}}},h_{\bar{\mathcal{K}}},
f^B_{\bar{\mathcal{K}}}, u^B_{\bar{\mathcal{K}}})$.
The dynamical vector fields
$\tilde{X}=X^B_{(\mathcal{K},\omega^B_{\mathcal{K}},
H_{\mathcal{K}}, F^B_{\mathcal{K}}, u^B_{\mathcal{K}})}$
and $\hat{X}= X^B_{(\bar{\mathcal{K}},\omega^B_{\bar{\mathcal{K}}},
h_{\bar{\mathcal{K}}}, f^B_{\bar{\mathcal{K}}}, u^B_{\bar{\mathcal{K}}})}$
are $\pi_{/G}$-related, that is,
$\hat{X}\cdot \pi_{/G}=T\pi_{/G}\cdot \tilde{X}.$
\end{defi}

For a given nonholonomic reducible CMH system
$(T^*Q,G,\omega^B,\mathcal{D},H,F,u)$ with the associated
distributional CMH system
$(\mathcal{K},\omega^B_{\mathcal {K}},H_{\mathcal{K}}, F^B_{\mathcal{K}}, u^B_{\mathcal{K}})$
and the nonholonomic reduced distributional CMH system
$(\bar{\mathcal{K}},\omega^B_{\bar{\mathcal{K}}},h_{\bar{\mathcal{K}}}, f^B_{\bar{\mathcal{K}}}, u^B_{\bar{\mathcal{K}}})$, the magnetic vector field $X^0= X^B_H-X_H, $
which is determined by the magnetic equation
$\mathbf{i}_{X^0}\omega=\mathbf{i}_{X^B_H}( \pi_Q^*B)$
on $T^*Q$. Note that vector fields $X^0$, $ X^B_H$, $X_H, $ and distribution
$\mathcal{K}$ are pushed down to
$\bar{\mathcal{M}}$, that is,  $X^0_{\bar{\mathcal{M}}}=T\pi_{/G}\cdot X^0_\mathcal{M}$,
$X^B_{\bar{\mathcal{M}}}=T\pi_{/G}\cdot X^B_\mathcal{M}$,
and $X_{\bar{\mathcal{M}}}=T\pi_{/G}\cdot X_\mathcal{M}$
where the vector field $X^0_\mathcal{M}$, $X^B_\mathcal{M}$
and $X_\mathcal{M}$ are the restrictions of $X^0$, $ X^B_H$ and $X_H$
on $\mathcal{M}$ and the
distribution $\mathcal{K}$ is pushed down to a distribution
$T\pi_{/G}\cdot \mathcal{K}$ on $\bar{\mathcal{M}}$,
and define the distribution $\bar{\mathcal{K}}$ by
$\bar{\mathcal{K}}= T\pi_{/G}\cdot \mathcal{U}$.
Denote by $X^0_{\bar{\mathcal {K}}}=\tau_{\bar{\mathcal {K}}}(X^0_{\bar{\mathcal{M}}})
= \tau_{\bar{\mathcal {K}}}(X^B_{\bar{\mathcal{M}}})- \tau_{\bar{\mathcal {K}}}(X_{\bar{\mathcal{M}}})
=X^B_{\bar{\mathcal {K}}}- X_{\bar{\mathcal {K}}},$
from the expression (7.4) of the dynamical
vector field of the nonholonomic reduced distributional CMH system
$(\bar{\mathcal{K}},\omega^B_{\bar{\mathcal{K}}},h_{\bar{\mathcal{K}}}, f^B_{\bar{\mathcal{K}}}, u^B_{\bar{\mathcal{K}}})$, we have that
\begin{align}\hat{X}
=X^B_{\bar{\mathcal{K}}}+ F^B_{\bar{\mathcal{K}}}+u^B_{\bar{\mathcal{K}}}
=X_{\bar{\mathcal{K}}}+ X^0_{\bar{\mathcal{K}}}+ F^B_{\bar{\mathcal{K}}}+u^B_{\bar{\mathcal{K}}}.
\label{7.5} \end{align}

If the vector fields $F^B_{\bar{\mathcal{K}}}$ and $u^B_{\bar{\mathcal{K}}}$ satisfy the following condition
\begin{equation}
 X^0_{\bar{\mathcal{K}}}+ F^B_{\bar{\mathcal{K}}}+u^B_{\bar{\mathcal{K}}}=0, \;\; \label{7.6}
\end{equation}
then from (7.5) we have that $X^B_{(\bar{\mathcal{K}},\omega^B_{\bar{\mathcal{K}}},
h_{\bar{\mathcal{K}}}, f^B_{\bar{\mathcal{K}}}, u^B_{\bar{\mathcal{K}}})}
=X_{\bar{\mathcal{K}}}, $ that is, in this case the dynamical vector
field of the nonholonomic reduced distributional CMH system is just the dynamical
vector field of the nonholonomic reduced canonical distributional Hamiltonian system
without the actions of magnetic, external force and control.
Thus, the condition (7.6) is called the magnetic vanishing condition for
the nonholonomic reduced distributional CMH system $(\bar{\mathcal{K}},\omega^B_{\bar{\mathcal{K}}},h_{\bar{\mathcal{K}}}, f^B_{\bar{\mathcal{K}}}, u^B_{\bar{\mathcal{K}}})$.\\

Since the non-degenerate and nonholonomic reduced distributional two-form
$\omega^B_{\bar{\mathcal{K}}}$ is not a "true two-form"
on a manifold, and it is not symplectic, and hence
the nonholonomic reduced distributional CMH system
$(\bar{\mathcal{K}},\omega^B_{\bar{\mathcal{K}}},h_{\bar{\mathcal{K}}}, f^B_{\bar{\mathcal{K}}}, u^B_{\bar{\mathcal{K}}})$ is not a Hamiltonian system,
and has no yet generating function,
and hence we can not describe the Hamilton-Jacobi equation for the nonholonomic reduced
distributional CMH system just like as in Theorem 1.1.
But, for a given nonholonomic reducible CMH system
$(T^*Q,G,\omega^B,\mathcal{D},H,F,u)$ with the associated
distributional CMH system
$(\mathcal{K},\omega^B_{\mathcal {K}},H_{\mathcal{K}}, F^B_{\mathcal{K}}, u^B_{\mathcal{K}})$
and the nonholonomic reduced distributional CMH system
$(\bar{\mathcal{K}},\omega^B_{\bar{\mathcal {K}}},h_{\bar{\mathcal{K}}}, f^B_{\bar{\mathcal{K}}}, u^B_{\bar{\mathcal{K}}})$, by using Lemma 3.4 and Lemma 6.3,
we can derive precisely
the geometric constraint conditions of the nonholonomic reduced distributional two-form
$\omega^B_{\bar{\mathcal{K}}}$ for the nonholonomic reducible dynamical vector field
$\tilde{X}=X^B_{(\mathcal{K},\omega^B_{\mathcal{K}},
H_{\mathcal{K}}, F^B_{\mathcal{K}}, u^B_{\mathcal{K}})}$,
that is, the two types of Hamilton-Jacobi equation for the
nonholonomic reduced distributional CMH system
$(\bar{\mathcal{K}},\omega^B_{\bar{\mathcal {K}}},h_{\bar{\mathcal{K}}}, f^B_{\bar{\mathcal{K}}}, u^B_{\bar{\mathcal{K}}})$.
At first, using the fact that the one-form $\gamma: Q
\rightarrow T^*Q $ satisfies the condition,
$\mathbf{d}\gamma=-B $ on $\mathcal{D}$ with respect to
$T\pi_Q: TT^* Q \rightarrow TQ, $
$\textmd{Im}(\gamma)\subset \mathcal{M}, $ and it is $G$-invariant,
as well as $ \textmd{Im}(T\gamma)\subset \mathcal{K}, $
we can prove the Type I of
Hamilton-Jacobi theorem for the nonholonomic reduced distributional
CMH system. For convenience, the maps involved in the
following theorem and its proof are shown in Diagram-7.
\begin{center}
\hskip 0cm \xymatrix{ & \mathcal{M} \ar[d]_{X^B_{\mathcal{K}}}
\ar[r]^{i_{\mathcal{M}}} & T^* Q \ar[d]_{X^B_{H}}
 \ar[r]^{\pi_Q}
  & Q \ar[d]_{\tilde{X}^\gamma} \ar[r]^{\gamma}
  & T^* Q \ar[d]_{\tilde{X}} \ar[r]^{\pi_{/G}}
  & T^* Q/G \ar[d]_{X^B_{h_{\bar{\mathcal{M}}}}}
  & \mathcal{\bar{M}} \ar[l]_{i_{\mathcal{\bar{M}}}} \ar[d]_{X^B_{\mathcal{\bar{K}}}}\\
  & \mathcal{K}
  & T(T^*Q) \ar[l]_{\tau_{\mathcal{K}}}
  & TQ \ar[l]_{T\gamma}
  & T(T^* Q) \ar[l]_{T\pi_Q} \ar[r]^{T\pi_{/G}}
  & T(T^* Q/G) \ar[r]^{\tau_{\mathcal{\bar{K}}}} & \mathcal{\bar{K}} }
\end{center}
$$\mbox{Diagram-7}$$

\begin{theo} (Type I of Hamilton-Jacobi Theorem for a Nonholonomic
Reduced Distributional CMH System)
For a given nonholonomic reducible CMH system
$(T^*Q,G,\omega^B,\mathcal{D},H,F,u)$ with the associated
distributional CMH system
$(\mathcal{K},\omega^B_{\mathcal {K}},H_{\mathcal{K}}, F^B_{\mathcal{K}}, u^B_{\mathcal{K}})$
and the nonholonomic reduced distributional CMH system
$(\bar{\mathcal{K}},\omega^B_{\bar{\mathcal{K}}},h_{\bar{\mathcal{K}}}, f^B_{\bar{\mathcal{K}}}, u^B_{\bar{\mathcal{K}}})$, assume that
$\gamma: Q \rightarrow T^*Q$ is an one-form on $Q$, and
$\tilde{X}^\gamma = T\pi_{Q}\cdot \tilde{X}\cdot \gamma$,
where $\tilde{X}=X^B_{(\mathcal{K},\omega^B_{\mathcal{K}},
H_{\mathcal{K}}, F^B_{\mathcal{K}}, u^B_{\mathcal{K}})}
=X^B_\mathcal {K}+ F^B_{\mathcal{K}}+u^B_{\mathcal{K}}$
is the dynamical vector field of the distributional CMH system
corresponding to the nonholonomic reducible CMH system with symmetry
$(T^*Q,G,\omega^B,\mathcal{D},H,F,u)$. Moreover,
assume that $\textmd{Im}(\gamma)\subset \mathcal{M}, $ and it is
$G$-invariant, $ \textmd{Im}(T\gamma)\subset \mathcal{K}, $ and
$\bar{\gamma}=\pi_{/G}(\gamma): Q \rightarrow T^* Q/G .$ If the
one-form $\gamma: Q \rightarrow T^*Q $ satisfies the condition,
$\mathbf{d}\gamma=-B $ on $\mathcal{D}$ with respect to
$T\pi_Q: TT^* Q \rightarrow TQ, $ then $\bar{\gamma}$ is a solution
of the equation $T\bar{\gamma}\cdot \tilde{X}^ \gamma =
X^B_{\bar{\mathcal{K}}}\cdot \bar{\gamma}. $ Here
$X^B_{\bar{\mathcal{K}}}$ is the nonholonomic reduced dynamical vector
field. The equation $T\bar{\gamma}\cdot \tilde{X}^ \gamma = X^B_{\bar{\mathcal{K}}}\cdot
\bar{\gamma},$ is called the Type I of Hamilton-Jacobi equation for
the nonholonomic reduced distributional CMH system
$(\bar{\mathcal{K}},\omega^B_{\bar{\mathcal{K}}},h_{\bar{\mathcal{K}}}, f^B_{\bar{\mathcal{K}}}, u^B_{\bar{\mathcal{K}}})$.
\end{theo}

\noindent{\bf Proof: } At first, for the dynamical vector field of the distributional CMH system
$(\mathcal{K},\omega^B_{\mathcal {K}}, H_{\mathcal{K}}, F^B_{\mathcal{K}}, u^B_{\mathcal{K}})$,
$\tilde{X}=X^B_{(\mathcal{K},\omega^B_{\mathcal{K}},
H_{\mathcal{K}}, F^B_{\mathcal{K}}, u^B_{\mathcal{K}})}
=X^B_\mathcal {K}+ F^B_{\mathcal{K}}+u^B_{\mathcal{K}}$,
and $F^B_{\mathcal{K}}=\tau_{\mathcal{K}}\cdot \textnormal{vlift}(F_{\mathcal{M}})X^B_H$,
and $u^B_{\mathcal{K}}=\tau_{\mathcal{K}}\cdot \textnormal{vlift}(u_{\mathcal{M}})X^B_H$,
note that $T\pi_{Q}\cdot \textnormal{vlift}(F_{\mathcal{M}})X^B_H=T\pi_{Q}\cdot \textnormal{vlift}(u_{\mathcal{M}})X^B_H=0, $
then we have that $T\pi_{Q}\cdot F^B_{\mathcal{K}}=T\pi_{Q}\cdot u^B_{\mathcal{K}}=0,$
and hence $T\pi_{Q}\cdot \tilde{X}\cdot \gamma=T\pi_{Q}\cdot X^B_{\mathcal{K}}\cdot \gamma. $
Moreover, from Theorem 6.4, we know that
$\gamma$ is a solution of the Type I of Hamilton-Jacobi equation
$T\gamma\cdot \tilde{X}^\gamma= X^B_{\mathcal{K}}\cdot \gamma .$ Next, we note that
$\textmd{Im}(\gamma)\subset \mathcal{M}, $ and it is $G$-invariant,
$ \textmd{Im}(T\gamma)\subset \mathcal{K}, $ and hence
$\textmd{Im}(T\bar{\gamma})\subset \bar{\mathcal{K}}, $ in this case,
$\pi^*_{/G}\cdot\omega^B_{\bar{\mathcal{K}}}\cdot\tau_{\bar{\mathcal{K}}}= \tau_{\mathcal{U}}\cdot
\omega^B_{\mathcal{M}}= \tau_{\mathcal{U}}\cdot i_{\mathcal{M}}^*\cdot
\omega^B, $ along $\textmd{Im}(T\bar{\gamma})$.
From the distributional Hamiltonian equation (7.1),
we have that $X^B_{\mathcal{K}}= \tau_{\mathcal{K}}\cdot X^B_H,$
and $\tau_{\mathcal{K}}\cdot X^B_{H}\cdot \gamma
= X^B_{\mathcal{K}}\cdot \gamma \in \mathcal{K}$.
Because the vector fields $X^B_{\mathcal{K}}$
and $X^B_{\bar{\mathcal{K}}}$ are $\pi_{/G}$-related,
$T\pi_{/G}(X^B_{\mathcal{K}})=X^B_{\bar{\mathcal{K}}}\cdot \pi_{/G}$,
and hence $\tau_{\bar{\mathcal{K}}}\cdot T\pi_{/G}(X^B_{\mathcal{K}}\cdot \gamma)
=\tau_{\bar{\mathcal{K}}}\cdot (T\pi_{/G}(X^B_{\mathcal{K}}))\cdot (\gamma)
= \tau_{\bar{\mathcal{K}}}\cdot (X^B_{\bar{\mathcal{K}}}\cdot \pi_{/G})\cdot (\gamma)
= \tau_{\bar{\mathcal{K}}}\cdot X^B_{\bar{\mathcal{K}}}\cdot \pi_{/G}(\gamma)
= X^B_{\bar{\mathcal{K}}}\cdot \bar{\gamma}.$
Thus, using the non-degenerate nonholonomic reduced distributional two-form
$\omega^B_{\bar{\mathcal{K}}}$, from Lemma 3.4(ii) and Lemma 6.3, if we take that
$v=X^B_{\mathcal{K}}\cdot \gamma \in \mathcal{K} (\subset \mathcal{F}),$
and for any $w \in \mathcal{F}, \; T\lambda(w)\neq 0, $ and
$\tau_{\bar{\mathcal{K}}}\cdot T\pi_{/G}\cdot w \neq 0, $ then we have that
\begin{align*}
& \omega^B_{\bar{\mathcal{K}}}(T\bar{\gamma} \cdot \tilde{X}^\gamma, \;
\tau_{\bar{\mathcal{K}}}\cdot T\pi_{/G} \cdot w)
= \omega^B_{\bar{\mathcal{K}}}(\tau_{\bar{\mathcal{K}}}\cdot T(\pi_{/G} \cdot
\gamma) \cdot \tilde{X}^\gamma, \; \tau_{\bar{\mathcal{K}}}\cdot T\pi_{/G} \cdot w )\\
& = \pi^*_{/G}
\cdot \omega^B_{\bar{\mathcal{K}}}\cdot\tau_{\bar{\mathcal{K}}}
(T\gamma \cdot T\pi_{Q}\cdot \tilde{X}\cdot \gamma, \; w) =
\tau_{\mathcal{U}}\cdot i^*_{\mathcal{M}} \cdot\omega^B(T\gamma \cdot
T\pi_{Q}\cdot X^B_{\mathcal{K}}\cdot \gamma, \; w)\\
& = \tau_{\mathcal{U}}\cdot i^*_{\mathcal{M}} \cdot
\omega^B(T(\gamma \cdot \pi_Q)\cdot X^B_{\mathcal{K}}\cdot \gamma, \; w) \\
& = \tau_{\mathcal{U}}\cdot i^*_{\mathcal{M}} \cdot
(\omega^B (X^B_{\mathcal{K}}\cdot \gamma, \; w-T(\gamma \cdot \pi_Q)\cdot w)
- (\mathbf{d}\gamma+B)(T\pi_{Q}
(X^B_{\mathcal{K}}\cdot \gamma), \; T\pi_{Q}(w)))\\
& = \tau_{\mathcal{U}}\cdot i^*_{\mathcal{M}} \cdot \omega^B (X^B_{\mathcal{K}} \cdot
\gamma, \; w) - \tau_{\mathcal{U}}\cdot i^*_{\mathcal{M}} \cdot
\omega^B (X^B_{\mathcal{K}}\cdot \gamma, \; T(\gamma \cdot \pi_Q) \cdot w)\\
& \;\;\;\;\;\;
- \tau_{\mathcal{U}}\cdot i^*_{\mathcal{M}} \cdot (\mathbf{d}\gamma+B)(T\pi_{Q}
(X^B_{\mathcal{K}}\cdot \gamma), \; T\pi_{Q}(w))\\
& =\pi^*_{/G}\cdot \omega^B_{\bar{\mathcal{K}}}\cdot\tau_{\bar{\mathcal{K}}}(X^B_{\mathcal{K}}\cdot \gamma, \;
w) - \pi^*_{/G}\cdot \omega^B_{\bar{\mathcal{K}}}\cdot\tau_{\bar{\mathcal{K}}}(X^B_{\mathcal{K}}\cdot \gamma,
\; T(\gamma \cdot \pi_Q) \cdot w)\\
& \;\;\;\;\;\; - \tau_{\mathcal{U}}\cdot i^*_{\mathcal{M}}
\cdot(\mathbf{d}\gamma+B)(T\pi_{Q}(X^B_{\mathcal{K}}\cdot \gamma), \; T\pi_{Q}(w))\\
& = \omega^B_{\bar{\mathcal{K}}}(\tau_{\bar{\mathcal{K}}}\cdot T\pi_{/G}(X^B_{\mathcal{K}}\cdot \gamma), \;
\tau_{\bar{\mathcal{K}}}\cdot T\pi_{/G} \cdot w) - \omega^B_{\bar{\mathcal{K}}}(\tau_{\bar{\mathcal{K}}}\cdot T\pi_{/G}(X^B_{\mathcal{K}}\cdot \gamma), \;
\tau_{\bar{\mathcal{K}}}\cdot T(\pi_{/G} \cdot\gamma) \cdot T\pi_{Q}(w))\\
& \;\;\;\;\;\; -\tau_{\mathcal{U}}\cdot i^*_{\mathcal{M}}
\cdot(\mathbf{d}\gamma+B)(T\pi_{Q}(X^B_{\mathcal{K}}\cdot \gamma), \; T\pi_{Q}(w))\\
& = \omega^B_{\bar{\mathcal{K}}}(X^B_{\bar{\mathcal{K}}} \cdot
\bar{\gamma}, \; \tau_{\bar{\mathcal{K}}}\cdot T\pi_{/G} \cdot w)-
\omega^B_{\bar{\mathcal{K}}}(X^B_{\bar{\mathcal{K}}} \cdot
\bar{\gamma}, \; T\bar{\gamma} \cdot T\pi_{Q}(w)) \\
& \;\;\;\;\;\; - \tau_{\mathcal{U}}\cdot
i^*_{\mathcal{M}} \cdot(\mathbf{d}\gamma+B)
(T\pi_{Q}(X^B_{\mathcal{K}}\cdot \gamma), \; T\pi_{Q}(w)),
\end{align*}
where we have used that $\tau_{\bar{\mathcal{K}}}\cdot T\pi_{/G}(X^B_{\mathcal{K}}\cdot \gamma)
=X^B_{\bar{\mathcal{K}}}\cdot \bar{\gamma}, $ and
$\tau_{\bar{\mathcal{K}}}\cdot T\bar{\gamma}=T\bar{\gamma}, $ since
$\textmd{Im}(T\bar{\gamma})\subset \bar{\mathcal{K}}. $
If the one-form $\gamma: Q \rightarrow T^*Q $ satisfies the condition,
$\mathbf{d}\gamma=-B $ on $\mathcal{D}$ with respect to
$T\pi_Q: TT^* Q \rightarrow TQ, $ then we have that
$(\mathbf{d}\gamma +B)(T\pi_{Q}(X^B_{\mathcal{K}}\cdot \gamma), \; T\pi_{Q}(w))=0, $
since $X^B_{\mathcal{K}}\cdot \gamma, \; w \in \mathcal{F},$ and
$T\pi_{Q}(X^B_{\mathcal{K}}\cdot \gamma), \; T\pi_{Q}(w) \in \mathcal{D}, $ and hence
$$
\tau_{\mathcal{U}}\cdot
i_{\mathcal{M}}^* \cdot (\mathbf{d}\gamma +B)(T\pi_{Q}(X^B_{\mathcal{K}}\cdot \gamma),
\; T\pi_{Q}(w))=0,
$$
and
\begin{equation}
\omega^B_{\bar{\mathcal{K}}}(T\bar{\gamma} \cdot \tilde{X}^\gamma, \;
\tau_{\bar{\mathcal{K}}}\cdot T\pi_{/G} \cdot w)
- \omega^B_{\bar{\mathcal{K}}}(X^B_{\bar{\mathcal{K}}} \cdot
\bar{\gamma}, \; \tau_{\bar{\mathcal{K}}}\cdot T\pi_{/G} \cdot w)
= -\omega^B_{\bar{\mathcal{K}}}(X^B_{\bar{\mathcal{K}}} \cdot
\bar{\gamma}, \; T\bar{\gamma} \cdot T\pi_{Q}(w)).
\label{7.7} \end{equation}
If $\bar{\gamma}$ satisfies the equation $
T\bar{\gamma}\cdot \tilde{X}^ \gamma
= X^B_{\bar{\mathcal{K}}}\cdot \bar{\gamma} ,$
from Lemma 3.4(i) we know that the right side of (7.7) becomes that
\begin{align*}
 -\omega^B_{\bar{\mathcal{K}}}(X^B_{\bar{\mathcal{K}}} \cdot
\bar{\gamma}, \; T\bar{\gamma} \cdot T\pi_{Q}(w))
& = -\omega^B_{\bar{\mathcal{K}}}\cdot\tau_{\bar{\mathcal{K}}}(T\bar{\gamma}\cdot \tilde{X}^\gamma, \; T\bar{\gamma} \cdot T\pi_{Q}(w))\\
& = -\bar{\gamma}^*\omega^B_{\bar{\mathcal{K}}}\cdot\tau_{\bar{\mathcal{K}}}
(T\pi_{Q} \cdot \tilde{X}\cdot \gamma, \; T\pi_{Q}(w))\\
& = - \gamma^* \cdot \pi^*_{/G}\cdot \omega^B_{\bar{\mathcal{K}}}\cdot\tau_{\bar{\mathcal{K}}}
(T\pi_{Q} \cdot X^B_{\mathcal{K}} \cdot \gamma, \; T\pi_{Q}(w))\\
& = - \gamma^* \cdot \tau_{\mathcal{U}}\cdot
i_{\mathcal{M}}^* \cdot \omega^B (T\pi_{Q}(X^B_{\mathcal{K}}\cdot\gamma), \; T\pi_{Q}(w))\\
& = -\tau_{\mathcal{U}}\cdot
i_{\mathcal{M}}^* \cdot\gamma^*\omega^B ( T\pi_{Q}(X^B_{\mathcal{K}}\cdot\gamma), \; T\pi_{Q}(w))\\
& = \tau_{\mathcal{U}}\cdot i_{\mathcal{M}}^* \cdot
(\mathbf{d}\gamma +B)(T\pi_{Q}( X_{\mathcal{K}}\cdot\gamma ), \; T\pi_{Q}(w))=0,
\end{align*}
where $\gamma^*\cdot \tau_{\mathcal{U}}\cdot i^*_{\mathcal{M}}
\cdot \omega^B= \tau_{\mathcal{U}}\cdot i^*_{\mathcal{M}}
\cdot\gamma^*\cdot \omega^B, $ because $\textmd{Im}(\gamma)\subset
\mathcal{M}. $
But, since the nonholonomic reduced distributional two-form
$\omega^B_{\bar{\mathcal{K}}}$ is non-degenerate,
the left side of (7.7) equals zero, only when
$\bar{\gamma}$ satisfies the equation $
T\bar{\gamma}\cdot \tilde{X}^ \gamma = X^B_{\bar{\mathcal{K}}}\cdot
\bar{\gamma} .$ Thus,
if the one-form $\gamma: Q \rightarrow T^*Q $ satisfies the condition,
$\mathbf{d}\gamma=-B $ on $\mathcal{D}$ with respect to
$T\pi_Q: TT^* Q \rightarrow TQ, $ then $\bar{\gamma}$ must be a solution of
the Type I of Hamilton-Jacobi equation
$T\bar{\gamma}\cdot \tilde{X}^ \gamma = X^B_{\bar{\mathcal{K}}}\cdot
\bar{\gamma}. $
\hskip 0.3cm $\blacksquare$\\

Next, for any $G$-invariant symplectic map $\varepsilon: T^* Q \rightarrow T^* Q $
with respect to $\omega^B$, we can prove
the following Type II of
Hamilton-Jacobi theorem for the nonholonomic reduced distributional CMH system.
For convenience, the maps involved in the following theorem and its
proof are shown in Diagram-8.
\begin{center}
\hskip 0cm \xymatrix{ & \mathcal{M} \ar[d]_{X^B_{\mathcal{K}}}
\ar[r]^{i_{\mathcal{M}}} & T^* Q \ar[d]_{X^B_{H\cdot \varepsilon}}
\ar[dr]^{\tilde{X}^\varepsilon} \ar[r]^{\pi_Q}
  & Q \ar[r]^{\gamma} & T^* Q \ar[d]_{\tilde{X}} \ar[r]^{\pi_{/G}}
  & T^* Q/G \ar[d]_{X^B_{h_{\bar{\mathcal{M}}}}}
  & \mathcal{\bar{M}} \ar[l]_{i_{\mathcal{\bar{M}}}} \ar[d]_{X^B_{\mathcal{\bar{K}}}}\\
  & \mathcal{K}
  & T(T^*Q) \ar[l]_{\tau_{\mathcal{K}}}
  & TQ \ar[l]_{T\gamma}
  & T(T^* Q) \ar[l]_{T\pi_Q} \ar[r]^{T\pi_{/G}}
  & T(T^* Q/G) \ar[r]^{\tau_{\mathcal{\bar{K}}}} & \mathcal{\bar{K}} }
\end{center}
$$\mbox{Diagram-8}$$

\begin{theo} (Type II of Hamilton-Jacobi Theorem for a Nonholonomic
Reduced Distributional CMH System)
For a given nonholonomic reducible CMH system
$(T^*Q,G,\omega^B,\mathcal{D},H,F,u)$ with the associated
distributional CMH system
$(\mathcal{K},\omega^B_{\mathcal {K}},H_{\mathcal{K}}, F^B_{\mathcal{K}}, u^B_{\mathcal{K}})$
and the nonholonomic reduced distributional CMH system
$(\bar{\mathcal{K}},\omega^B_{\bar{\mathcal{K}}},h_{\bar{\mathcal{K}}},
f^B_{\bar{\mathcal{K}}}, u^B_{\bar{\mathcal{K}}})$, assume that
$\gamma: Q \rightarrow T^*Q$ is an one-form on $Q$, and $\lambda=
\gamma \cdot \pi_{Q}: T^* Q \rightarrow T^* Q, $ and for any $G$-invariant
symplectic map $\varepsilon: T^* Q \rightarrow T^* Q $
with respect to $\omega^B$, denote by
$\tilde{X}^\varepsilon = T\pi_{Q}\cdot \tilde{X}\cdot \varepsilon$,
where $\tilde{X}=X^B_{(\mathcal{K},\omega^B_{\mathcal{K}},
H_{\mathcal{K}}, F^B_{\mathcal{K}}, u^B_{\mathcal{K}})}
=X^B_\mathcal {K}+ F^B_{\mathcal{K}}+u^B_{\mathcal{K}}$
is the dynamical vector field of the distributional CMH system
corresponding to the nonholonomic reducible CMH system with symmetry
$(T^*Q,G,\omega^B,\mathcal{D},H,F,u)$. Moreover,
assume that $\textmd{Im}(\gamma)\subset \mathcal{M}, $ and it is
$G$-invariant, $\varepsilon(\mathcal{M})\subset \mathcal{M}$,
$ \textmd{Im}(T\gamma)\subset \mathcal{K}, $ and
$\bar{\gamma}=\pi_{/G}(\gamma): Q \rightarrow T^* Q/G $, and
$\bar{\lambda}=\pi_{/G}(\lambda): T^* Q \rightarrow T^* Q/G, $ and
$\bar{\varepsilon}=\pi_{/G}(\varepsilon): T^* Q \rightarrow T^* Q/G. $ Then
$\varepsilon$ and $\bar{\varepsilon}$ satisfy the equation
$\tau_{\bar{\mathcal{K}}}\cdot T\bar{\varepsilon}\cdot X^B_{h_{\bar{\mathcal{K}}}\cdot
\bar{\varepsilon}}= T\bar{\lambda} \cdot \tilde{X}\cdot \varepsilon, $ if and only if they satisfy the
equation $T\bar{\gamma}\cdot \tilde{X}^ \varepsilon =
X^B_{\bar{\mathcal{K}}}\cdot \bar{\varepsilon}. $ Here
$ X^B_{h_{\bar{\mathcal{K}}} \cdot\bar{\varepsilon}}$ is the magnetic Hamiltonian
vector field of the function $h_{\bar{\mathcal{K}}}\cdot \bar{\varepsilon}: T^* Q\rightarrow
\mathbb{R}, $ and $X^B_{\bar{\mathcal{K}}}$ is the nonholonomic
reduced dynamical vector field. The equation $
T\bar{\gamma}\cdot \tilde{X}^\varepsilon = X^B_{\bar{\mathcal{K}}}\cdot
\bar{\varepsilon},$ is called the Type II of Hamilton-Jacobi equation for the
nonholonomic reduced distributional CMH system
$(\bar{\mathcal{K}},\omega^B_{\bar{\mathcal{K}}},h_{\bar{\mathcal{K}}}, f^B_{\bar{\mathcal{K}}}, u^B_{\bar{\mathcal{K}}})$.
\end{theo}

\noindent{\bf Proof: } In the same way, for the dynamical vector field of the distributional CMH system
$(\mathcal{K},\omega^B_{\mathcal {K}}, H_{\mathcal{K}}, F^B_{\mathcal{K}}, u^B_{\mathcal{K}})$,
$\tilde{X}=X^B_{(\mathcal{K},\omega^B_{\mathcal{K}},
H_{\mathcal{K}}, F^B_{\mathcal{K}}, u^B_{\mathcal{K}})}
=X^B_\mathcal {K}+ F^B_{\mathcal{K}}+u^B_{\mathcal{K}}$,
and $F^B_{\mathcal{K}}=\tau_{\mathcal{K}}\cdot \textnormal{vlift}(F_{\mathcal{M}})X^B_H$,
and $u^B_{\mathcal{K}}=\tau_{\mathcal{K}}\cdot \textnormal{vlift}(u_{\mathcal{M}})X^B_H$,
note that $T\pi_{Q}\cdot \textnormal{vlift}(F_{\mathcal{M}})X^B_H=T\pi_{Q}\cdot \textnormal{vlift}(u_{\mathcal{M}})X^B_H=0, $
then we have that $T\pi_{Q}\cdot F^B_{\mathcal{K}}=T\pi_{Q}\cdot u^B_{\mathcal{K}}=0,$
and hence $T\pi_{Q}\cdot \tilde{X}\cdot \varepsilon
=T\pi_{Q}\cdot X^B_{\mathcal{K}}\cdot \varepsilon. $
Next, we note that
$\textmd{Im}(\gamma)\subset \mathcal{M}, $ and it is $G$-invariant,
$ \textmd{Im}(T\gamma)\subset \mathcal{K}, $ and hence
$\textmd{Im}(T\bar{\gamma})\subset \bar{\mathcal{K}}, $ in this case,
$\pi^*_{/G}\cdot\omega^B_{\bar{\mathcal{K}}}\cdot\tau_{\bar{\mathcal{K}}}= \tau_{\mathcal{U}}\cdot
\omega^B_{\mathcal{M}}= \tau_{\mathcal{U}}\cdot i_{\mathcal{M}}^*\cdot
\omega^B, $ along $\textmd{Im}(T\bar{\gamma})$.
Moreover, from the distributional magnetic Hamiltonian equation (7.1),
we have that $X^B_{\mathcal{K}}= \tau_{\mathcal{K}}\cdot X^B_H.$
Note that $\varepsilon(\mathcal{M})\subset \mathcal{M},$ and
$T\pi_{Q}(X^B_H\cdot \varepsilon(q,p))\in
\mathcal{D}_{q}, \; \forall q \in Q, \; (q,p) \in \mathcal{M}(\subset T^* Q), $
and hence $X^B_H\cdot \varepsilon \in \mathcal{F}$ along $\varepsilon$,
and $\tau_{\mathcal{K}}\cdot X^B_{H}\cdot \varepsilon
= X^B_{\mathcal{K}}\cdot \varepsilon \in \mathcal{K}$.
Because the vector fields $X^B_{\mathcal{K}}$
and $X^B_{\bar{\mathcal{K}}}$ are $\pi_{/G}$-related, then
$T\pi_{/G}(X^B_{\mathcal{K}})=X^B_{\bar{\mathcal{K}}}\cdot \pi_{/G}$,
and hence $\tau_{\bar{\mathcal{K}}}\cdot T\pi_{/G}(X^B_{\mathcal{K}}\cdot \varepsilon)
=\tau_{\bar{\mathcal{K}}}\cdot (T\pi_{/G}(X^B_{\mathcal{K}}))\cdot (\varepsilon)
= \tau_{\bar{\mathcal{K}}}\cdot (X^B_{\bar{\mathcal{K}}}\cdot \pi_{/G})\cdot (\varepsilon)
= \tau_{\bar{\mathcal{K}}}\cdot X^B_{\bar{\mathcal{K}}}\cdot \pi_{/G}(\varepsilon)
= X^B_{\bar{\mathcal{K}}}\cdot \bar{\varepsilon}.$
Thus, using the non-degenerate and nonholonomic reduced distributional two-form
$\omega^B_{\bar{\mathcal{K}}}$, from Lemma 3.4 and Lemma 6.3, if we take that
$v=X^B_{\mathcal{K}}\cdot \varepsilon \in \mathcal{K} (\subset \mathcal{F}),$
and for any $w \in \mathcal{F}, \; T\lambda(w)\neq 0 $,
$\tau_{\bar{\mathcal{K}}}\cdot T\pi_{/G}\cdot w \neq 0, $
and $\tau_{\bar{\mathcal{K}}}\cdot T\pi_{/G}\cdot T\lambda(w) \neq 0, $
then we have that
\begin{align*}
& \omega^B_{\bar{\mathcal{K}}}(T\bar{\gamma} \cdot \tilde{X}^\varepsilon, \;
\tau_{\bar{\mathcal{K}}}\cdot T\pi_{/G} \cdot w)
= \omega^B_{\bar{\mathcal{K}}}(\tau_{\bar{\mathcal{K}}}\cdot T(\pi_{/G} \cdot
\gamma) \cdot \tilde{X}^\varepsilon, \; \tau_{\bar{\mathcal{K}}}\cdot T\pi_{/G} \cdot w )\\
& = \pi^*_{/G}\cdot \omega^B_{\bar{\mathcal{K}}}
\cdot\tau_{\bar{\mathcal{K}}}(T\gamma \cdot \tilde{X}^\varepsilon, \; w) =
\tau_{\mathcal{U}}\cdot i^*_{\mathcal{M}} \cdot \omega^B
(T\gamma \cdot T\pi_Q \cdot \tilde{X} \cdot \varepsilon, \; w)\\
& = \tau_{\mathcal{U}}\cdot i^*_{\mathcal{M}} \cdot
\omega^B (T(\gamma \cdot \pi_Q)\cdot X^B_{\mathcal{K}}\cdot \varepsilon, \; w) \\
& = \tau_{\mathcal{U}}\cdot i^*_{\mathcal{M}} \cdot
(\omega^B (X^B_{\mathcal{K}}\cdot \varepsilon, \; w-T(\gamma \cdot \pi_Q)\cdot w)
- (\mathbf{d}\gamma+B)(T\pi_{Q}(X^B_{\mathcal{K}}\cdot \varepsilon), \; T\pi_{Q}(w)))\\
& = \tau_{\mathcal{U}}\cdot i^*_{\mathcal{M}} \cdot \omega^B (X^B_{\mathcal{K}}\cdot
\varepsilon, \; w) - \tau_{\mathcal{U}}\cdot i^*_{\mathcal{M}} \cdot
\omega^B (X^B_{\mathcal{K}}\cdot \varepsilon, \; T\lambda \cdot w) \\
& \;\;\;\;\;\;
- \tau_{\mathcal{U}}\cdot i^*_{\mathcal{M}} \cdot (\mathbf{d}\gamma+B)
(T\pi_{Q}(X^B_{\mathcal{K}}\cdot \varepsilon), \; T\pi_{Q}(w))\\
& =\pi^*_{/G}\cdot \omega^B_{\bar{\mathcal{K}}}\cdot \tau_{\bar{\mathcal{K}}}
(X^B_{\mathcal{K}}\cdot \varepsilon, \;
w) - \pi^*_{/G}\cdot \omega^B_{\bar{\mathcal{K}}}\cdot \tau_{\bar{\mathcal{K}}}
(X^B_{\mathcal{K}}\cdot \varepsilon,
\; T\lambda \cdot w)+ \tau_{\mathcal{U}}\cdot i^*_{\mathcal{M}}
\cdot \lambda^* \omega^B (X^B_{\mathcal{K}}\cdot \varepsilon, \; w)\\
& = \omega^B_{\bar{\mathcal{K}}}(\tau_{\bar{\mathcal{K}}}\cdot T\pi_{/G}
(X^B_{\mathcal{K}}\cdot \varepsilon), \;
\tau_{\bar{\mathcal{K}}}\cdot T\pi_{/G} \cdot w)
- \omega^B_{\bar{\mathcal{K}}}(\tau_{\bar{\mathcal{K}}}\cdot T\pi_{/G}(X^B_{\mathcal{K}}\cdot
\varepsilon), \; \tau_{\bar{\mathcal{K}}}\cdot T(\pi_{/G} \cdot\lambda) \cdot w)\\
& \;\;\;\;\;\; +\pi^*_{/G}\cdot \omega^B_{\bar{\mathcal{K}}}\cdot
\tau_{\bar{\mathcal{K}}}(T\lambda\cdot X^B_{\mathcal{K}}\cdot \varepsilon, \;
T\lambda \cdot w)\\
& = \omega^B_{\bar{\mathcal{K}}}(X^B_{\bar{\mathcal{K}}} \cdot
\bar{\varepsilon}, \; \tau_{\bar{\mathcal{K}}}\cdot T\pi_{/G} \cdot w)-
\omega^B_{\bar{\mathcal{K}}}(X^B_{\bar{\mathcal{K}}}\cdot
\bar{\varepsilon}, \;  \tau_{\bar{\mathcal{K}}}\cdot T\bar{\lambda} \cdot w)+ \omega^B_{\bar{\mathcal{K}}}
(T\bar{\lambda}\cdot X^B_{\mathcal{K}}\cdot \varepsilon,
\;  \tau_{\bar{\mathcal{K}}}\cdot T\bar{\lambda} \cdot w),
\end{align*}
where we have used that $\tau_{\bar{\mathcal{K}}}\cdot T\pi_{/G}(X^B_{\mathcal{K}}\cdot \varepsilon)
=X^B_{\bar{\mathcal{K}}}\cdot \bar{\varepsilon}, $
and $\tau_{\bar{\mathcal{K}}}\cdot T\bar{\gamma}=T\bar{\gamma},$
and $\tau_{\bar{\mathcal{K}}}\cdot T\pi_{/G}\cdot T\lambda
=\tau_{\bar{\mathcal{K}}}\cdot T\bar{\lambda}=T\bar{\lambda},$ since
$\textmd{Im}(T\bar{\gamma})\subset \bar{\mathcal{K}}.$
From the nonholonomic reduced distributional magnetic Hamiltonian equation (7.2),
$\mathbf{i}_{X^B_{\bar{\mathcal{K}}}}\omega^B_{\bar{\mathcal{K}}} =
\mathbf{d}h_{\bar{\mathcal{K}}},$ we have that $X^B_{\bar{\mathcal{K}}}
=\tau_{\bar{\mathcal{K}}}\cdot X^B_{h_{\bar{\mathcal{K}}}},$
where $ X^B_{h_{\bar{\mathcal{K}}}}$ is the magnetic Hamiltonian vector field of
the function $h_{\bar{\mathcal{K}}}: \bar{M}(\subset T^* Q/G)\rightarrow \mathbb{R}.$
Note that
$\varepsilon: T^* Q \rightarrow T^* Q $ is symplectic
with respect to $\omega^B$, and
$\bar{\varepsilon}=\pi_{/G}(\varepsilon): T^* Q \rightarrow T^* Q/G$
is also symplectic along $\bar{\varepsilon}$, and
hence $X^B_{h_{\bar{\mathcal{K}}}}\cdot \bar{\varepsilon}
= T\bar{\varepsilon} \cdot X^B_{h_{\bar{\mathcal{K}}} \cdot
\bar{\varepsilon}}, $ along $\bar{\varepsilon}$, and hence
$X^B_{\bar{\mathcal{K}}}\cdot \bar{\varepsilon}
=\tau_{\bar{\mathcal{K}}}\cdot X^B_{h_{\bar{\mathcal{K}}}} \cdot\bar{\varepsilon}
= \tau_{\bar{\mathcal{K}}}\cdot T\bar{\varepsilon} \cdot X^B_{h_{\bar{\mathcal{K}}}
\cdot \bar{\varepsilon}}, $ along $\bar{\varepsilon}$.
Note that
$T\bar{\lambda} \cdot X^B_\mathcal{K}\cdot \varepsilon
=T\pi_{/G}\cdot T\lambda \cdot X^B_\mathcal{K}\cdot \varepsilon
=T\pi_{/G}\cdot T\gamma \cdot T\pi_Q\cdot X^B_\mathcal{K}\cdot \varepsilon
=T\pi_{/G}\cdot T\gamma \cdot T\pi_Q\cdot \tilde{X}\cdot \varepsilon
=T\pi_{/G}\cdot T\lambda\cdot \tilde{X}\cdot \varepsilon
=T\bar{\lambda} \cdot \tilde{X}\cdot \varepsilon.$
Then we have that
\begin{align*}
& \omega^B_{\bar{\mathcal{K}}}(T\bar{\gamma} \cdot \tilde{X}^\varepsilon, \;
\tau_{\bar{\mathcal{K}}}\cdot T\pi_{/G}\cdot w)-
\omega^B_{\bar{\mathcal{K}}}(X^B_{\bar{\mathcal{K}}}\cdot \bar{\varepsilon},
\; \tau_{\bar{\mathcal{K}}}\cdot T\pi_{/G} \cdot w) \nonumber \\
&= - \omega^B_{\bar{\mathcal{K}}}(X^B_{\bar{\mathcal{K}}}\cdot
\bar{\varepsilon}, \;  \tau_{\bar{\mathcal{K}}}\cdot T\bar{\lambda} \cdot w)+ \omega^B_{\bar{\mathcal{K}}}
(T\bar{\lambda}\cdot X^B_{\mathcal{K}}\cdot \varepsilon,
\;  \tau_{\bar{\mathcal{K}}}\cdot T\bar{\lambda} \cdot w) \\
& =-\omega^B_{\bar{\mathcal{K}}}(\tau_{\bar{\mathcal{K}}} \cdot X^B_{h_{\bar{\mathcal{K}}}}\cdot
\bar{\varepsilon}, \; \tau_{\bar{\mathcal{K}}}\cdot T\bar{\lambda} \cdot w)+ \omega^B_{\bar{\mathcal{K}}}
(T\bar{\lambda}\cdot \tilde{X} \cdot \varepsilon,
\; \tau_{\bar{\mathcal{K}}}\cdot T\bar{\lambda} \cdot w)\\
& = -\omega^B_{\bar{\mathcal{K}}}(\tau_{\bar{\mathcal{K}}} \cdot T\bar{\varepsilon}
\cdot X^B_{h_{\bar{\mathcal{K}}} \cdot \bar{\varepsilon}}, \;
 \tau_{\bar{\mathcal{K}}}\cdot T\bar{\lambda} \cdot w)+ \omega^B_{\bar{\mathcal{K}}}
(T\bar{\lambda}\cdot \tilde{X} \cdot \varepsilon,
\;  \tau_{\bar{\mathcal{K}}}\cdot T\bar{\lambda} \cdot w)\\
& = \omega^B_{\bar{\mathcal{K}}}(T\bar{\lambda}\cdot \tilde{X} \cdot \varepsilon- \tau_{\bar{\mathcal{K}}}\cdot T\bar{\varepsilon} \cdot X^B_{h_{\bar{\mathcal{K}}}\cdot \bar{\varepsilon}},
\;  \tau_{\bar{\mathcal{K}}}\cdot T\bar{\lambda} \cdot w).
\end{align*}
Because the nonholonomic reduced distributional two-form
$\omega^B_{\bar{\mathcal{K}}}$ is non-degenerate, it follows that the equation
$T\bar{\gamma}\cdot \tilde{X}^\varepsilon = X^B_{\bar{\mathcal{K}}}\cdot
\bar{\varepsilon},$ is equivalent to the equation $T\bar{\lambda}\cdot \tilde{X} \cdot \varepsilon
= \tau_{\bar{\mathcal{K}}}\cdot T\bar{\varepsilon} \cdot X^B_{h_{\bar{\mathcal{K}}}
\cdot \bar{\varepsilon}}. $
Thus, $\varepsilon$ and $\bar{\varepsilon}$ satisfy the equation
$T\bar{\lambda}\cdot \tilde{X} \cdot \varepsilon
= \tau_{\bar{\mathcal{K}}}\cdot T\bar{\varepsilon} \cdot X^B_{h_{\bar{\mathcal{K}}}
\cdot \bar{\varepsilon}},$ if and only if they satisfy
the Type II of Hamilton-Jacobi equation
$T\bar{\gamma}\cdot \tilde{X}^\varepsilon
= X^B_{\bar{\mathcal{K}}}\cdot \bar{\varepsilon}.$
\hskip 0.3cm $\blacksquare$\\

For a given nonholonomic reducible CMH system
$(T^*Q,G,\omega^B,\mathcal{D},H,F,u)$ with the associated
the distributional CMH system $(\mathcal{K},\omega^B_{\mathcal{K}},
H_{\mathcal{K}}, F^B_{\mathcal{K}}, u^B_{\mathcal{K}})$
and the nonholonomic reduced distributional CMH system
$(\bar{\mathcal{K}},\omega^B_{\bar{\mathcal{K}}},h_{\bar{\mathcal{K}}},
f^B_{\bar{\mathcal{K}}}, u^B_{\bar{\mathcal{K}}})$,
we know that the nonholonomic dynamical vector field
$X^B_{\mathcal{K}}$ and the nonholonomic reduced dynamical vector field
$X^B_{\bar{\mathcal{K}}}$ are $\pi_{/G}$-related, that is,
$X^B_{\bar{\mathcal{K}}}\cdot \pi_{/G}=T\pi_{/G}\cdot X^B_{\mathcal{K}}.$
Then we can prove the following Theorem 7.4,
which states the relationship between the solutions of Type II of
Hamilton-Jacobi equations and nonholonomic reduction.

\begin{theo}
For a given nonholonomic reducible CMH system
$(T^*Q,G,\omega^B,\mathcal{D},H,F,u)$ with the associated
the distributional CMH system $(\mathcal{K},\omega^B_{\mathcal{K}},
H_{\mathcal{K}}, F^B_{\mathcal{K}}, u^B_{\mathcal{K}})$
and the nonholonomic reduced distributional CMH system
$(\bar{\mathcal{K}},\omega^B_{\bar{\mathcal{K}}},h_{\bar{\mathcal{K}}},
f^B_{\bar{\mathcal{K}}}, u^B_{\bar{\mathcal{K}}})$, assume that
$\gamma: Q \rightarrow T^*Q$ is an one-form on $Q$,
and $\varepsilon: T^* Q \rightarrow T^* Q $ is a $G$-invariant
symplectic map with respect to $\omega^B$,
and $\bar{\gamma}=\pi_{/G}(\gamma): Q \rightarrow T^* Q/G $, and
$\bar{\varepsilon}=\pi_{/G}(\varepsilon): T^* Q \rightarrow T^* Q/G. $
Under the hypotheses and notations of Theorem 7.3, then we have that
$\varepsilon$ is a solution of the Type II of Hamilton-Jacobi equation, $T\gamma\cdot
\tilde{X}^\varepsilon= X^B_{\mathcal{K}}\cdot \varepsilon, $ for the distributional
CMH system $(\mathcal{K},\omega^B_{\mathcal{K}},H_{\mathcal{K}},
F^B_{\mathcal{K}}, u^B_{\mathcal{K}})$, if and only if
$\varepsilon$ and $\bar{\varepsilon}$ satisfy the Type II of
Hamilton-Jacobi equation $T\bar{\gamma}\cdot \tilde{X}^\varepsilon =
X^B_{\bar{\mathcal{K}}}\cdot \bar{\varepsilon}, $ for the nonholonomic reduced
distributional CMH system $ (\bar{\mathcal{K}},
\omega^B_{\bar{\mathcal{K}}}, h_{\bar{\mathcal{K}}},
f^B_{\bar{\mathcal{K}}}, u^B_{\bar{\mathcal{K}}} ). $
\end{theo}

\noindent{\bf Proof: }At first, for the dynamical vector field of the distributional CMH system
$(\mathcal{K},\omega^B_{\mathcal {K}}, H_{\mathcal{K}},
F^B_{\mathcal{K}}, u^B_{\mathcal{K}})$,
$\tilde{X}=X^B_{(\mathcal{K},\omega^B_{\mathcal{K}},
H_{\mathcal{K}}, F^B_{\mathcal{K}}, u^B_{\mathcal{K}})}
=X^B_\mathcal {K}+ F^B_{\mathcal{K}}+u^B_{\mathcal{K}}$,
and $F^B_{\mathcal{K}}=\tau_{\mathcal{K}}\cdot \textnormal{vlift}(F_{\mathcal{M}})X^B_H$,
and $u^B_{\mathcal{K}}=\tau_{\mathcal{K}}\cdot \textnormal{vlift}(u_{\mathcal{M}})X^B_H$,
note that $T\pi_{Q}\cdot \textnormal{vlift}(F_{\mathcal{M}})X^B_H=T\pi_{Q}\cdot \textnormal{vlift}(u_{\mathcal{M}})X^B_H=0, $
then we have that $T\pi_{Q}\cdot F^B_{\mathcal{K}}=T\pi_{Q}\cdot u^B_{\mathcal{K}}=0,$
and hence $T\pi_{Q}\cdot \tilde{X}\cdot \varepsilon
=T\pi_{Q}\cdot X^B_{\mathcal{K}}\cdot \varepsilon. $
Next, under the hypotheses and notations of Theorem 7.3,
$\textmd{Im}(\gamma)\subset \mathcal{M},$ and
it is $G$-invariant, $\textmd{Im}(T\gamma)\subset \mathcal{K}, $
and hence $\textmd{Im}(T\bar{\gamma})\subset \bar{\mathcal{K}}, $ in
this case, $\pi^*_{/G}\cdot\omega^B_{\bar{\mathcal{K}}}\cdot \tau_{\bar{\mathcal{K}}}= \tau_{\mathcal{U}}\cdot
\omega^B_{\mathcal{M}}= \tau_{\mathcal{U}}\cdot i_{\mathcal{M}}^*\cdot
\omega^B, $ along $\textmd{Im}(T\bar{\gamma})$.
In addition, from the distributional magnetic Hamiltonian equation (7.1),
we have that $X^B_{\mathcal{K}}=\tau_{\mathcal{K}}\cdot X^B_H, $
and from the nonholonomic reduced distributional magnetic Hamiltonian equation (7.2),
we have that $X^B_{\bar{\mathcal{K}}}
=\tau_{\bar{\mathcal{K}}}\cdot X^B_{h_{\bar{\mathcal{K}}}},$
and the nonholonomic dynamical vector field $X^B_{\mathcal{K}}$ and the
nonholonomic reduced dynamical vector field $X^B_{\bar{\mathcal{K}}}$ are $\pi_{/G}$-related,
that is, $X^B_{\bar{\mathcal{K}}}\cdot \pi_{/G}=T\pi_{/G}\cdot
X^B_{\mathcal{K}}. $ Note that $\varepsilon(\mathcal{M})\subset \mathcal{M},$
and hence $X^B_H\cdot \varepsilon \in \mathcal{F}$ along $\varepsilon$,
and $\tau_{\mathcal{K}}\cdot X^B_{H}\cdot \varepsilon
= X^B_{\mathcal{K}}\cdot \varepsilon \in \mathcal{K}$.
Then $\tau_{\bar{\mathcal{K}}}\cdot T\pi_{/G}(X^B_{\mathcal{K}}\cdot \varepsilon)
=\tau_{\bar{\mathcal{K}}}\cdot (T\pi_{/G}(X^B_{\mathcal{K}}))\cdot (\varepsilon)
= \tau_{\bar{\mathcal{K}}}\cdot (X^B_{\bar{\mathcal{K}}}\cdot \pi_{/G})\cdot (\varepsilon)
= \tau_{\bar{\mathcal{K}}}\cdot X^B_{\bar{\mathcal{K}}}\cdot \pi_{/G}(\varepsilon)
= X^B_{\bar{\mathcal{K}}}\cdot \bar{\varepsilon}.$
Thus, using the non-degenerate and nonholonomic reduced distributional two-form
$\omega^B_{\bar{\mathcal{K}}}$, note that $\tau_{\bar{\mathcal{K}}}\cdot T\bar{\gamma} =T\bar{\gamma},$
for any $w \in \mathcal{F}, \; \tau_{\mathcal{K}} \cdot w\neq 0, $ and
$\tau_{\bar{\mathcal{K}}}\cdot T\pi_{/G}\cdot w \neq 0, $ then
we have that
\begin{align*}
& \omega^B_{\bar{\mathcal{K}}}(T\bar{\gamma} \cdot \tilde{X}^\varepsilon
- X^B_{\bar{\mathcal{K}}}\cdot \bar{\varepsilon}, \; \tau_{\bar{\mathcal{K}}}\cdot T\pi_{/G}\cdot w)\\
& = \omega^B_{\bar{\mathcal{K}}}(T\bar{\gamma} \cdot \tilde{X}^\varepsilon, \;
\tau_{\bar{\mathcal{K}}}\cdot T\pi_{/G}\cdot w)-
\omega^B_{\bar{\mathcal{K}}}(X^B_{\bar{\mathcal{K}}}\cdot \bar{\varepsilon},
\; \tau_{\bar{\mathcal{K}}}\cdot T\pi_{/G} \cdot w) \\
& = \omega^B_{\bar{\mathcal{K}}}(\tau_{\bar{\mathcal{K}}}\cdot T\bar{\gamma}\cdot \tilde{X}^
\varepsilon, \; \tau_{\bar{\mathcal{K}}}\cdot T\pi_{/G}\cdot w)
-\omega^B_{\bar{\mathcal{K}}}(\tau_{\bar{\mathcal{K}}}\cdot X^B_{\bar{\mathcal{K}}}
\cdot \pi_{/G}\cdot \varepsilon,
\; \tau_{\bar{\mathcal{K}}}\cdot T\pi_{/G}\cdot w)\\
& = \omega^B_{\bar{\mathcal{K}}}\cdot \tau_{\bar{\mathcal{K}}}(T\pi_{/G}\cdot T\gamma \cdot \tilde{X}^
\varepsilon, \; T\pi_{/G} \cdot w)
- \omega^B_{\bar{\mathcal{K}}}\cdot \tau_{\bar{\mathcal{K}}}(T\pi_{/G}\cdot
X^B_{\mathcal{K}}\cdot \varepsilon, \; T\pi_{/G}\cdot w)\\
& = \pi^*_{/G}\cdot\omega^B_{\bar{\mathcal{K}}}\cdot \tau_{\bar{\mathcal{K}}}(T\gamma \cdot \tilde{X}^
\varepsilon, \; w)
- \pi^*_{/G}\cdot\omega^B_{\bar{\mathcal{K}}}\cdot \tau_{\bar{\mathcal{K}}}(X^B_{\mathcal{K}} \cdot
\varepsilon, \; w)\\
& = \tau_{\mathcal{U}}\cdot i_{\mathcal{M}}^* \cdot
\omega^B (T\gamma \cdot \tilde{X}^\varepsilon, \; w)
- \tau_{\mathcal{U}}\cdot i_{\mathcal{M}}^* \cdot
\omega^B (X^B_{\mathcal{K}} \cdot
\varepsilon, \; w).
\end{align*}
In the case we note that $\tau_{\mathcal{U}}\cdot i_{\mathcal{M}}^* \cdot
\omega^B=\tau_{\mathcal{K}}\cdot i_{\mathcal{M}}^* \cdot
\omega^B= \omega^B_{\mathcal{K}}\cdot \tau_{\mathcal{K}}, $
and
$\tau_{\mathcal{K}}\cdot T\gamma =T\gamma, \; \tau_{\mathcal{K}} \cdot X^B_{\mathcal{K}}
= X^B_{\mathcal{K}}$, since $\textmd{Im}(\gamma)\subset
\mathcal{M}, $ and $\textmd{Im}(T\gamma)\subset \mathcal{K}. $
Thus, we have that
\begin{align*}
& \omega^B_{\bar{\mathcal{K}}}(T\bar{\gamma} \cdot \tilde{X}^\varepsilon
- X^B_{\bar{\mathcal{K}}}\cdot \bar{\varepsilon}, \; \tau_{\bar{\mathcal{K}}}\cdot T\pi_{/G}\cdot w)\\
& = \omega^B_{\mathcal{K}}\cdot \tau_{\mathcal{K}}(T\gamma \cdot \tilde{X}^
\varepsilon, \; w)- \omega^B_{\mathcal{K}}\cdot \tau_{\mathcal{K}}(X^B_{\mathcal{K}} \cdot \varepsilon, \; w)\\
& = \omega^B_{\mathcal{K}}(\tau_{\mathcal{K}} \cdot T\gamma \cdot \tilde{X}^
\varepsilon, \; \tau_{\mathcal{K}} \cdot w)
- \omega^B_{\mathcal{K}}(\tau_{\mathcal{K}} \cdot X^B_{\mathcal{K}}\cdot \varepsilon,
\; \tau_{\mathcal{K}} \cdot w)\\
& = \omega^B_{\mathcal{K}}(T\gamma \cdot \tilde{X}^
\varepsilon- X^B_{\mathcal{K}}\cdot \varepsilon, \; \tau_{\mathcal{K}} \cdot w).
\end{align*}
Because the distributional two-form $\omega^B_{\mathcal{K}}$
and the nonholonomic reduced distributional
two-form $\omega^B_{\bar{\mathcal{K}}}$ are both non-degenerate,
it follows that the equation
$T\bar{\gamma}\cdot \tilde{X}^\varepsilon=
X^B_{\bar{\mathcal{K}}}\cdot \bar{\varepsilon}, $
is equivalent to the equation $T\gamma\cdot \tilde{X}^\varepsilon= X^B_{\mathcal{K}}\cdot \varepsilon. $
Thus, $\varepsilon$ is a solution of the Type II of Hamilton-Jacobi equation
$T\gamma\cdot \tilde{X}^\varepsilon= X^B_{\mathcal{K}}\cdot \varepsilon, $ for the distributional
CMH system $(\mathcal{K},\omega^B_{\mathcal {K}},H_{\mathcal{K}},
F^B_{\mathcal{K}}, u^B_{\mathcal{K}})$, if and only if
$\varepsilon$ and $\bar{\varepsilon} $ satisfy the Type II of Hamilton-Jacobi
equation $T\bar{\gamma}\cdot \tilde{X}^\varepsilon=
X^B_{\bar{\mathcal{K}}}\cdot \bar{\varepsilon}, $ for the
nonholonomic reduced distributional CMH system
$(\bar{\mathcal{K}},\omega^B_{\bar{\mathcal{K}}},h_{\bar{\mathcal{K}}}, f^B_{\bar{\mathcal{K}}}, u^B_{\bar{\mathcal{K}}})$.
\hskip 0.3cm $\blacksquare$ \\

\begin{rema}
It is worthy of noting that,
the Type I of Hamilton-Jacobi equation
$T\bar{\gamma}\cdot \tilde{X}^ \gamma = X^B_{\bar{\mathcal{K}}}\cdot
\bar{\gamma} ,$ is the equation of
the reduced differential one-form $\bar{\gamma}$; and
the Type II of Hamilton-Jacobi equation $T\bar{\gamma}\cdot \tilde{X}^\varepsilon
= X^B_{\bar{\mathcal{K}}}\cdot \bar{\varepsilon},$ is the equation of the symplectic
diffeomorphism map $\varepsilon$ and the reduced symplectic
diffeomorphism map $\bar{\varepsilon}$.
If the nonholonomic CMH system with symmetry we considered
has not any the external force and control, that is, $F=0 $ and $W=\emptyset$,
in this case, the nonholonomic CMH system with symmetry
$(T^*Q,G,\omega^B,\mathcal{D},H,F,W)$ is just the nonholonomic magnetic Hamiltonian system
with symmetry $(T^*Q,G,\omega^B,\mathcal{D},H)$,
and with the magnetic symplectic form $\omega^B$ on $T^*Q$.
From the above Type I and Type II of Hamilton-Jacobi theorems, that is,
Theorem 7.2 and Theorem 7.3, we can get the Theorem 5.2 and Theorem 5.3 in Wang \cite{wa21c}.
It shows that Theorem 7.2 and Theorem 7.3 can be regarded as an extension of two types of
Hamilton-Jacobi theorem for the nonholonomic magnetic Hamiltonian system with symmetry given in
Wang \cite{wa21c} to that for the system with the external force and control.
In particular, in this case, if $B=0$, then
the magnetic symplectic form $\omega^B$
is just the canonical symplectic form $\omega$ on $T^*Q$,
from the proofs of Theorem 7.2 and Theorem 7.3, we can also get the Theorem 4.2
and Theorem 4.3 in Le\'{o}n and Wang \cite{lewa15}.
It shows that Theorem 7.2 and Theorem 7.3 can be regarded as an extension of two types of
Hamilton-Jacobi theorem for the nonholonomic Hamiltonian system with symmetry given in
Le\'{o}n and Wang \cite{lewa15} to that for the system with the magnetic,
external force and control.
\end{rema}

The theory of controlled mechanical system is a very important subject,
its research gathers together some separate areas of research such as
mechanics, differential geometry and nonlinear control theory, etc.,
and the emphasis of this research on geometry is motivated by the
aim of understanding the structure of equations of motion of the
system in a way that helps both analysis and design. Thus, it is
natural to study the controlled mechanical systems by combining with the
analysis of dynamical systems and the geometric reduction theory of
Hamiltonian and Lagrangian systems. Following
the theoretical development of geometric mechanics, a lot
of important problems about this subject are being explored and
studied, see Le\'{o}n and Wang \cite{lewa15},
Marsden et al. \cite{mawazh10}, Ratiu and Wang \cite{rawa12},
Wang \cite{wa18, wa15a, wa17, wa13d, wa21a, wa21c, wa22a, wa20a, wa13e},
and Wang and Zhang \cite{wazh12}.
These research works from the geometrical point of view
reveal the internal relationships of
the geometrical structures of phase spaces, symmetric reductions,
constraints, dynamical vector fields and controls of a
mechanical system and its regular reduced systems.
In particular, it is the key thought of the researches of geometrical mechanics
of the professor Jerrold E. Marsden to explore and reveal the deeply internal
relationship between the geometrical structure of phase space and the dynamical
vector field of a mechanical system. It is also our goal of pursuing and inheriting.

\end{document}